\documentclass[10pt]{amsart}
\usepackage{amsmath}
\usepackage{amssymb}
\usepackage{amsfonts}
\usepackage{amsthm}
\usepackage{enumerate}
\usepackage[all]{xy}
\usepackage[latin1]{inputenc}
\usepackage{graphicx}
\usepackage{latexsym}
\usepackage{accents}
\usepackage{mathrsfs}
\usepackage{indentfirst}
\usepackage{fancyhdr}
\usepackage{color}
\usepackage{graphicx}
\usepackage{epstopdf}
\input xy

\usepackage{tikz}

\makeatletter
    \newtheoremstyle{definition}% name of the style to be used
        {5pt}% measure of space to leave above the theorem. E.g.: 3pt
        {3pt}% measure of space to leave below the theorem. E.g.: 3pt
        {}% name of font to use in the body of the theorem
        {0pt}% measure of space to indent
        {\scshape}% name of head font
        {.}% punctuation between head and body
        {5pt}% space after theorem head; " " = normal interword space
        {\thmname{#1} \thmnumber{#2} \thmnote{[#3]}} % Manually specify head

\newtheoremstyle{theorems}% name of the style to be used
        {5pt}% measure of space to leave above the theorem. E.g.: 3pt
        {3pt}% measure of space to leave below the theorem. E.g.: 3pt
        {\itshape}% name of font to use in the body of the theorem
        {0pt}% measure of space to indent
        {\scshape}% name of head font
        {.}% punctuation between head and body
        {5pt}% space after theorem head; " " = normal interword space
        {\thmname{#1} \thmnumber{#2}\thmnote{[#3]}} % Manually specify head

\swapnumbers

\theoremstyle{theorems}

\newtheorem{Theorem}{Theorem}[section]
\newtheorem{Lemma}[Theorem]{Lemma}
\newtheorem{Definition-Lemma}[Theorem]{Definition-Lemma}

\newtheorem{Definition}[Theorem]{Definition}

\newtheorem{Proposition}[Theorem]{Proposition}
\newtheorem{Example}[Theorem]{Example}
\newtheorem{Remark}[Theorem]{Remark}
\newtheorem{Conjecture}[Theorem]{Conjecture}

\begin{document}

\title [Positivity for quantum cluster algebras from unpunctured orbifolds]
 {Positivity for quantum cluster algebras from unpunctured orbifolds}

\author{Min Huang $\;\;\;\;\;\;$}
\address{Min Huang
\newline D\'{e}partement de math\'{e}matiques, Universit\'{e} de Sherbrooke, Sherbrooke, Qu\'{e}bec,
J1K 2R1, Canada}
\email{minhuang1989@hotmail.com}

\date{version of \today}



\thanks{}

\maketitle

\vspace{-25pt}

\begin{abstract}

We give the quantum Laurent expansion formula for the quantum cluster algebras from unpunctured orbifolds with arbitrary coefficients and quantization. As an application, positivity for such class of quantum cluster algebras is given. For technical reasons, it will always be assumed that the weights of the orbifold points are $2$.
\end{abstract}

\smallskip

\tableofcontents

\section{Introduction}

Cluster algebras are commutative algebras introduced by Fomin and Zelevinsky around the year 2000. The quantum cluster algebras were later introduced in \cite{BZ}. The theory of cluster algebras is related to numerous other fields including Lie theory, representation theory of algebras, the periodicity issue, Teichm$\ddot{u}$ller theory and mathematical physics.

\medskip

A cluster algebra is a subalgebra of rational function field with a distinguished set of generators, called cluster variables. Different cluster variables are related by an iterated procedure, called mutation. By construction, cluster variables are rational functions. In \cite{fz1}, Fomin and Zelevinsky in fact proved that they are Laurent polynomials of initial cluster variables, which is known as Laurent phenomenon. It was proved that these Laurent polynomials have non-negative coefficients, known as positivity, see \cite{LS,GHKK}.

\medskip

Quantum cluster algebras are quantum deformations of cluster algebras. It was proved in \cite{BZ} that the Laurent phenomenon have a quantum version in the quantum setting. More precisely, the quantum cluster variables are quantum Laurent polynomials of initial quantum cluster variables. It was conjectured that the coefficients are in $\mathbb N[q^{\pm 1/2}]$, known as positivity conjecture for quantum cluster algebras, where $q$ is the quantum parameter. Kimura and Qin proved the positivity conjecture for the acyclic skew-symmetric quantum cluster algebras in \cite{KQ}, Davison proved this conjecture for skew-symmetric case  in \cite{D}. However, the positivity for skew-symmetrizable quantum cluster algebras has not been considered.

\medskip

The original motivation of Fomin and Zelevinsky is to provide a combinatorial characterization to the canonical bases in quantum groups (see \cite{L,K}) and the total positivity in algebraic groups. They conjectured that the cluster structure should serve as an algebraic framework for the study of the ``dual canonical bases" in various coordinate rings and their $q$-deformations. Particularly, they conjectured all cluster monomials belong to the dual canonical bases, this was proved recently in \cite{KKKO,Q}. Generally it can be very hard to write the dual canonical bases explicitly. From this point of view, it is important to give the explicit (quantum) Laurent polynomial of a (quantum) cluster variable with respect to any (quantum) cluster.

\medskip

(Quantum) cluster algebras from orbifolds is an important class of skew-symmetrizable (quantum) cluster algebras. Almost all (quantum) cluster algebras of finite mutation type are in this class, see \cite{FST3}. Using the unfolding method, the positivity for cluster algebras from orbifolds can be deduced by the positivity for the cluster algebras from surfaces, see \cite{FST3}. However, it should be emphasized that the unfolding method does not give the information of the $q$-coefficients. Thus, we do not know the quantum Laurent expansion formula for quantum cluster algebras from orbifolds even that is known for quantum cluster algebras form surfaces. The aim of this study is to solve the positivity conjecture for such class quantum cluster algebras by giving the quantum Laurent polynomial of a (quantum) cluster variable with respect to any (quantum) cluster. We would generalize the methods in \cite{H,H1} from the surface case to the orbifold case.

\medskip

The paper is organized as follows. We first give background on cluster algebras and cluster algebras from unpunctured orbifolds in Section \ref{pre}, then the background of the quantum version is given in Section \ref{pre2}. The main results, Theorems \ref{partition bi}, \ref{expansion-comm} for cluster algebras, and Theorems \ref{mainthm}, \ref{expansion} for quantum cluster algebras are stated in Section \ref{cle} and Section \ref{qle}, respectively. In Section \ref{compare}, the necessary preparations are made for the proof of Theorems \ref{partition bi}, \ref{mainthm}. We finally prove Theorem \ref{partition bi} and Theorem \ref{mainthm} in Section \ref{main1} and Section \ref{main2}, respectively.

\medskip

Throughout this study, we denote by $E(G)$ the edges set of a graph $G$ and by $|S|$ the cardinality of a set $S$.

\section{Preliminaries on cluster algebras}\label{pre}

\subsection{Commutative cluster algebras}

Herein, we recall the definition of cluster algebra in \cite{fz1}.

\medskip

A triple $(\mathbb P,\oplus,\cdot)$ is called a \emph{semifield} if $(\mathbb P,\cdot)$ is an abelian multiplicative group and $(\mathbb P,\oplus)$ is a commutative semigroup such that $``\oplus"$ is distributive with respect to $``\cdot"$. The \emph{tropical semifield} ${\rm Trop}(u_1,\cdots,u_l)$ is a semifield freely generated by $u_1,\cdots,u_l$ as abelian groups with $\oplus$ defined by $\prod_ju_j^{a_j}\oplus \prod_ju_j^{b_j}=\prod_ju_j^{min(a_j,b_j)}$. Let $(\mathbb P,\oplus,\cdot)$ be a semifield. The group ring $\mathbb {ZP}$ will be used as \emph{ground ring}. For a given integer $n$, let $\mathcal F$ be the rational function field in $n$ independent variables, with coefficients in $\mathbb {QP}$.

\medskip

A seed $t$ in $\mathcal F$ is a triple $(x(t),y(t),B(t))$, where

\begin{enumerate}[$(1)$]

  \item $x(t)=\{x_1(t),\cdots,x_n(t)\}$ such that $\mathcal F$ is freely generated by $x(t)$ over $\mathbb {QP}$.

  \item $y(t)=\{y_1(t),\cdots,y_n(t)\}\subseteq \mathbb P$.

  \item $B(t)=(b_{ij})$ is an $n\times n$ skew-symmetrizable integer matrix.

\end{enumerate}

\medskip

Given a seed $t$ in $\mathcal F$, for any $k\in [1,n]$, we define the \emph{mutation} of $t$ at the $k$-th direction to be the new seed $t'=\mu_k(t)=(x(t'),y(t'),B(t'))$, where

\begin{enumerate}[$(1)$]

  \item \[\begin{array}{ccl} x_i(t') &=&
         \left\{\begin{array}{ll}
              x_i(t), &\mbox{if $i\neq k$},  \\
             \frac{y_k(t)\prod x_i(t)^{[b_{ik}]_{+}}+\prod x_i(t)^{[-b_{ik}]_{+}}}{(y_k(t)\oplus 1)x_k(t)}, &\mbox{otherwise}.
         \end{array}\right.
        \end{array}\]

  \item \[\begin{array}{ccl} y_i(t') &=&
         \left\{\begin{array}{ll}
              y^{-1}_k(t), &\mbox{if $i=k$},  \\
              y_j(t)y_k(t)^{[b_{kj}]_{+}}(1\oplus y_k(t))^{-b_{kj}}, &\mbox{otherwise}.
         \end{array}\right.
        \end{array}\]

  \item $B(t')=(b'_{ij})$ is determined by $B(t)=(b_{ij})$:

  \[\begin{array}{ccl} b'_{ij} &=&
         \left\{\begin{array}{ll}
              -b_{ij}, &\mbox{if $i=k$ or $j=k$},  \\
              b_{ij}+[b_{ik}]_{+}[b_{kj}]_{+}-[-b_{ik}]_{+}[-b_{kj}]_{+}, &\mbox{otherwise}.
         \end{array}\right.
        \end{array}\]

\end{enumerate}
here $[a]_{+}=max(a,0)$ for any $a\in\mathbb Z$.

\medskip

We define a \emph{cluster algebra} $\mathcal A$ (of rank $n$) over $\mathbb P$ as following:
\begin{enumerate}[$(1)$]
  \item Choose an initial seed $t_0=(x(t_0),y(t_0),B(t_0))$.
  \item All the seeds $t$ are obtained from $t_0$ by iterated mutations at directions $k\in[1,n]$.
  \item $\mathcal A=\mathbb{ZP}[x_i(t)]_{i\in [1,n],t}$.
  \item $x(t)$ is called a \emph{cluster} of $\mathcal A$ for any $t$.
  \item $x_i(t)$ is called a \emph{cluster variable} of $\mathcal A$ for any $i\in [1,n]$ and $t$.
  \item A monomial in $x(t)$ is called a \emph{cluster monomial} of $\mathcal A$ for any $t$.
  \item $y(t)$ is called a \emph{coefficient tuple} of $\mathcal A$ for any $t$.
  \item $B(t)$ is called an \emph{exchange matrix} of $\mathcal A$ for any $t$.
\end{enumerate}
In particular, $\mathcal A$ is called of \emph{geometric type} if $\mathbb P$ is a tropical semifield.

\medskip

When $\mathbb P={\rm Trop}(u_1,\cdots,u_l)$, let $m=n+l$. For a seed $t$ in $\mathcal F$, we have $y_j(t)=\prod u_i^{a_{ij}}$ for some integers $a_{ij}$. We can write $t$ as $(\widetilde x(t),\widetilde B(t))$, where

\begin{enumerate}[$(1)$]

  \item $\widetilde x(t)=\{x_1(t),\cdots,x_n(t),x_{n+1}(t)=u_1,\cdots,x_m(t)=u_l\}$.

  \item $\widetilde B(t)=(b_{ij})$ is an $m\times n$ with $b_{ij}=a_{i-n,j}$ for $i\in [n+1,m]$.

\end{enumerate}
In this case, the mutation of $t$ at direction $k$ is $t'=\mu_k(t)=(\widetilde x(t'),\widetilde B(t'))$, where

\begin{enumerate}[$(1)$]

  \item \[\begin{array}{ccl} x_i(t') &=&
         \left\{\begin{array}{ll}
              x_i(t), &\mbox{if $i\neq k$},  \\
             \frac{\prod x_i(t)^{[b_{ik}]_{+}}+\prod x_i(t)^{[-b_{ik}]_{+}}}{x_k(t)}, &\mbox{otherwise}.
         \end{array}\right.
        \end{array}\]

  \item $\widetilde B(t')=(b'_{ij})$ is determined by $\widetilde B(t)=(b_{ij})$:

  \[\begin{array}{ccl} b'_{ij} &=&
         \left\{\begin{array}{ll}
              -b_{ij}, &\mbox{if $i=k$ or $j=k$},  \\
              b_{ij}+[b_{ik}]_{+}[b_{kj}]_{+}-[-b_{ik}]_{+}[-b_{kj}]_{+}, &\mbox{otherwise}.
         \end{array}\right.
        \end{array}\]
We also denote $\widetilde B(t')=(b'_{ij})$ by $\mu_k\widetilde B(t)$.
\end{enumerate}

\medskip

\subsection{Cluster algebras from orbifolds} Herein, we recall some combinatorial notation related to the orbifolds in \cite{FST3,FST4}. Let $\mathcal O$ be a connected Riemann surface with boundary. Fix a non-empty set $M$ of marked points in the closure of $\mathcal O$ with at least one marked point on each boundary component. Fix a finite set $U$ in the interior of $\mathcal O$ such that $U\cap M=\emptyset$. We call the triple $(\mathcal O, M, U)$ an \emph{orbifold}. Marked points in the interior of $\mathcal O$ are called \emph{punctures}. The points in $U$ are called \emph{orbifold points.} Each orbifold point of $\mathcal O$ comes with a weight $w=1/2$ or $2$. For technical reasons, we always assume that the weights of orbifold pints are $2$.

\medskip

In this paper, we consider the case that $(\mathcal O,M,U)$ without punctures, and we refer $(\mathcal O,M,U)$, or simply $\mathcal O$ if there is no confusion, as an \emph{unpunctured orbifold}.

\medskip

An \emph{arc} $\gamma$ in $(\mathcal O,M,U)$ is a curve (up to isotopy of $\mathcal O\setminus (M\cup U)$) in $\mathcal O$ such that: the endpoints are in $M$ or one endpoint in $M$ the another in $U$; $\gamma$ does not cross itself, except its endpoints may coincide; except for the endpoints, $\gamma$ is disjoint from $M\cup U$ and from the boundary of $\mathcal O$; and if $\gamma$ cuts out a monogon then this monogon contains at least two points of $U$. If the endpoints of $\gamma$ are in $M$, then $\gamma$ is called an \emph{ordinary arc}; if one endpoint of $\gamma$ is in $M$ and another is in $U$, then $\gamma$ is called a \emph{pending arc}. A pending arc incident to an orbifold point with weight $w$ is assigned with the same weight $w$. An ordinary arc is assigned with the weight $w=1$.

\medskip

For two arcs $\gamma, \gamma'$ in $(\mathcal O,M,U)$, if $\gamma$ is an ordinary arc, then the \emph{crossing number} $N(\gamma,\gamma')$ of $\gamma$ with $\gamma'$ is the minimum of the numbers of crossings of arcs $\alpha$ and $\alpha'$, where $\alpha$ is isotopic to $\gamma$ and $\alpha'$ is isotopic to $\gamma'$; if $\gamma$ is a pending arc which incident an orbifold point $o$, denote by $l(\gamma)$ the loop cutting out a monogon which contains only $o$, then the \emph{crossing number} $N(\gamma,\gamma')$ of $\gamma$ with $\gamma'$ is defined to be $N(l(\gamma),\gamma')$. It should be noted that $N(\gamma,\gamma')\neq N(\gamma',\gamma)$ generally.

\medskip

We call two arcs $\gamma$ and $\gamma'$ are \emph{compatible} if $N(\gamma,\gamma')=0$. A \emph{triangulation} is a maximal collection of compatible arcs. Given a triangulation $T$ and a non-boundary arc $\tau$ in $T$, there exists a unique arc $\tau'$ such that $(T\setminus \{\tau\})\cup \{\tau'\}$ is a new triangulation. We denote $(T\setminus \{\tau\})\cup \{\tau'\}$ by $\mu_{\tau}(T)$. For an arc $\gamma$ and a triangulation $T$, the \emph{crossing number} $N(\gamma,T)$ of $\gamma$ and $T$ is defined as $\sum_{\tau\in T}N(\gamma,\tau)$. We call a triangulation $T$ an \emph{indexed triangulation} if the order of the arcs in $T$ is fixed. For $\tau\in T$, denote by $w(\tau)$ the weight of $\tau$.

\medskip

\centerline{\includegraphics{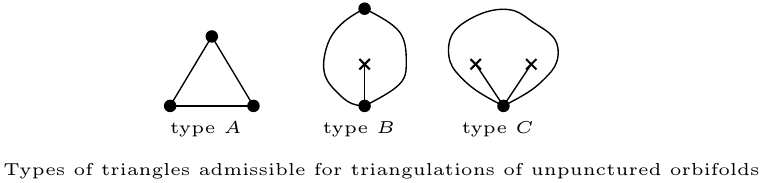}}

For two non-boundary arcs $\tau$, $\tau'$ in an indexed triangulation $T=\{\tau_1,\cdots,\tau_n,\cdots,\tau_l\}$ and a triangle $\Delta$ of $T$, we may assume that $\tau_1,\cdots,\tau_n$ are the non-boundary arcs, define
\[\begin{array}{ccl} b^{T,\Delta}_{\tau\tau'}=
         \left\{\begin{array}{ll}
              w(\tau'), &\mbox{if $\tau,\tau'$ are sides of $\Delta$ and $\tau'$ follows $\tau$ in the clockwise order},  \\
              -w(\tau'), &\mbox{if $\tau,\tau'$ are sides of $\Delta$ and $\tau$ follows $\tau'$ in the clockwise order}, \\
              0, &\mbox{otherwise.}
         \end{array}\right.
 \end{array}\]
and $b^T_{\tau\tau'}=\sum_{\Delta}b^{T,\Delta}_{\tau\tau'}$.
The $n\times n$ matrix $B^T=(b^T_{ij})$ with $b^T_{ij}=b^T_{\tau_i\tau_j}, 1\leq i,j\leq n$ is called the \emph{signed adjacency matrix} of $T$, see \cite{FST3,FST4}. $B^T$ is skew-symmetrizable. In fact, let $D^T=diag(w(\tau))_{\tau\in T}$, we have $D^TB^T$ is skew-symmetric.

\medskip

For any non-boundary arc $\tau\in T$, we have $B^{\mu_{\tau}(T)}=\mu_{\tau}(B^T)$.

\medskip

A cluster algebra $\mathcal A$ is called \emph{coming from $(\mathcal O, M, U)$} if there exists a triangulation $T$ such that $B^{T}$ is an exchange matrix of $\mathcal A$.

\medskip

Throughout this paper, let $\mathcal O$ be an unpunctured orbifold and $T=\{\tau_1,\cdots,\tau_n,\cdots,\tau_l\}$ be an indexed triangulation, and let $\gamma$ be an oriented arc in $\mathcal O$. When an arc $\tau\in T$ is fixed, let $T'=\mu_{\tau}(T)$ and $\tau'$ be the arc obtained from $T$ by flip at $\tau$. We always assume that $\tau_1,\cdots,\tau_n$ are the non-boundary arcs in $T$.

\medskip

\subsection{Snake graphs and Perfect matchings}

In this subsection, fix a triangulation $T$, we construct a snake graph $G_{T,\gamma}$ for each arc $\gamma$. For the surface case, see \cite[Section 4]{MSW}, \cite{CS,H,H1}.

\medskip

If $\gamma$ is a pending arc which incident an orbifold point $o$, denote by $l(\gamma)$ the loop cutting out a monogon which contains only $o$. Let $\zeta=\gamma$ if $\gamma$ is an ordinary arc and $\zeta=l(\gamma)$ if $\gamma$ is a pending arc. Let $p_0$ be the starting point of $\zeta$, and let $p_{d+1}$ be its endpoint. We assume that $\zeta$ crosses $T$ at $p_1,\cdots,p_d$ in order. Let $\tau_{i_j}$ be the arc in $T$ containing $p_j$. If $\tau_{i_j}$ is an ordinary arc, let $\Delta_{j-1}$ and $\Delta_{j}$ be the two ideal triangles in $T$ on either side of $\tau_{i_j}$. If $\tau_{i_j}$ is a pending arc, let $\Delta_{j-1}=\Delta_j$ be the unique triangle with an edge $\tau_{i_j}$.

\medskip

We associate each $p_j$ with a \emph{tile} $G(p_j)$ as follows. Define $\Delta_1^j$ and $\Delta_2^j$ to be two triangles with edges labeled as in $\Delta_{j-1}$ and $\Delta_{j}$, further, the orientations of $\Delta_1^j$ and $\Delta_2^j$ both agree with those of $\Delta_{j-1}$ and $\Delta_{j}$ if $j$ is odd; the orientations of $\Delta_1^j$ and $\Delta_2^j$ both disagree with those of $\Delta_{j-1}$ and $\Delta_{j}$ otherwise. We glue $\Delta_1^j$ and $\Delta_2^j$ at the edge labeled $\tau_{i_j}$, so that the orientations of $\Delta_1^j$ and $\Delta_2^j$ both either agree or disagree with those of $\Delta_{j-1}$ and $\Delta_{j}$. We call the edge labeled $\tau_{i_j}$ the \emph{diagonal} of $G(p_j)$.

\medskip

The two arcs $\tau_{i_j}$ and $\tau_{i_{j+1}}$ form two edges of the triangle $\Delta_j$. We denote the third edge of $\Delta_j$ by $\tau_{[\zeta_j]}$. After glue the tiles $G(p_j)$ and $G(p_{j+1})$ at the edge labeled $\tau_{[\zeta_j]}$ for $1\leq j<d-1$ step by step, we obtain a graph, denote as $\overline{G_{T,\gamma}}$. Let $G_{T,\gamma}$ be the graph obtained from $\overline{G_{T,\gamma}}$ by removing the diagonal of each tile. See the following figures for example. In particular, when $\gamma\in T$, let $G_{T,\gamma}$ be the graph with one only edge labeled $\gamma$. Similarly, $G_{T,\gamma}$ can be constructed for any oriented curve $\gamma$ connecting two marked points in $\mathcal O$.

\centerline{\includegraphics{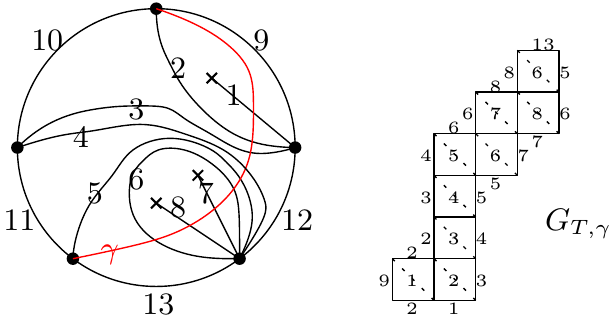}}

\medskip

\centerline{\includegraphics{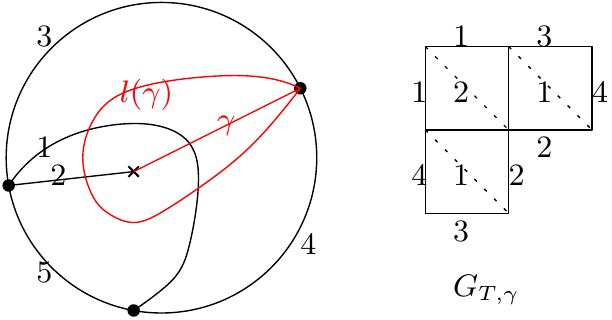}}

\medskip

\begin{Definition}\cite[Definition 4.6]{MSW}\label{perctvec}
A \emph{perfect matching} of a graph $G$ is a subset $P$ of the edges of $G$ such that each vertex of $G$ is incident to exactly one edge of $P$. We denote the set of all perfect matchings of $G$ by $\mathcal P(G)$.

\end{Definition}

\medskip

We have the following observation according to the definition.

\medskip

\begin{Lemma}\label{inone}\cite[Lemma 3.3]{H}
Let $G_i$ and $G_{i+1}$ be two consecutive tiles of $G_{T,\gamma}$ sharing same edge $a$. If $b$ (respectively $c$) is an edge of $G_i$ (respectively $G_{i+1}$) which is incident to $a$, then $b$ and $c$ can not in a perfect matching of $G$ at the same time.

\end{Lemma}

\medskip

\begin{Definition} \cite[Definition 4.7]{MSW}
Let $a_1$ and $a_2$ be the two edges of $\overline{G_{T,\gamma}}$ which lie in the counterclockwise direction from the diagonal of $G(p_1)$. Then the \emph{minimal matching} $P_{-}(G_{T,\gamma})$ is defined as the unique perfect matching which contains only boundary edges and does not contain edges $a_1$ or $a_2$. The \emph{maximal matching} $P_{+}(G_{T,\gamma})$ is the other perfect matching with only boundary edges.

\end{Definition}

\medskip

We have the following observation according to the definition.

\medskip

\begin{Lemma}\label{max-min}\cite[Lemma 3.5]{H}
Let $a$ be an edge of the tile $G(p_j)$. If $a$ is in the maximal/minimal perfect matching of $G_{T,\gamma}$, then $a$ lies in the counterclockwise/clockwise direction from the diagonal of $G(p_j)$ if $j$ is odd and lies in the clockwise/counterclockwise direction from the diagonal of $G(p_j)$ if $j$ is even.

\end{Lemma}

\medskip

\begin{Definition}\cite[Definition 3.6]{H}
Let $P$ be a perfect matching of $G_{T,\gamma}$.

\begin{enumerate}[$(1)$]

  \item \cite{MSW1} We call $P$ can \emph{twist} on a tile $G(p)$ if $G(p)$ has two edges belong to $P$, in such case we define the \emph{twist} $\mu_pP$ of $P$ on $G(p)$ to be the perfect matching obtained from $P$ by replacing the edges in $G(p)$ by the remaining two edges.

  \item In case $P$ can do twist on $G(p)$ with diagonal labeled by $\tau\in T$, we call the pair of the edges of $G(p)$ lying in $P$ a \emph{$\tau$-mutable edges pair in $P$}, any other edge in $P$ is called \emph{non-$\tau$-mutable edge in $P$}.

\end{enumerate}

\end{Definition}

\medskip

\begin{Lemma}\label{transitive}\cite[Lemma 6.3]{H}
For any pair of perfect matchings $P,Q\in \mathcal P(G_{T,\gamma})$, $Q$ can be obtained from $P$ by a sequence of twists.

\end{Lemma}

\medskip

Let $G$ be a snake graph consisting tiles $G(p_1),\cdots,G(p_d)$. Suppose that we truncate $G$ as subgraphs $G_1,\cdots,G_k$ at some non-boundary edges $u_1,\cdots,u_{k-1}$. Particularly, $G$ can be obtained from $G_i,i\in [1,k]$ by gluing $G_i$ and $G_{i+1}$ at the edges $u_i,i\in [1,k-1]$, called \emph{gluing edges}.

\medskip

We denote the edges incident to $u_i$ by $v_i,w_i$ and $v'_i,w'_i$, as shown in the figure below, where $v_i,w_i$ are edges of $G_i$ and $v'_i,w'_i$ are edges of $G_{i+1}$.

\centerline{\begin{tikzpicture}
\draw[-] (1,0) -- (2,0);
\draw[-] (1,-1) -- (2,-1);
\draw[-] (2,0) -- (3,0);
\draw[-] (2,-1) -- (3,-1);
\draw[-] (3,0) -- (3,-1);
\node[above] at (1.5,-0.05) {$v_i$};
\node[below] at (1.5,-0.95) {$w_i$};
\node[right] at (1.75,-0.5) {$u_i$};
\node[above] at (2.5,-0.05) {$v'_i$};
\node[below] at (2.5,-0.95) {$w'_i$};
\draw [-] (1, 0) -- (1,-1);
\draw [-] (2, 0) -- (2,-1);
\draw [fill] (1,0) circle [radius=.05];
\draw [fill] (2,0) circle [radius=.05];
\draw [fill] (1,-1) circle [radius=.05];
\draw [fill] (2,-1) circle [radius=.05];
\draw [fill] (3,0) circle [radius=.05];
\draw [fill] (3,-1) circle [radius=.05];
\end{tikzpicture}}

Given $P\in \mathcal P(G)$, by Lemma \ref{inone}, $u_i\in P$ or $v_i\in P$ or $w'_i\in P$. For each $j\in [1,k]$, we associate $P$ with a perfect matching $\iota_j(P)\in \mathcal P(G_j)$ as follows:

 \[\begin{array}{ccl} \iota_j(P) &=&

         \left\{\begin{array}{ll}

             P\cap E(G_{T,\zeta_j}), &\mbox{if $v_{j-1},w'_{j}\notin P$}, \\

             P\cap E(G_{T,\zeta_j})\cup\{u_{j-1}\}, &\mbox{if $v_{j-1}\in P$ and $w'_{j}\notin P$}, \\

             P\cap E(G_{T,\zeta_j})\cup\{u_{j}\}, &\mbox{if $v_{j-1}\notin P$ and $w'_{j}\in P$}, \\

             P\cap E(G_{T,\zeta_j})\cup\{u_{j-1},u_j\}, &\mbox{if $v_{j-1},w'_{j}\in P$},

         \end{array}\right.

 \end{array}\]
Particularly,

 \[\begin{array}{ccl} \iota_1(P) &=&

         \left\{\begin{array}{ll}

             P\cap E(G_{T,\gamma_1}), &\mbox{if $w'_{1}\notin P$}, \\

             P\cap E(G_{T,\gamma_1})\cup\{u_1\}, &\mbox{if $w'_{1}\in P$}.

         \end{array}\right.

 \end{array}\]

 \[\begin{array}{ccl} \iota_k(P) &=&

         \left\{\begin{array}{ll}

             P\cap E(G_{T,\gamma_1}), &\mbox{if $v_{k-1}\notin P$}, \\

             P\cap E(G_{T,\gamma_1})\cup\{u_{k-1}\}, &\mbox{if $v_{k-1}\in P$}.

         \end{array}\right.

 \end{array}\]

It can be seen that $\iota_j(P)$ is a perfect matching of $G_{j}$ by definition for $j\in [1, k]$.

\medskip

\begin{Proposition}\label{decompose}\cite[Proposition 5.1]{H}
With the same notation as above, we have a bijection
$$\iota: \mathcal P(G)\rightarrow \{(P_i)_{i} \mid P_i\in \mathcal P(G_{i}), \forall i, u_i\in P_i\cup P_{i+1}\},\;\;\;P\rightarrow (\iota_i(P))_{i}.$$
\end{Proposition}

\section{Preliminaries on quantum cluster algebras}\label{pre2}

\subsection{Quantum cluster algebras}

Herein, we recall the definition of quantum cluster algebras in \cite{BZ}. We follow the convention in \cite{Q}. Fix two integers $n\leq m$. Let $\widetilde B$ be an $m\times n$ integer matrix. Let $\Lambda$ be an $m\times m$ skew-symmetric integer matrix. We call $(\widetilde B,\Lambda)$ \emph{compatible} if $(\widetilde B)^t \Lambda=(D\;\;0)$ for some diagonal matrix $D$ with positive entries, where $(\widetilde B)^t$ is the transpose of $\widetilde B$. It should be noted that in this case, the upper $n\times n$ submatrix of $\widetilde B$ is skew-symmetrizable, and $\widetilde B$ is full rank.

\medskip

Let $q$ be the quantum parameter. A \emph{quantum seed} $t$ consists a compatible pair $(\widetilde B,\Lambda)$ and a collection of indeterminate $X_i(t),i\in [1,m]$, called \emph{quantum cluster variables}. Let $\{e_i\}$ be the standard basis of $\mathbb Z^m$ and $X(t)^{e_i}=X_i(t)$. We define the corresponding \emph{quantum torus} $\mathcal T(t)$ to be the algebra which is freely generated by $X(t)^{\vec{a}},\vec{a}\in \mathbb Z^m$ as $\mathbb Z[q^{\pm 1/2}]$-module, with multiplication on these elements defined by
$$X(t)^{\vec{a}}X(t)^{\vec{b}}=q^{\Lambda(t)(\vec{a},\vec{b})/2}X(t)^{\vec{a}+\vec{b}},$$
where $\Lambda(t)(,)$ is the bilinear form on $\mathbb Z^m$ such that
$$\Lambda(t)(e_i,e_j)=\Lambda(t)_{ij}.$$

For any $k\in [1,n]$, we define the \emph{mutation} of $t$ at the $k$-th direction to be the new seed $t'=\mu_k(t)=((X_i(t')_{i\in [1,m]}), \widetilde B(t'),\Lambda(t'))$, where
\begin{enumerate}[$(1)$]

\item $X_i(t')=X_i(t)$ for $i\neq t$,

\item $X_k(t')=X(t)^{-e_k+\sum_i[b_{ik}]_{+}e_i}+X(t)^{-e_k+\sum_i[-b_{ik}]_{+}e_i}$.

\item $\widetilde{B}(t')=\mu_k\widetilde B(t)$.

\item $\Lambda(t')$ is skew-symmetric and satisfies:
\[\begin{array}{ccl} \Lambda(t')_{ij} &=&
         \left\{\begin{array}{ll}
              \Lambda(t)_{ij}, &\mbox{if $i,j\neq k$},  \\
              \Lambda(t)(e_i,-e_k+\sum_{l}[b_{lk}]_{+}e_l), &\mbox{if $i\neq k=j$},
         \end{array}\right.
 \end{array}\]

\end{enumerate}

\medskip

It can be seen that $(\widetilde B(t'), \Lambda(t'))$ is compatible since $\widetilde B^t(t)\Lambda(t)=\widetilde B^t(t')\Lambda(t')$.

\medskip

The quantum torus $\mathcal T(t')$ for the new seed $t'$ is defined similarly.

\medskip

A \emph{quantum cluster algebra} $\mathcal A_q$ is defined as following:
\begin{enumerate}[$(1)$]

  \item Choose an initial seed $t_0=((X_1,\cdots,X_m), B,\Lambda)$.

  \item All the seeds $t$ are obtained from $t_0$ by iterated mutations at directions $k\in [1,n]$.

  \item $\mathcal A_q=\mathbb Z[q^{\pm1/2}]\langle X_i(t)\rangle_{t,i\in [1,m]}$.

  \item $X_{n+1},\cdots, X_m$ are called \emph{frozen variables} or \emph{coefficients}.

  \item A quantum cluster variable in $t$ is called a \emph{quantum cluster variable} of $\mathcal A_q$.

  \item $X(t)^{\vec{a}}$ for some $t$ and $\vec{a}\in \mathbb N^m$ is called a \emph{quantum cluster monomial}.

\end{enumerate}

\medskip

\begin{Theorem}\cite[Quantum Laurent Phenomenon]{BZ}
Let $\mathcal A_q$ be a quantum cluster algebra and $t$ be a seed. For any quantum cluster variable $X$, we have $X\in \mathcal T(t)$.

\end{Theorem}

\medskip

\begin{Conjecture}\cite[Positivity Conjecture]{BZ}
Let $\mathcal A_q$ be a quantum cluster algebra and $t$ be a seed. For any quantum cluster variable $X$ of $\mathcal A_q$,
$$X\in \mathbb N[q^{\pm1/2}]\langle X(t)^{\vec{a}}\mid \vec{a}\in \mathbb Z^m\rangle.$$

\end{Conjecture}

This conjecture was recently proved by \cite{D} in the skew-symmetric case.

\begin{Remark}

In this study, to distinguish commutative cluster algebras and quantum cluster algebras, we use the notation $x$ to denote the commutative cluster variables and $X$ to denote the quantum cluster variables.

\end{Remark}

\subsection{Quantum cluster algebras from unpunctured orbifolds}
We first fix some notation for the rest of this paper. Let $(\mathcal O,M,U)$ be an unpunctured orbifold and $T=\{\tau_1,\cdots,\tau_n,\cdots,\tau_l\}$ be an indexed triangulation. Let $B^T$ be the signed adjacency matrix of $T$. A quantum cluster algebra $\mathcal A_q$ is called \emph{coming from} $\mathcal O$ if there is a quantum seed $t$ such that the upper $n\times n$-submatrix of $\widetilde B(t)$ is $B^T$. By \cite[Theorem 6.1]{BZ} and \cite{FST3}, the quantum seeds/quantum cluster variables of $\mathcal A_q$ and the triangulations/arcs of $\mathcal O$ are one to one correspondence. Let $\Sigma_q^T=(X^{T},\widetilde B^T, \Lambda^T)$ be the quantum seed of $\mathcal A_q$ associate with $T$, where $X^T=\{X^T_{\alpha}\mid \alpha\in T\}\cup \{X_i\mid i\in [n+1,m]\}$ and $X_i,i\in [n+1,m]$ are the coefficients. Let $X_{\gamma}$ be the quantum cluster variable associate with $\gamma$. Let $X_{\gamma}=1$ if $\gamma$ is a boundary arc. For $i\in [1,n]$, denote by $b^{T}_{\tau_i}$ the $i$-th column of $\widetilde B^T$, let $(b^{T}_{\tau_i})_{+}=([b^{T}_{ji}]_{+})_j$ be the positive part of $b^{T}_{\tau_i}$. Dually, let $(b^{T}_{\tau_i})_{-}=([-b^{T}_{ji}]_{+})_j$. Clearly, $b^{T}_{\tau}=(b^{T}_{\tau})_{+}-(b^{T}_{\tau})_{-}$. Set $e_{\tau_i}=e_i\in \mathbb Z^m$ with the $i$-th coordinate 1 and others 0 for $i\in [1,n]$.

\medskip

We denote $q^{-\frac{1}{2}\sum_{i<j}\Lambda^T_{ij}}(X^{T}_{\tau_1})^{a_1}\cdots (X^{T}_{\tau_n})^{a_n}X_{n+1}^{a_{n+1}}\cdots X_{m}^{a_{m}}$ by $(X^T)^{\vec{a}}$ for any $\vec{a}=(a_1,\cdots,a_m)\in \mathbb Z^m$. Therefore the quantum cluster monomials of $\mathcal A_q$ are the forms $(X^T)^{\vec{a}}$ for indexed triangulation $T$ and $\vec{a}\in \mathbb N^m$.

\medskip

We assume $(\widetilde B^T)^t\Lambda=(D^T\;\;0)$ for some $D^T=diag(d^T(\tau))_{\tau\in T}$. As $(\widetilde B^T)^t \Lambda^T=(\widetilde B^{T'})^t \Lambda^{T'}$ and the weight does not change after flip, $d^T(\tau)=d^{T'}(\tau')$ and $d^T(\alpha)=d^{T'}(\alpha)$ for $\alpha\neq \tau$. Let $\gamma$ be an oriented arc in $\mathcal O$ and $\zeta$ be the corresponding ordinary arc which crosses $T$ with points $p_1,\cdots, p_d$ in order. We assume that $p_1,\cdots,p_d$ belong to the arcs $\tau_{i_1},\cdots,\tau_{i_d}$, respectively in $T$. For any $s\in [1,d]$, denote by $m_{p_s}^{+}(\tau_{i_s},\gamma)$ and $m_{p_s}^{-}(\tau_{i_s},\gamma)$ the numbers of $\tau_{i_s}$ in the sets $\{\tau_{i_t}\mid t>s\}$ and $\{\tau_{i_t}\mid t<s\}$, respectively.

\medskip

Let $P$ be a perfect matching of $G_{T,\gamma}$ with edges labeled $\tau_{j_1},\cdots, \tau_{j_r}$ in order. When $P$ can twist on $G(p_s)$, assume that $\tau_{j_t},\tau_{j_{t+1}}$ are edges of $G(p_s)$. We denote by $n_{p_s}^{+}(\tau_{i_s},P)$ and $n_{p_s}^{-}(\tau_{i_s},P)$ the numbers of $\tau_{i_s}$ in $\{\tau_{j_u}\mid u>t+1\}$ and $\{\tau_{j_u}\mid u<t\}$, respectively.

\medskip

Fix $s\in [1,d]$, we assume that $a_{1_s},a_{4_s},\tau_{i_s}$ and $a_{2_s},a_{3_s},\tau_{i_s}$ are two triangles in $T$ such that $a_{1_s},a_{3_s}$ follow $\tau_{i_s}$ in the clockwise order and $a_{2_s},a_{4_s}$ follow $\tau_{i_s}$ in the counterclockwise order.

\medskip

\begin{Definition}\label{omega}

Suppose that $P$ can twist on $G(p_s)$. If the edges labeled $a_{2_s},a_{4_s}$ of $G(p_s)$ are in $P$, let
$$\Omega(p_s, P)=[n_{p_s}^{+}(\tau_{i_s}, P)-m_{p_s}^{+}(\tau_{i_s}, \gamma)-n_{p_s}^{-}(\tau_{i_s}, P)+m_{p_s}^{-}(\tau_{i_s}, \gamma)]d^T(\tau_{i_s}),$$ otherwise, let $$\Omega(p_s, P)=-[n_{p_s}^{+}(\tau_{i_s}, P)-m_{p_s}^{+}(\tau_{i_s}, \gamma)-n_{p_s}^{-}(\tau_{i_s}, P)+m_{p_s}^{-}(\tau_{i_s}, \gamma)]d^T(\tau_{i_s}).$$
\end{Definition}

Clearly, we have $\Omega(p_s,\mu_{p_s}P)=-\Omega(p_s,P)$ if $P$ can twist on $G(p_s)$.

\medskip

\section{Commutative Laurent expansions}\label{cle}

In this section, the main result on commutative cluster algebras from unpunctured orbifolds is stated. More precisely, the Laurent expansion formula of a cluster variable with respect to arbitrary cluster is given.

\medskip

Throughout this section, $(\mathcal O,M,U)$ is an unpunctured orbifold, $T$ is an indexed triangulation of $\mathcal O$. Let $\mathcal A$ be a cluster algebra from $\mathcal O$ with semi-field $\mathbb P$. By \cite{FST3}, the seeds/cluster variables of $\mathcal A$ and the triangulations/arcs of $\mathcal O$ are one to one correspondence. Let $\Sigma^T=(x^T,y^T,B^T)$ be the seed associated with $T$, where $x^{T}=\{x^{T}_{\alpha}\mid \alpha\in T\}$, $y^{T}=\{y^{T}_{\alpha}\mid \alpha\in T\}\subset \mathbb P$. Let $x_{\gamma}$ be the cluster variable associate with $\gamma$. Set $x_{\gamma}=1$ if $\gamma$ is a boundary arc. Fix $\tau\in T$, let $\mu_\tau(T)=T'$.

\medskip

\begin{Definition-Lemma}\cite{MSW}

\begin{enumerate}[$(1)$]

  \item If $\gamma$ is an oriented arc and $\tau_{i_1},\cdots, \tau_{i_d}$ is the sequence of arcs in $T$ which $\gamma$ crosses, then we define the \emph{crossing monomial} of $\gamma$ with respect to $T$ to be $$c(\gamma, T)=\textstyle\prod_{j=1}^dx^T_{\tau_{i_j}}.$$

  \item Let $P\in \mathcal P(G_{T,\gamma})$. If the edges of $P$ are labeled $\tau_{j_1},\cdots, \tau_{j_r}$, then we define the \emph{weight} $w^T(P)$ of $P$ to be $x^T_{\tau_{j_1}}\cdots x^T_{\tau_{j_r}}$.

  \item Let $P\in \mathcal P(G_{T,\gamma})$. The set $(P_{-}(G_{T,\gamma})\cup P)\setminus (P_{-}(G_{T,\gamma})\cap P)$ is the set of boundary edges of a (possibly disconnected) subgraph $G_P$ of $G_{T,\gamma}$, which is a union of cycles. These cycles enclose a set of tiles $\bigcup_{j\in J}G(p_{i_j})$, where $J$ is a finite index set. We define the \emph{height monomial} $y^T(P)$ of $P$ by
      $$y^T(P)=\textstyle\prod_{k=1}^n (y^{T}_{\tau_k})^{m_k},$$
      where $m_k$ is the number of tiles in $\bigcup_{j\in J}G(p_{i_j})$ whose diagonal is labeled $\tau_k$.

\end{enumerate}

\end{Definition-Lemma}

\medskip

It should be noted that $x_{\alpha}=1$ if $\alpha$ is a boundary arc.

\medskip

\begin{Definition}\cite[Definition 4.2]{H1}\label{com-mon}
Let $Q\in \mathcal P(G_{T,\gamma})$. We define the cluster monomial $x^T(Q)$ associated with $Q$ to be
$$x^T(Q)=\frac{w^T(Q)\cdot y^T(Q)}{c(\gamma, T)\cdot \bigoplus_{P\in \mathcal P(G_{T,\gamma})}y^T(P)},$$
it should be noted that the operation $``\oplus"$ in $\bigoplus_{P\in \mathcal P(G_{T,\gamma})}y^T(P)$ is taken in $\mathbb P$.
\end{Definition}

\medskip

We have the following crucial lemma according to the definition of $y^T(P)$.

\medskip

\begin{Lemma}\cite[Lemma 4.3]{H1}\label{brick1}
Let $P\in \mathcal P(G_{T,\gamma})$. Suppose that $P$ can twist on a tile $G(p)$ with diagonal labeled $a$. We assume that the labels of the edges of $G(p)$ are $u_1,u_2,u_3,u_4$ with $u_1,u_3$ are clockwise to $a$ and $u_2,u_3$ are counterclockwise to $a$ in $T$. If the edges labeled $u_1,u_3$ of $G(p)$ are in $P$, then $\frac{y^T(\mu_pP)}{y^T(P)}=y^T_{a}$.

\end{Lemma}

\medskip

\begin{Definition}\cite{H,H1} Let $S,S'$ be two sets. A subset $\mathcal{S}$ of the power set of $S$ is called \emph{a partition} of $S$ if $\cup_{R\in \mathcal S}R=S$ and $R\cap R'=\emptyset$ if $R\neq R'\in \mathcal S$. \emph{A partition map} from a set $S$ to a set $S'$ is a map from a partition of $S$ to a partition of $S'$. A partition map is called \emph{a partition bijection} if it is a bijection.

\end{Definition}

\medskip

\begin{Remark}

To give a partition bijection from $S$ to $S'$ is equivalent to associate each $s\in S$ with a non-empty subset $\pi(s)\subset S'$ which satisfies $\pi(s_1)=\pi(s_2)$ if $\pi(s_1)\cap \pi(s_2)\neq \emptyset$ for $s_1,s_2\in S$ and $\cup_{s\in S}\pi(s)=S'$.

\end{Remark}

\medskip

After the above preparations, we now state the first main result on commutative cluster algebras from unpunctured orbifolds.

\medskip

\begin{Theorem}\label{partition bi}

Let $(\mathcal O,M,U)$ be an unpunctured orbifold and $T$ be an indexed triangulation. We assume that the weights of the orbifold points in $\mathcal O$ are $2$. Let $\gamma$ be an oriented arc in $\mathcal O$. For any $\tau\in T$, there are partitions $\mathfrak P$ and $\mathfrak P'$ of $\mathcal P(G_{T,\gamma})$ and $\mathcal P(G_{T',\gamma})$, respectively, and a bijection $\pi: \mathfrak P\rightarrow \mathfrak P'$ such that

\begin{enumerate}[$(1)$]

  \item $|S|=1$ or $|\pi(S)|=1$ for any $S\in \mathfrak P$.

  \item $\sum_{P\in S}x^T(P)=\sum_{P'\in \pi(S)}x^{T'}(P')$ for any $S\in \mathfrak P$.

\end{enumerate}

\end{Theorem}

\medskip

As a corollary of Theorem \ref{partition bi}, the commutative Laurent expansion formula can be given.

\medskip

\begin{Theorem}\label{expansion-comm}

Let $(\mathcal O,M,U)$ be an unpuntured orbifold and $T$ be an indexed triangulation. We assume that the weights of the orbifold points in $\mathcal O$ are $2$. If $\gamma$ is an oriented arc in $\mathcal O$, then the Laurent expansion of $x_{\gamma}$ with respect to the cluster $x^{T}$ is $$x_{\gamma}=\textstyle\sum_{P\in \mathcal P(G_{T,\gamma})} x^T(P).$$

\end{Theorem}

\begin{proof}

The proof is the same as that of \cite[Theorem 4.8]{H1} by Theorem \ref{partition bi} and induction.
\end{proof}

\medskip

\section{Positivity of the quantum cluster algebras from orbifolds}\label{qle}

In this section, we present the main result on quantum cluster algebras. More precisely, the quantum Laurent expansion of a quantum cluster variable with respect to arbitrary quantum cluster of $\mathcal A_q$ is given. As an application, the positivity for such class of quantum cluster algebras is proved.

\medskip

When specialize $q=1$, we obtain a commutative cluster algebra $\mathcal A_q\mid_{q=1}$. Let $P\in \mathcal P(G_{T,\gamma})$. For the element $x^T(P)$ associate with $P$ in Definition \ref{com-mon}, there clearly exists a unique $\vec{a}(P)\in \mathbb Z^m$ such that $(X^T)^{\vec{a}(P)}\mid_{q=1}=x^T(P)$. We denote $(X^T)^{\vec{a}(P)}$ by $X^T(P)$ and the quantum cluster variable associate with $\gamma$ by $X_{\gamma}$.

\medskip

The following is the main result on quantum cluster algebras from unpunctured orbifolds.

\medskip

\begin{Theorem}\label{mainthm}

Let $(\mathcal O,M,U)$ be an unpunctured orbifold and $T$ be an indexed triangulation. We assume that the weights of the orbifold points in $\mathcal O$ are $2$. If $\gamma$ is an oriented arc in $\mathcal O$, then

\begin{enumerate}[$(1)$]

  \item  there uniquely exists a \emph{valuation} map $v:\mathcal P(G_{T,\gamma})\rightarrow \mathbb Z$ such that

\begin{enumerate}[$(a)$]

  \item (initial conditions) $v(P_{+}(G_{T,\gamma}))=v(P_{-}(G_{T,\gamma}))=0$.

  \item (iterated relation) If $P\in \mathcal P(G_{T,\gamma}))$ can twist on a tile $G(p)$, then $$v(P)-v(\mu_pP)=\Omega(p,P).$$

\end{enumerate}

  \item Furthermore, if $\tau\in T$ and $\pi$ is the partition bijection from $\mathcal P(G_{T,\gamma})$ to $\mathcal P(G_{\mu_\tau T, \gamma})$ given in Theorem \ref{partition bi}, then for any $S\in \mathfrak P$,
      $$\textstyle\sum_{P\in S}q^{v(P)/2}X^T(P)=\textstyle\sum_{P'\in \pi(S)}q^{v(P')/2}X^{T'}(P').$$

\end{enumerate}

\end{Theorem}

\medskip

As a corollary of Theorem \ref{mainthm}, the quantum Laurent expansion formula is given.

\medskip

\begin{Theorem}\label{expansion}

Let $(\mathcal O,M)$ be an unpuntured orbifold and $T$ be an indexed triangulation. We assume that the weights of the orbifold points are $2$. If $\gamma$ is an oriented arc in $\mathcal O$, then the quantum Laurent expansion of $X_{\gamma}$ with respect to the quantum cluster $X^{T}$ is $$X_{\gamma}=\textstyle\sum_{P\in \mathcal P(G_{T,\gamma})} q^{v(P)/2}X^T(P).$$

\end{Theorem}

\begin{proof}

The proof is the same as that of \cite[Theorem 5.2]{H1} by Theorems \ref{partition bi}, \ref{mainthm} and induction.
\end{proof}

\medskip

As an corollary of Theorem \ref{expansion}, the positivity can be verified.

\medskip

\begin{Theorem}\label{positive}

Let $\mathcal O$ be a orbifold without punctures. Assume the weights of the orbifold points are $2$. The positivity holds for the quantum cluster algebra $\mathcal A_q(\mathcal O)$.

\end{Theorem}

\begin{proof}

Since $q^{v(P)/2}\in \mathbb N[q^{\pm 1}]$ for any $P\in \mathcal P(G_{T,\gamma})$, by Theorem \ref{expansion}, the result follows.
\end{proof}

\medskip

\begin{Example}\label{example}

Let $(\mathcal O, M, U)$, $T$ and $\gamma$ be as shown in the following figure. We assume that the extended matrix associate with $T$ is $\widetilde B^T=((B^T)^t\;\;I_3)^t$ and $\Lambda^T=\left(\begin{array}{cc}
0 & -D\\
D & -DB^T
\end{array}\right)$, where $D=diag(2,2,1)$ the $3\times 3$ diagonal matrix with diagonal entries $2,2,1$.

\medskip

\centerline{\includegraphics{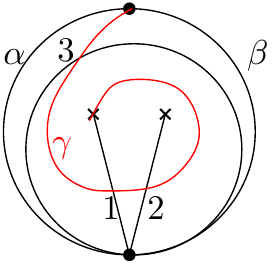}}

The arcs are labeled $1,2,3,\alpha,\beta$. We assume that the frozen variables are $x_4,x_5$ and $x_6$. Thus, $G_{T,\gamma}$ is the following graph and $v:\mathcal P(G_{T,\gamma})\rightarrow \mathbb Z$ is shown in Figure 1.

\centerline{\begin{tikzpicture}
\draw[-] (0,-2) -- (1,-2);
\draw[-] (0,-2) -- (0,-1);
\draw[-] (1,-2) -- (1,-1);
\draw[-] (4,1) -- (5,1);
\draw[-] (4,0) -- (4,1);
\draw[-] (5,0) -- (5,1);
\draw[-] (0,0) -- (1,0);
\draw[-] (0,0) -- (0,-1);
\draw[-] (0,-1) -- (1,-1);
\draw[-] (1,0) -- (2,0);
\draw[-] (1,-1) -- (2,-1);
\draw[-] (2,0) -- (3,0);
\draw[-] (2,-1) -- (3,-1);
\draw[-] (3,0) -- (3,-1);
\draw[-] (3,0) -- (4,0);
\draw[-] (3,-1) -- (4,-1);
\draw[-] (4,0) -- (4,-1);
\draw[-] (4,0) -- (5,0);
\draw[-] (4,-1) -- (5,-1);
\draw[-] (5,0) -- (5,-1);
\draw[dashed] (1,0) -- (2,-1);
\draw[dashed] (0,0) -- (1,-1);
\draw[dashed] (0,-1) -- (1,-2);
\draw[dashed] (2,0) -- (3,-1);
\draw[dashed] (3,0) -- (4,-1);
\draw[dashed] (4,0) -- (5,-1);
\draw[dashed] (4,1) -- (5,0);
\node[above] at (0.5,-0.05) {$2$};
\node[above] at (1.5,-0.05) {$1$};
\node[above] at (2.5,-0.05) {$2$};
\node[above] at (3.5,-0.05) {$1$};
\node[above] at (4.5,-0.25) {$2$};
\node[below] at (4.5,-0.95) {$2$};
\node[below] at (0.5,-0.75) {$2$};
\node[below] at (1.5,-0.95) {$1$};
\node[below] at (2.5,-0.95) {$2$};
\node[below] at (3.5,-0.95) {$1$};
\node[right] at (1.4,-0.45) {$2$};
\node[right] at (3.4,-0.45) {$2$};
\node[right] at (2.4,-0.5) {$1$};
\node[right] at (4.4,-0.5) {$1$};
\node[right] at (4.4,0.5) {$3$};
\node[left] at (0.8,-0.5) {$1$};
\node[left] at (0.8,-1.5) {$3$};
\node[left] at (0.05,-1.5) {$\beta$};
\node[right] at (0.95,-1.5) {$1$};
\node[below] at (0.5,-1.95) {$\alpha$};
\node[left] at (4.05,0.5) {$1$};
\node[right] at (4.95,0.5) {$\beta$};
\node[above] at (4.5,0.95) {$\alpha$};
\node[left] at (0.05,-0.5) {$3$};
\node[left] at (1.25,-0.5) {$3$};
\node[left] at (2.25,-0.5) {$3$};
\node[left] at (3.25,-0.5) {$3$};
\node[left] at (4.25,-0.5) {$3$};
\node[right] at (4.95,-0.5) {$3$};
\draw [-] (1, 0) -- (1,-1);
\draw [-] (2, 0) -- (2,-1);
\draw [fill] (0,0) circle [radius=.05];
\draw [fill] (0,-1) circle [radius=.05];
\draw [fill] (1,0) circle [radius=.05];
\draw [fill] (2,0) circle [radius=.05];
\draw [fill] (1,-1) circle [radius=.05];
\draw [fill] (2,-1) circle [radius=.05];
\draw [fill] (3,0) circle [radius=.05];
\draw [fill] (3,-1) circle [radius=.05];
\draw [fill] (4,0) circle [radius=.05];
\draw [fill] (4,-1) circle [radius=.05];
\draw [fill] (5,0) circle [radius=.05];
\draw [fill] (5,-1) circle [radius=.05];
\draw [fill] (0,-2) circle [radius=.05];
\draw [fill] (1,-2) circle [radius=.05];
\draw [fill] (4,1) circle [radius=.05];
\draw [fill] (5,1) circle [radius=.05];
\end{tikzpicture}}

\medskip

\begin{figure}[h]
\centering
\includegraphics{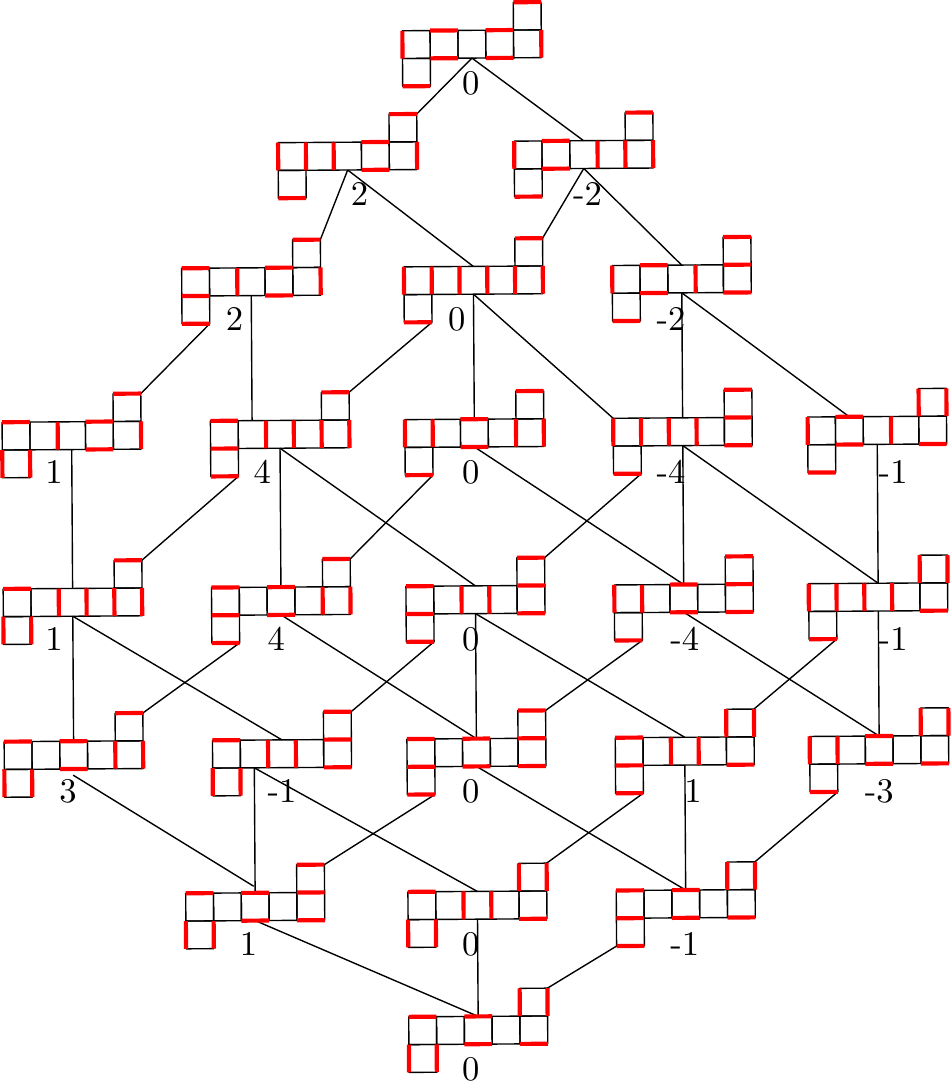}

 {\rm Figure 1 for Example \ref{example}}
\end{figure}

Therefore, $$
\begin{array}{rcl} X_{\gamma}

& = & x^{e_1-2e_2+3e_4+2e_5+2e_6}+(q^{2/2}+q^{-2/2})x^{-e_1-2e_2+2e_3+3e_4+e_5+2e_6} \vspace{2pt}  \\

& + & (q^{2/2}+q^{-2/2})x^{-e_1+2e_4+e_5+2e_6}+x^{-3e_1-2e_2+4e_3+3e_4+2e_6} \vspace{2pt}  \\

& + & (q^{1/2}+q^{-1/2})x^{-e_2+2e_4+e_5+e_6}+(q^{4/2}+1+q^{-4/2})x^{-3e_1+2e_3+2e_4+2e_6}  \vspace{2pt}  \\

& + & (q^{1/2}+q^{-1/2})x^{-2e_1-e_2+2e_3+2e_4+e_6}+(q^{4/2}+1+q^{-4/2})x^{-3e_1+2e_2+e_4+2e_6} \vspace{2pt}  \\

& + & (q^{3/2}+q^{1/2}+q^{-1/2}+q^{-3/2})x^{-2e_1+e_2+e_4+e_6}+x^{-3e_1+4e_2-2e_3+2e_6} \vspace{2pt}  \\

& + & (q^{1/2}+q^{-1/2})x^{-2e_1+3e_2-2e_3+e_6}+x^{-e_1+e_4}+x^{-e_1+2e_2-2e_3}.
\end{array}$$

\end{Example}

\medskip

\section{Partition bijection $\psi_{\rho}$ between $\mathcal P(G_{T,\rho})$ and $\mathcal P(G_{T',\rho})$}\label{compare}

Let $T$ be an indexed triangulation of $\mathcal O$ and $\zeta$ be an oriented ordinary curve in $\mathcal O$. We assume that $\tau$ is a non-boundary arc in $T$ and denote by $\tau'$ the arc obtained from $T$ by flip at $\tau$. Throughout this section, denote the labels of edges of the tile with diagonal labeled $\tau$ by $a_1,a_2,a_3,a_4$ so that $a_1,a_3$ follow $\tau$ in the clockwise order and $a_2,a_4$ follow $\tau$ in the counterclockwise order in $\mathcal O$.

\medskip

Let $p_0$ be the starting point of $\zeta$, and $p_{d+1}$ be its endpoint. We assume that $\zeta$ crosses $T$ at $p_1,\cdots, p_d$ in order. Thus, $p_j,j\in [1,d]$ divide $\zeta$ into some segments. We assume $p_j\in \tau_{i_j}\in T$. For $j\in [1,d-1]$, $\tau_{i_j}$ and $\tau_{i_{j+1}}$ form two edges of a triangle $\Delta_j$ of $T$ such that the segment connecting $p_j$ and $p_{j+1}$ lies inside of $\Delta_j$. We denote the third edge of $\Delta_j$ by $\tau_{[\zeta_j]}$. We choose a point on the segment connecting $p_j$ and $p_{j+1}$ if one of the following two cases happens

\begin{enumerate}[$(1)$]

  \item one of $\tau_{i_j}$ and $\tau_{i_{j+1}}$ is in a same triangle with $\tau$, the other one is not.

  \item $\tau_{i_j},\tau_{i_{j+1}}\neq \tau$, $\tau_{i_j}$ and $\tau_{i_{j+1}}$ are in a same triangle with $\tau$ but $\tau_{[\zeta_j]}\neq \tau$.

\end{enumerate}

\medskip

We denote the chosen points by $o_1,o_2,\cdots,o_{k-1}$ in order and let $p_0=o_0, p_{d+1}=o_{k}$. We denote the subcurve connecting $o_{i-1}$ and $o_i$ of $\zeta$ by $\zeta_{i}$ for $i\in [1,k]$. See Figure 2 for example.

\medskip

\begin{figure}[h]\centering

\includegraphics{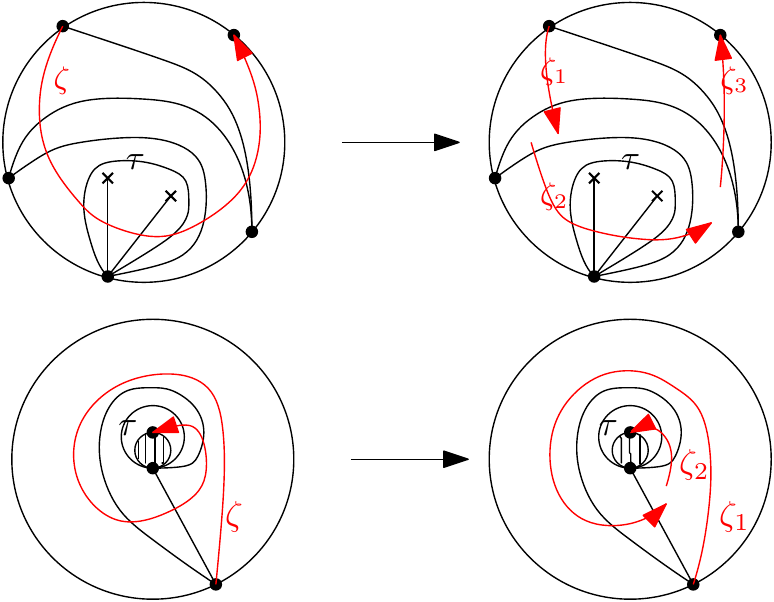}

{\rm Figure 2}
\end{figure}

\medskip

Let $\rho$ be one of $\zeta_i,i\in [1,k]$. In this section, we compare $\mathcal P(G_{T,\rho})$ with $\mathcal P(G_{\mu_{\tau}(T),\rho})$. It should be noted that $G_{T,\rho}=G_{\mu_{\tau}(T),\rho}$ if $\rho$ crosses no arc which is in the same triangle with $\tau$. Herein, we may assume that $\tau'\neq \rho\notin T$ and $\rho$ crosses at least one arc which is in the same triangle with $\tau$.

\medskip

We denote the first tile and the last tile of $G_{T,\rho}$ by $G_1$ and $G_s$, respectively. Then $G_{T,\zeta}$ can be obtained by gluing some graph $G$ left or below to $G_1$ and some graph $G'$ right or upper to $G_s$. It should be noted that $G$ is empty if and only if $\rho$ and $\zeta$ have the same starting point, $G'$ is empty if and only if $\rho$ and $\zeta$ have the same endpoint. When $G$ is not empty, we say that the left almost or lower almost edge of $G_1$ is the \emph{first gluing edge} of $G_{T,\rho}$, and when $G'$ is not empty, we say that the right almost or upper almost edge of $G_s$ is the \emph{last gluing edge} of $G_{T,\rho}$. Similarly, we can define the first/last gluing edge of $G_{T',\rho}$. We assume that $\rho$ crosses $T$ at $p_j,\cdots,p_r$ for some $1\leq j\leq r\leq d$. It can be seen the labels of the first/last gluing edges of $G_{T,\rho}$ and $G_{\mu_{\tau}(T),\rho}$ are the same, more precisely, the first gluing edges are labeled $\tau_{[\zeta_{j-1}]}$ and the last gluing edges are labeled $\tau_{[\zeta_r]}$.

\medskip

Two edges of $G_{T,\rho}$ labeled $\tau$ are called \emph{equivalent} if they are incident to a same diagonal of a tile of $G_{T,\rho}$. Each equivalence class is called a \emph{$\tau$-equivalence class} in $G_{T,\rho}$. We divide the edges of $G_{T,\rho}$ labeled $\tau$ into the following types, as shown in the following graphs (up to a relative orientation).

\begin{enumerate}[$(I)$]

  \item there are two non-incident edges in its equivalence class.

  \item there are two incident edges in its equivalence class.

  \item there is one edge in its equivalence class and it is incident to a diagonal.

  \item there is one edge in its equivalence class and it is not incident to a diagonal.

\end{enumerate}

Type (I):

\centerline{\begin{tikzpicture}
\draw[-] (0,0) -- (1,0);
\draw[-] (0,0) -- (0,-1);
\draw[-] (0,-1) -- (1,-1);
\draw[-] (1,0) -- (2,0);
\draw[-] (1,-1) -- (2,-1);
\draw[-] (2,0) -- (3,0);
\draw[-] (2,-1) -- (3,-1);
\draw[-] (3,0) -- (3,-1);
\draw[dashed] (1,0) -- (2,-1);
\draw[dashed] (0,0) -- (1,-1);
\draw[dashed] (2,0) -- (3,-1);
\node[right] at (1.4,-0.45) {$\tau$};
\node[above] at (0.5,0) {$\tau$};
\node[below] at (2.5,-1) {$\tau$};
\draw [-] (1, 0) -- (1,-1);
\draw [-] (2, 0) -- (2,-1);
\draw [fill] (0,0) circle [radius=.05];
\draw [fill] (0,-1) circle [radius=.05];
\draw [fill] (1,0) circle [radius=.05];
\draw [fill] (2,0) circle [radius=.05];
\draw [fill] (1,-1) circle [radius=.05];
\draw [fill] (2,-1) circle [radius=.05];
\draw [fill] (3,0) circle [radius=.05];
\draw [fill] (3,-1) circle [radius=.05];
%\draw[-] (5,0) -- (6,0);
%\draw[-] (5,-1) -- (6,-1);
%\draw [-] (5, 0) -- (5,-1);
%\draw [-] (6, 0) -- (6,-1);
%\draw[-] (5,0) -- (5,1);
%\draw[-] (5,-1) -- (5,-2);
%\draw[-] (6,0) -- (6,1);
%\draw[-] (6,-1) -- (6,-2);
%\draw[-] (5,1) -- (6,1);
%\draw[-] (5,-2) -- (6,-2);
%\draw[dashed] (5,0) -- (6,-1);
%\draw[dashed] (5,1) -- (6,0);
%\draw[dashed] (5,-1) -- (6,-2);
%\node[right] at (6,-1.5) {$\tau$};
%\node[left] at (5.9,-0.4) {$\tau$};
%\node[left] at (5,0.5) {$\tau$};
%\draw [fill] (5,0) circle [radius=.05];
%\draw [fill] (6,0) circle [radius=.05];
%\draw [fill] (5,-1) circle [radius=.05];
%\draw [fill] (6,-1) circle [radius=.05];
%\draw [fill] (6,-2) circle [radius=.05];
%\draw [fill] (5,-2) circle [radius=.05];
%\draw [fill] (5,1) circle [radius=.05];
%\draw [fill] (6,1) circle [radius=.05];
\end{tikzpicture}}

Type (II):

\centerline{\begin{tikzpicture}
\draw[-] (1,0) -- (2,0);
\draw[-] (1,-1) -- (2,-1);
\draw[-] (2,0) -- (3,0);
\draw[-] (2,-1) -- (3,-1);
\draw[-] (3,0) -- (3,-1);
\draw[-] (1,-1) -- (1,-2);
\draw[-] (2,-1) -- (2,-2);
\draw[-] (1,-2) -- (2,-2);
\draw[dashed] (1,0) -- (2,-1);
\draw[dashed] (2,0) -- (3,-1);
\draw[dashed] (1,-1) -- (2,-2);
\node[right] at (1.4,-0.45) {$\tau$};
\node[right] at (1.95,-1.5) {$\tau$};
\node[below] at (2.5,-1) {$\tau$};
\draw [-] (1, 0) -- (1,-1);
\draw [-] (2, 0) -- (2,-1);
\draw [fill] (1,0) circle [radius=.05];
\draw [fill] (2,0) circle [radius=.05];
\draw [fill] (1,-1) circle [radius=.05];
\draw [fill] (2,-1) circle [radius=.05];
\draw [fill] (1,-2) circle [radius=.05];
\draw [fill] (2,-2) circle [radius=.05];
\draw [fill] (3,0) circle [radius=.05];
\draw [fill] (3,-1) circle [radius=.05];
%\draw[-] (6,0) -- (7,0);
%\draw[-] (6,-1) -- (7,-1);
%\draw[-] (7,0) -- (7,-1);
%\draw[-] (7,-1) -- (7,-2);
%\draw[-] (6,-2) -- (7,-2);
%\draw[-] (5,-1) -- (6,-1);
%\draw [-] (6, 0) -- (6,-1);
%\draw[-] (6,-1) -- (7,-1);
%\draw[-] (5,-2) -- (6,-2);
%\draw[-] (5,-1) -- (5,-2);
%\draw[-] (6,-1) -- (6,-2);
%\draw[dashed] (5,-1) -- (6,-2);
%\draw[dashed] (6,0) -- (7,-1);
%\draw[dashed] (6,-1) -- (7,-2);
%\node[above] at (5.5,-1.05) {$\tau$};
%\node[left] at (6.05,-0.5) {$\tau$};
%\node[right] at (6.4,-1.45) {$\tau$};
%\draw [fill] (7,0) circle [radius=.05];
%\draw [fill] (7,-2) circle [radius=.05];
%\draw [fill] (7,-1) circle [radius=.05];
%\draw [fill] (6,0) circle [radius=.05];
%\draw [fill] (5,-1) circle [radius=.05];
%\draw [fill] (6,-1) circle [radius=.05];
%\draw [fill] (5,-2) circle [radius=.05];
%\draw [fill] (6,-2) circle [radius=.05];
\end{tikzpicture}}

Type (III)

\centerline{\begin{tikzpicture}
\draw[-] (1,0) -- (2,0);
\draw[-] (1,-1) -- (2,-1);
\draw[-] (1,-1) -- (1,-2);
\draw[-] (2,-1) -- (2,-2);
\draw[-] (1,-2) -- (2,-2);
\draw[dashed] (1,0) -- (2,-1);
\draw[dashed] (1,-1) -- (2,-2);
\node[right] at (1.4,-0.45) {$\tau$};
\node[right] at (1.95,-1.5) {$\tau$};
\draw [-] (1, 0) -- (1,-1);
\draw [-] (2, 0) -- (2,-1);
\draw [fill] (1,0) circle [radius=.05];
\draw [fill] (2,0) circle [radius=.05];
\draw [fill] (1,-1) circle [radius=.05];
\draw [fill] (2,-1) circle [radius=.05];
\draw [fill] (1,-2) circle [radius=.05];
\draw [fill] (2,-2) circle [radius=.05];
%\draw[-] (3,0) -- (4,0);
%\draw[-] (3,-1) -- (4,-1);
%\draw[-] (4,0) -- (5,0);
%\draw[-] (4,-1) -- (5,-1);
%\draw[-] (5,0) -- (5,-1);
%\draw[dashed] (3,0) -- (4,-1);
%\draw[dashed] (4,0) -- (5,-1);
%\node[right] at (3.4,-0.45) {$\tau$};
%\node[below] at (4.5,-1) {$\tau$};
%\draw [-] (3, 0) -- (3,-1);
%\draw [-] (4, 0) -- (4,-1);
%\draw [fill] (3,0) circle [radius=.05];
%\draw [fill] (4,0) circle [radius=.05];
%\draw [fill] (3,-1) circle [radius=.05];
%\draw [fill] (4,-1) circle [radius=.05];
%\draw [fill] (5,0) circle [radius=.05];
%\draw [fill] (5,-1) circle [radius=.05];
\draw[-] (6,0) -- (7,0);
\draw[-] (6,-1) -- (7,-1);
\draw[-] (7,0) -- (7,-1);
\draw[-] (7,-1) -- (7,-2);
\draw[-] (6,-2) -- (7,-2);
\draw [-] (6, 0) -- (6,-1);
\draw[-] (6,-1) -- (7,-1);
\draw[-] (6,-1) -- (6,-2);
\draw[dashed] (6,0) -- (7,-1);
\draw[dashed] (6,-1) -- (7,-2);
\node[left] at (6.05,-0.5) {$\tau$};
\node[right] at (6.4,-1.45) {$\tau$};
\draw [fill] (7,0) circle [radius=.05];
\draw [fill] (7,-2) circle [radius=.05];
\draw [fill] (7,-1) circle [radius=.05];
\draw [fill] (6,0) circle [radius=.05];
\draw [fill] (6,-1) circle [radius=.05];
\draw [fill] (6,-2) circle [radius=.05];
%\draw[-] (9,-1) -- (10,-1);
%\draw[-] (10,-1) -- (10,-2);
%\draw[-] (9,-2) -- (10,-2);
%\draw[-] (8,-1) -- (9,-1);
%\draw[-] (9,-1) -- (10,-1);
%\draw[-] (8,-2) -- (9,-2);
%\draw[-] (8,-1) -- (8,-2);
%\draw[-] (9,-1) -- (9,-2);
%\draw[dashed] (8,-1) -- (9,-2);
%\draw[dashed] (9,-1) -- (10,-2);
%\node[above] at (8.5,-1.05) {$\tau$};
%\node[right] at (9.4,-1.45) {$\tau$};
%\draw [fill] (10,-2) circle [radius=.05];
%\draw [fill] (10,-1) circle [radius=.05];
%\draw [fill] (8,-1) circle [radius=.05];
%\draw [fill] (9,-1) circle [radius=.05];
%\draw [fill] (8,-2) circle [radius=.05];
%\draw [fill] (9,-2) circle [radius=.05];
\end{tikzpicture}}

Type (IV)

\centerline{\begin{tikzpicture}
\draw[-] (1,0) -- (2,0);
\draw[-] (1,-1) -- (2,-1);
\draw[-] (1,0) -- (1,-1);
\draw[-] (2,0) -- (2,-1);
\draw[dashed] (1,0) -- (2,-1);
\draw [fill] (1,0) circle [radius=.05];
\draw [fill] (2,0) circle [radius=.05];
\draw [fill] (1,-1) circle [radius=.05];
\draw [fill] (2,-1) circle [radius=.05];
\node[left] at (1.05,-0.5) {$\tau$};
%\draw[-] (3,0) -- (4,0);
%\draw[-] (3,-1) -- (4,-1);
%\draw[-] (3,0) -- (3,-1);
%\draw[-] (4,0) -- (4,-1);
%\draw[dashed] (3,0) -- (4,-1);
%\draw [fill] (3,0) circle [radius=.05];
%\draw [fill] (4,0) circle [radius=.05];
%\draw [fill] (3,-1) circle [radius=.05];
%\draw [fill] (4,-1) circle [radius=.05];
%\node[right] at (3.98,-0.45) {$\tau$};
\draw[-] (5,0) -- (6,0);
\draw[-] (5,-1) -- (6,-1);
\draw[-] (5,0) -- (5,-1);
\draw[-] (6,0) -- (6,-1);
\draw[dashed] (5,0) -- (6,-1);
\draw [fill] (5,0) circle [radius=.05];
\draw [fill] (6,0) circle [radius=.05];
\draw [fill] (5,-1) circle [radius=.05];
\draw [fill] (6,-1) circle [radius=.05];
\node[above] at (5.5,-0.05) {$\tau$};
%\draw[-] (7,0) -- (8,0);
%\draw[-] (7,-1) -- (8,-1);
%\draw[-] (7,0) -- (7,-1);
%\draw[-] (8,0) -- (8,-1);
%\draw[dashed] (7,0) -- (8,-1);
%draw [fill] (7,0) circle [radius=.05];
%\draw [fill] (8,0) circle [radius=.05];
%\draw [fill] (7,-1) circle [radius=.05];
%\draw [fill] (8,-1) circle [radius=.05];
%\node[below] at (7.5,-0.95) {$\tau$};
\end{tikzpicture}}

\medskip

We denote by $n^{\tau}(T,\rho)$ the number of $\tau$-equivalence classes in $G_{T,\rho}$ and list the $\tau$-equivalence classes according to the order of the tiles. Let $S_{T,\rho}$ be the set containing all sequences $\nu=(\nu_{1},\cdots,\nu_{n^{\tau}(T,\rho)})$ which satisfies $\nu_{j}\in \{-1,0,1\}$ if the $j$-th $\tau$-equivalence class is of type (I), $\nu_{j}\in \{-1,0\}$ if the $j$-th $\tau$-equivalent class is of type (II) or (III), $\nu_{j}\in \{0,1\}$ if the $j$-th $\tau$-equivalence class is of type (IV). Given a sequence $\nu_=(\nu_{1},\cdots,\nu_{n^{\tau}(T,\rho)})\in S_{T,\rho}$, let $\mathcal P^{\tau}_{\nu}(G_{T,\rho})$ containing all perfect matching $P$ which contains $\nu_{j}+1$ edges of the $j$-th $\tau$-equivalence class if the $j$-th $\tau$-equivalence class is of type (I) (II) or (III), and $\nu_{j}$ edges of the $j$-th $\tau$-equivalence class if the $j$-th $\tau$-equivalence class is of type (IV).

\medskip

It should be noted that if $\nu_j=-1$, then the $j$-th $\tau$-equivalence class corresponds to a diagonal of a tile $G(p)$, moreover, any $P\in \mathcal P^{\tau}_{\nu}(G_{T,\rho})$ can twist on $G(p)$.

\medskip

It can be easily seen that the following result is true.

\medskip

\begin{Lemma}\label{localdec}

$\mathcal P(G_{T,\rho})=\bigsqcup_{\nu\in S_{T,\rho}}\mathcal P^{\tau}_{\nu}(G_{T,\rho})$.

\end{Lemma}

\medskip

\subsection{``Locally" unfolding}\label{local1}

The concept of unfolding of a skew-symmetrizable matrix was introduced by Zelevinsky and the authors of \cite{FST3} and first appeared in \cite{FST3}. It turns out that unfolding is a useful method to study skew-symmetrizable cluster algebras, see \cite{FST3,D,HL}. We do not give the detail definition in this paper. We use the idea of unfolding to compare $\mathcal P(G_{T,\rho})$ and $\mathcal P(G_{T',\rho})$. Herein, it can be seen that it suffices to consider the case that the surfaces without orbifold points.

\medskip

Let $\mathcal O_0$ be the subsurface of $\mathcal O$ formed by the triangle or two triangles containing $\tau$. We can lift $\mathcal O_0$ to a surface $\widetilde{\mathcal O}_0$, as show in the following tables.

\medskip

\centerline{\includegraphics{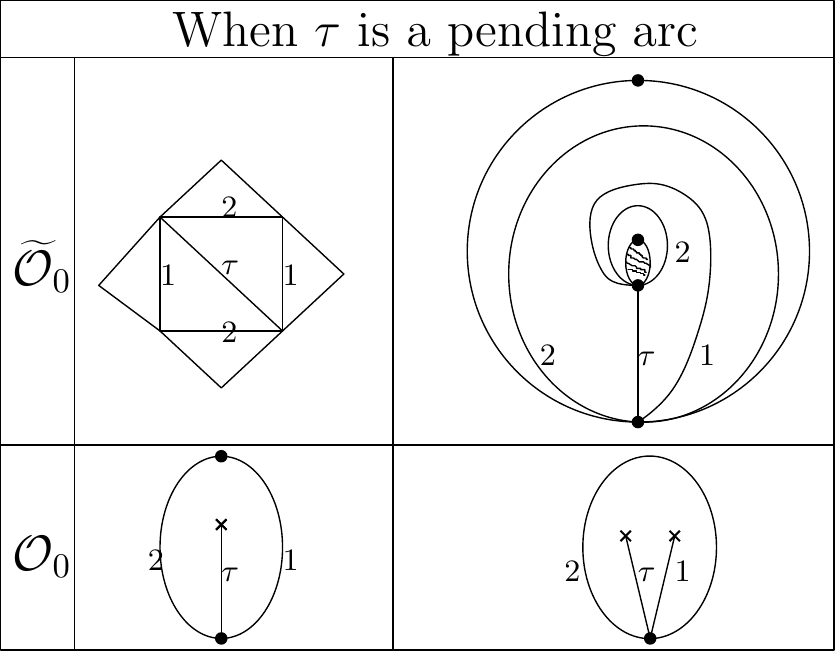}}

\medskip

\centerline{\includegraphics{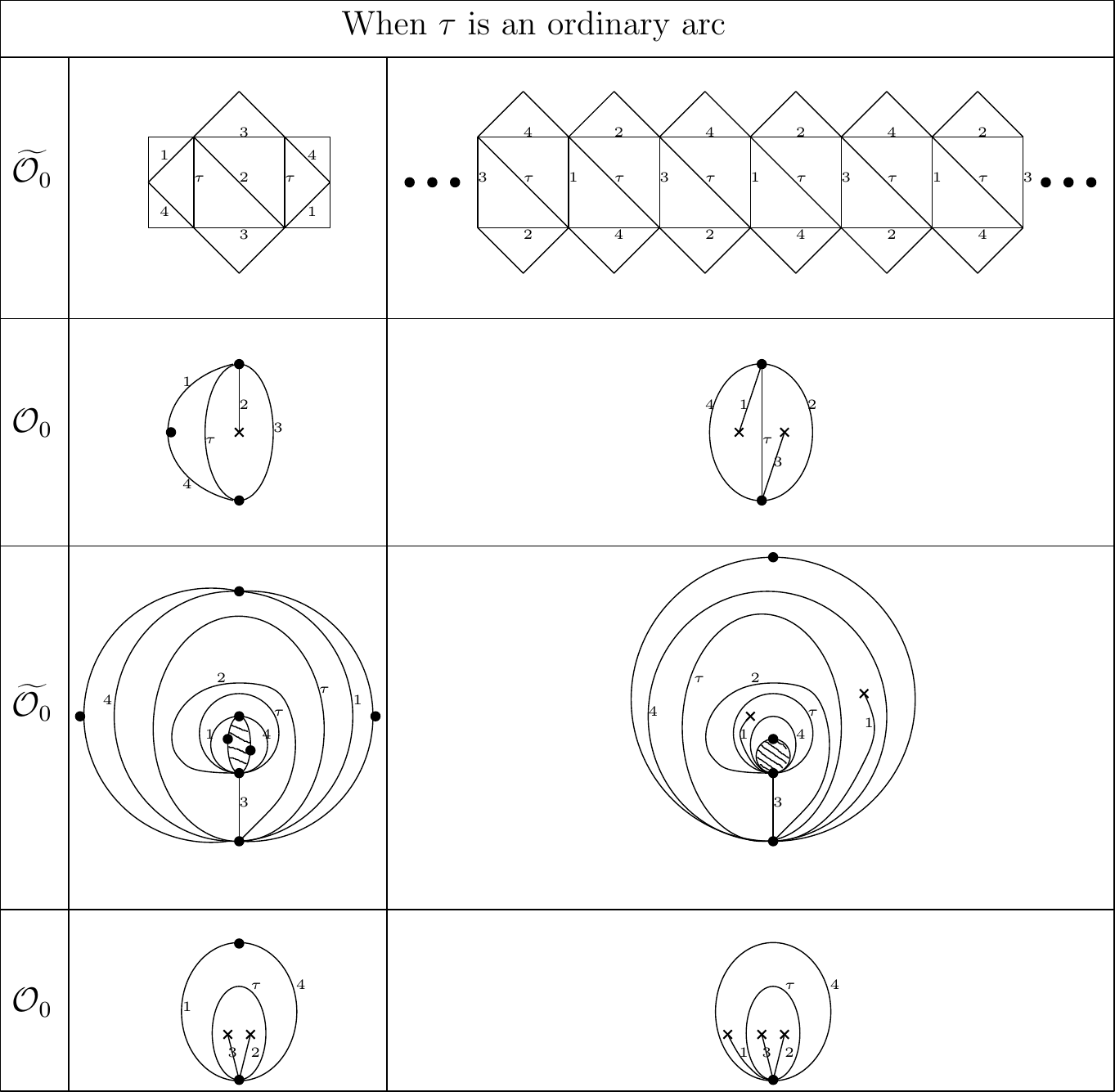}}

\medskip

Obviously, the restriction of $T$ on $\mathcal O_0$ gives a triangulation of $\mathcal O_0$, denote by $T_0$. It can be lift to a triangulation of $\widetilde{\mathcal O}_0$, denote by $\widetilde T_0$, as shown in the above tables. We use the notation $\widetilde \gamma$ to represent an arc in $\widetilde T_0$ and denote the corresponding arc in $T_0$ by $\chi(\widetilde\gamma)$. It should be noted that we do not consider the boundary arcs of $\widetilde{\mathcal O}_0$. For all above cases, any arc in $\mathcal O_0$ can be lift to at most two arcs in $\widetilde{\mathcal O}_0$.

\medskip

For a given curve $\gamma\in \mathcal O_0$, we assume that $\gamma$ crosses the arcs $\gamma_1,\cdots,\gamma_d$ in $T_0$ in order. Then $\gamma$ can be lift to a curve $\widetilde \gamma\in \widetilde{\mathcal O}_0$ which crosses the arcs $\widetilde {\gamma}_1,\cdots,\widetilde{\gamma}_d$ in $\widetilde T_0$ in order with $\chi(\widetilde{\gamma}_i)=\gamma_i$ for $i\in [1,d]$. It should be noted that we do not know that whether $\widetilde{\gamma}$ cross itself even if $\gamma$ does not cross itself.

\medskip

Let $\rho$ be the curve defined at the beginning of this section and $\widetilde \rho$ be a curve in $\widetilde{\mathcal O}_0$ lift by $\rho$. Let $\widetilde T'_0$ be the triangulation of $\widetilde {\mathcal O}_0$ obtained from $\widetilde T_0$ by flip at the arcs $\widetilde \tau$ with $\chi(\widetilde \tau)=\tau$. Thus, $G_{T,\rho}\cong G_{\widetilde T_0,\widetilde \rho}$ and $G_{T',\rho}\cong G_{\widetilde T'_0,\widetilde \rho}$ as graphs. Therefore, we can compare $\mathcal P(G_{T,\rho})$ and $\mathcal P(G_{T',\rho})$ via the comparison of $\mathcal P(G_{\widetilde T_0,\widetilde \rho})$ and $\mathcal P(G_{\widetilde T'_0,\widetilde \rho})$.

\subsection{$\mathcal O_0$ contains no orbifold points.}\label{local2}
Herein, we consider the case that $\mathcal O_0$ contains no orbifold points. Thus, $\tau$ is an ordinary arc. We recall that $a_1,a_4,\tau$ and $a_2,a_3,\tau$ form two triangles in $T$ such that $a_1,a_3$ are clockwise to $\tau$ and $a_2,a_4$ are counterclockwise to $\tau$, as shown in the figure below. As $\mathcal O_0$ contains no punctures, $T_0$ has no self-folded triangle, $a_2\neq a_3$ and $a_1\neq a_4$. We similarly know $a_1\neq a_2$ and $a_3\neq a_4$ after do flip at $\tau$.

\centerline{\begin{tikzpicture}
\draw[-] (1,0) -- (2,0);
\draw[-] (1,-1) -- (2,-1);
\draw[-] (1,0) -- (2,-1);
\node[right] at (1.4,-0.45) {$\tau$};
\node[above] at (1.5,-0.05) {$a_2$};
\node[left] at (1,-0.5) {$a_1$};
\node[below] at (1.5,-0.95) {$a_4$};
\node[right] at (2,-0.5) {$a_3$};
%\node[right] at (1.2,-0.5) {$p_{i_1}$};
\draw [-] (1, 0) -- (1,-1);
\draw [-] (2, 0) -- (2,-1);
\draw [fill] (1,0) circle [radius=.05];
\draw [fill] (2,0) circle [radius=.05];
\draw [fill] (1,-1) circle [radius=.05];
\draw [fill] (2,-1) circle [radius=.05];
\end{tikzpicture}}

\subsubsection{When $a_1\neq a_3$ and $a_2\neq a_4$.}

By the construction of $\rho$ and $\tau'\neq \rho\notin T$, we have the following possibilities (the addition operation below is taken in $\mathbb Z_4$):

\begin{enumerate}[$(1)$]

  \item $\rho$ crosses $\tau,a_i$ for $i=1,2,3,4$.

  \item $\rho$ crosses $a_i,\tau, a_{i+1}$ for $i=1,3$.

  \item $\rho$ crosses $a_i,\tau,\tau',a_{i+2}$ for $i=1,2$.

  \item $\rho$ crosses $\tau',a_i$ for $i=1,2,3,4$.

  \item $\rho$ crosses $a_i,\tau',a_{i+1}$ for $i=2,4$.

\end{enumerate}

Because case (1) and case (4) are dual, case (2) and case (5) are dual, we shall only consider cases (1), (2) and (3).

\medskip

In case (1), $\rho$ and $\zeta$ have the same starting point. We may assume that $\rho$ crosses $\tau, a_2$. Then up to a difference of relative orientation, $G_{T,\rho}$ and $G_{T',\rho}$ are the following graphs, respectively,

\centerline{\begin{tikzpicture}
\draw[-] (1,0) -- (2,0);
\draw[-] (1,-1) -- (2,-1);
\draw[-] (2,0) -- (3,0);
\draw[-] (2,-1) -- (3,-1);
\draw[-] (3,0) -- (3,-1);
\draw[dashed] (1,0) -- (2,-1);
\draw[dashed] (2,0) -- (3,-1);
\node[right] at (1.4,-0.45) {$\tau$};
\node[above] at (1.5,-0.05) {$a_2$};
\node[left] at (1,-0.5) {$a_1$};
\node[below] at (1.5,-1) {$a_4$};
\node[right] at (1.75,-0.5) {$a_3$};
\node[right] at (2.4,-0.5) {$a_2$};
\node[above] at (2.5,-0.05) {$a_6$};
%\node[left] at (3,-0.5) {$a_5$};
\node[below] at (2.5,-1) {$\tau$};
\node[right] at (3,-0.5) {$a_5$};
%\node[right] at (1.2,-0.5) {$p_{i_1}$};
\draw [-] (1, 0) -- (1,-1);
\draw [-] (2, 0) -- (2,-1);
\draw [fill] (1,0) circle [radius=.05];
\draw [fill] (2,0) circle [radius=.05];
\draw [fill] (1,-1) circle [radius=.05];
\draw [fill] (2,-1) circle [radius=.05];
\draw [fill] (3,0) circle [radius=.05];
\draw [fill] (3,-1) circle [radius=.05];
\draw[-] (5,0) -- (6,0);
\draw[-] (5,-1) -- (6,-1);
\draw[dashed] (5,0) -- (6,-1);
\node[right] at (5.4,-0.5) {$a_2$};
\node[above] at (5.5,-0.05) {$a_5$};
\node[left] at (5,-0.5) {$a_1$};
\node[below] at (5.5,-0.95) {$\tau'$};
\node[right] at (6,-0.5) {$a_6$};
%\node[right] at (1.2,-0.5) {$p_{i_1}$};
\draw [-] (5, 0) -- (5,-1);
\draw [-] (6, 0) -- (6,-1);
\draw [fill] (5,0) circle [radius=.05];
\draw [fill] (6,0) circle [radius=.05];
\draw [fill] (5,-1) circle [radius=.05];
\draw [fill] (6,-1) circle [radius=.05];
\end{tikzpicture}}

In this case, $G_{T,\rho}$ (respectively $G_{T',\rho}$) has one $\tau$-(respectively $\tau'$-)equivalence class of type (III) (respectively (IV)). Thus $n^{\tau}(T,\rho)=n^{\tau'}(T',\rho)$. Moreover, we have a partition bijection from $\mathcal P(G_{T,\rho})$ to $\mathcal P(G_{T',\rho})$, as shown in the figure below.

\centerline{\includegraphics{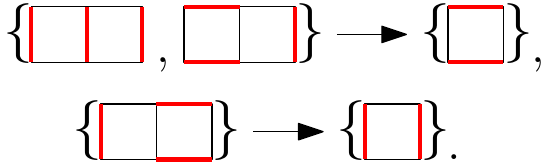}}

Under the partition bijection, it is easy to see the right upper edge labeled $a_5$/$a_6$ of $G_{T,\rho}$ is in the left set if and only if the right upper edge labeled $a_5$/$a_6$ of $G_{T',\rho}$ is in the right set. Thus, the last gluing edge of $G_{T,\rho}$ is in the left set if and only if the last gluing edge of $G_{T',\rho}$ in the right set. Since $\rho$ and $\zeta$ have the same starting point, $G_{T,\rho}$ and $G_{T',\rho}$ do not have first gluing edges.

\medskip

In case (2), we may assume that $\rho$ crosses $a_1, \tau, a_2$. Then up to a difference of relative orientation, $G_{T,\rho}$ and $G_{T',\rho}$ are the following graphs, respectively,

\centerline{\begin{tikzpicture}
\draw[-] (1,0) -- (2,0);
\draw[-] (1,-1) -- (2,-1);
\draw[-] (2,0) -- (3,0);
\draw[-] (2,-1) -- (3,-1);
\draw[-] (3,0) -- (3,-1);
\draw[-] (1,-1) -- (1,-2);
\draw[-] (2,-1) -- (2,-2);
\draw[-] (1,-2) -- (2,-2);
%\draw[-] (1,-1) -- (2,-1);
\draw[dashed] (1,0) -- (2,-1);
\draw[dashed] (2,0) -- (3,-1);
\draw[dashed] (1,-1) -- (2,-2);
\node[right] at (1.4,-0.45) {$\tau$};
\node[right] at (1.4,-1.5) {$a_1$};
\node[above] at (1.5,-0.05) {$a_2$};
\node[left] at (1.05,-0.5) {$a_1$};
\node[left] at (1.05,-1.5) {$a_8$};
\node[below] at (1.5,-0.8) {$a_4$};
\node[below] at (1.5,-1.95) {$a_7$};
\node[right] at (1.75,-0.5) {$a_3$};
\node[right] at (1.95,-1.5) {$\tau$};
\node[right] at (2.4,-0.5) {$a_2$};
\node[above] at (2.5,-0.05) {$a_6$};
%\node[left] at (3,-0.5) {$a_5$};
\node[below] at (2.5,-1) {$\tau$};
\node[right] at (2.95,-0.5) {$a_5$};
%\node[right] at (1.2,-0.5) {$p_{i_1}$};
\draw [-] (1, 0) -- (1,-1);
\draw [-] (2, 0) -- (2,-1);
\draw [fill] (1,0) circle [radius=.05];
\draw [fill] (2,0) circle [radius=.05];
\draw [fill] (1,-1) circle [radius=.05];
\draw [fill] (2,-1) circle [radius=.05];
\draw [fill] (1,-2) circle [radius=.05];
\draw [fill] (2,-2) circle [radius=.05];
\draw [fill] (3,0) circle [radius=.05];
\draw [fill] (3,-1) circle [radius=.05];
\draw[-] (5,0) -- (6,0);
\draw[-] (5,-1) -- (6,-1);
\draw[-] (5,-2) -- (6,-2);
\draw[-] (5,-1) -- (5,-2);
\draw[-] (6,-1) -- (6,-2);
\draw[dashed] (5,-1) -- (6,-2);
\draw[dashed] (5,0) -- (6,-1);
\node[right] at (5.4,-1.5) {$a_1$};
\node[right] at (5.4,-0.5) {$a_2$};
\node[above] at (5.5,-0.05) {$a_5$};
\node[left] at (5.05,-0.5) {$a_1$};
\node[below] at (5.5,-0.7) {$\tau'$};
\node[right] at (5.95,-0.5) {$a_6$};
\node[right] at (5.95,-1.5) {$a_2$};
\node[left] at (5.05,-1.5) {$a_8$};
%\node[below] at (1.5,-0.8) {$a_4$};
\node[below] at (5.5,-1.95) {$a_7$};
%\node[right] at (1.2,-0.5) {$p_{i_1}$};
\draw [-] (5, 0) -- (5,-1);
\draw [-] (6, 0) -- (6,-1);
\draw [fill] (5,0) circle [radius=.05];
\draw [fill] (6,0) circle [radius=.05];
\draw [fill] (5,-1) circle [radius=.05];
\draw [fill] (6,-1) circle [radius=.05];
\draw [fill] (5,-2) circle [radius=.05];
\draw [fill] (6,-2) circle [radius=.05];
\end{tikzpicture}}

In this case, $G_{T,\rho}$ (respectively $G_{T',\rho}$) has one $\tau$-(respectively $\tau'$-)equivalence class of type (II) (respectively (IV)). Thus $n^{\tau}(T,\rho)=n^{\tau'}(T',\rho)$. Moreover, we have a partition bijection from $\mathcal P(G_{T,\rho})$ to $\mathcal P(G_{T',\rho})$, as shown in the following.

\centerline{\includegraphics{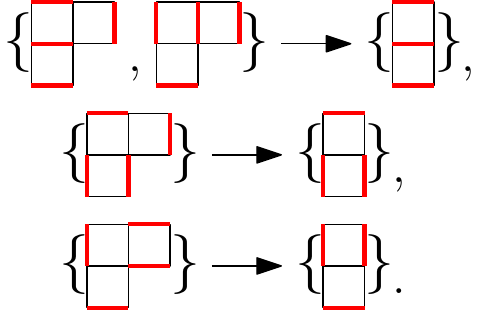}}

As case (1), it can be seen that the first/last gluing edge of $G_{T,\rho}$ is in the left set if and only if the first/last gluing edge of $G_{T',\rho}$ in the right set.

\medskip

In case (3), we may assume that $\rho$ crosses $a_2, \tau, \tau', a_4$. Then up to a difference of relative orientation, $G_{T,\rho}$ and $G_{T',\rho}$ are the following graphs, respectively,

\centerline{\begin{tikzpicture}
\draw[-] (0,0) -- (1,0);
\draw[-] (0,0) -- (0,-1);
\draw[-] (0,-1) -- (1,-1);
\draw[-] (1,0) -- (2,0);
\draw[-] (1,-1) -- (2,-1);
\draw[-] (2,0) -- (3,0);
\draw[-] (2,-1) -- (3,-1);
\draw[-] (3,0) -- (3,-1);
%\draw[-] (1,-1) -- (1,-2);
%\draw[-] (2,-1) -- (2,-2);
%\draw[-] (1,-2) -- (2,-2);
%\draw[-] (1,-1) -- (2,-1);
\draw[dashed] (1,0) -- (2,-1);
\draw[dashed] (0,0) -- (1,-1);
\draw[dashed] (2,0) -- (3,-1);
%\draw[dashed] (1,-1) -- (2,-2);
\node[right] at (1.4,-0.45) {$\tau$};
\node[above] at (1.5,0) {$a_2$};
\node[above] at (0.5,0) {$\tau$};
\node[left] at (1.35,-0.5) {$a_1$};
\node[below] at (1.5,-1) {$a_4$};
\node[below] at (0.5,-1) {$a_7$};
\node[right] at (1.75,-0.5) {$a_3$};
\node[right] at (2.4,-0.5) {$a_2$};
\node[left] at (0.8,-0.5) {$a_4$};
\node[above] at (2.5,0) {$a_6$};
\node[below] at (2.5,-1) {$\tau$};
\node[right] at (3,-0.5) {$a_5$};
\node[left] at (0,-0.5) {$a_8$};
\draw [-] (1, 0) -- (1,-1);
\draw [-] (2, 0) -- (2,-1);
\draw [fill] (0,0) circle [radius=.05];
\draw [fill] (0,-1) circle [radius=.05];
\draw [fill] (1,0) circle [radius=.05];
\draw [fill] (2,0) circle [radius=.05];
\draw [fill] (1,-1) circle [radius=.05];
\draw [fill] (2,-1) circle [radius=.05];
%\draw [fill] (1,-2) circle [radius=.05];
%\draw [fill] (2,-2) circle [radius=.05];
\draw [fill] (3,0) circle [radius=.05];
\draw [fill] (3,-1) circle [radius=.05];
\draw[-] (5,0) -- (6,0);
\draw[-] (5,-1) -- (6,-1);
\draw [-] (5, 0) -- (5,-1);
\draw [-] (6, 0) -- (6,-1);
\draw[-] (5,0) -- (5,1);
\draw[-] (5,-1) -- (5,-2);
\draw[-] (6,0) -- (6,1);
\draw[-] (6,-1) -- (6,-2);
\draw[-] (5,1) -- (6,1);
\draw[-] (5,-2) -- (6,-2);
\draw[dashed] (5,0) -- (6,-1);
\draw[dashed] (5,1) -- (6,0);
\draw[dashed] (5,-1) -- (6,-2);
\node[left] at (5.9,0.5) {$a_2$};
\node[above] at (5.5,-0.25) {$a_1$};
\node[right] at (6,-0.5) {$a_2$};
\node[right] at (6,0.5) {$a_5$};
\node[right] at (6,-1.5) {$\tau'$};
\node[right] at (5.3,-1.5) {$a_4$};
\node[below] at (5.5,-0.8) {$a_3$};
\node[left] at (5.9,-0.4) {$\tau'$};
\node[left] at (5,-0.5) {$a_4$};
\node[left] at (5,-1.5) {$a_8$};
\node[left] at (5,0.5) {$\tau'$};
\node[below] at (5.5,-1.95) {$a_7$};
\node[above] at (5.5,0.95) {$a_6$};
\draw [fill] (5,0) circle [radius=.05];
\draw [fill] (6,0) circle [radius=.05];
\draw [fill] (5,-1) circle [radius=.05];
\draw [fill] (6,-1) circle [radius=.05];
\draw [fill] (6,-2) circle [radius=.05];
\draw [fill] (5,-2) circle [radius=.05];
\draw [fill] (5,1) circle [radius=.05];
\draw [fill] (6,1) circle [radius=.05];
\end{tikzpicture}}

In this case, $G_{T,\rho}$ (respectively $G_{T',\rho}$) has one $\tau$-(respectively $\tau'$-)equivalence class of type (I). Thus $n^{\tau}(T,\rho)=n^{\tau'}(T',\rho)$. Moreover, we have a partition bijection from $\mathcal P(G_{T,\rho})$ to $\mathcal P(G_{T',\rho})$, as shown in the following figure.

\centerline{\includegraphics{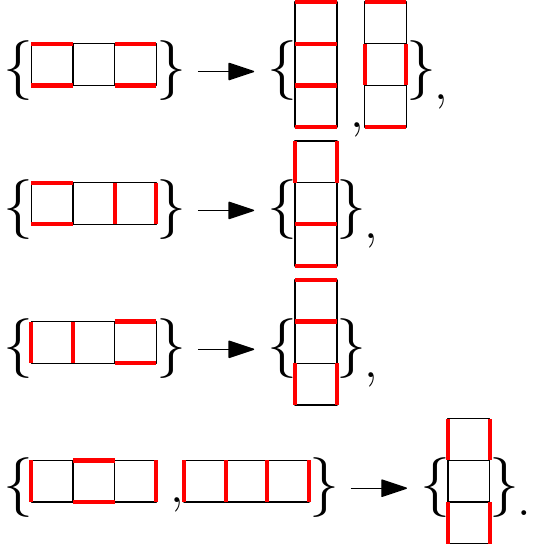}}

As cases (1) and (2), it can be seen that the first/last gluing edge of $G_{T,\rho}$ is in the left set if and only if the first/last gluing edge of $G_{T',\rho}$ in the right set.

\subsubsection{When $a_1= a_3$ or $a_2= a_4$.} By Lemma 4.2 in \cite{H}, $a_1=a_3$ and $a_2=a_4$ can not hold at the same time. We may assume $a_1=a_3$. Thus $a_2\neq a_4$.

\medskip

Consider the universal covering of $\mathcal O_0$, we have the following possibilities (we illustrate cases (1)-(8) in the figure below, the other cases are dual):

\medskip

\centerline{\includegraphics{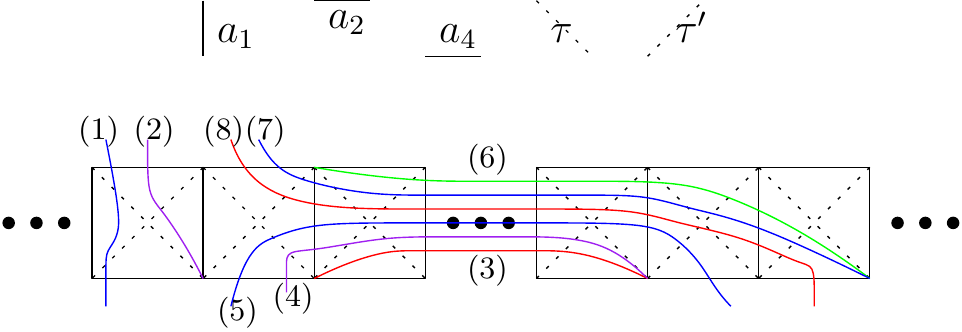}}

\begin{enumerate}[$(1)$]

  \item $\rho$ crosses $a_2,\tau',\tau,a_4$ in order.

  \item $\rho$ crosses $a_i,\tau$ in order for $i=2,4$.

  \item $\rho$ crosses $(\tau,\overset{s}{\overbrace{a_1,\tau'),\cdots,(\tau,a_1}},\tau')$ in order for $s\geq 1$.

  \item $\rho$ crosses $a_i,(\tau,\overset{s}{\overbrace{a_1,\tau'),\cdots,(\tau,a_1}},\tau')$ in order for $i=2,4$ and $s\geq 1$.

  \item $\rho$ crosses $a_i,(\tau,\overset{s}{\overbrace{a_1,\tau'),\cdots,(\tau,a_1}},\tau'),a_i$ in order for $i=2,4$ and $s\geq 1$.

  \item $\rho$ crosses $\tau',a_1,(\tau',\overset{s}{\overbrace{\tau,a_1),\cdots,(\tau',\tau}},a_1),\tau'$ in order for $s\geq 0$.

  \item $\rho$ crosses $a_i,\tau',a_1,(\tau',\overset{s}{\overbrace{\tau,a_1),\cdots,(\tau',\tau}},a_1),\tau'$ in order for $i=2,4$ and $s\geq 0$.

  \item $\rho$ crosses $a_2,\tau',a_1,(\tau',\overset{s}{\overbrace{\tau,a_1),\cdots,(\tau',\tau}},a_1),\tau',a_4$ in order for $s\geq 0$.

  \item $\rho$ crosses $a_i,\tau'$ in order for $i=2,4$.

  \item $\rho$ crosses $\tau,a_1,(\tau,\overset{s}{\overbrace{\tau',a_1),\cdots,(\tau,\tau'}},a_1),\tau$ in order for $s\geq 0$.

  \item $\rho$ crosses $a_i,\tau,a_1,(\tau,\overset{s}{\overbrace{\tau',a_1),\cdots,(\tau,\tau'}},a_1),\tau$ in order for $i=2,4$ and $s\geq 0$.

  \item $\rho$ crosses $a_2,\tau,a_1,(\tau,\overset{s}{\overbrace{\tau',a_1),\cdots,(\tau,\tau'}},a_1),\tau,a_4$ in order for $s\geq 0$.

\end{enumerate}

\medskip

Since cases (2) and (9), (6) and (10), (7) and (11), (8) and (12) are dual respectively, we shall only discuss cases (1)-(8). It should be noted that $\rho$ crosses itself in cases (3,4,5).

\medskip

We omit the discussion of case (1) and case (2) because they are similar to the case (3) and case (1), respectively, when $a_1\neq a_3$ and $a_2\neq a_4$.

\medskip

For the remaining cases, we should introduce more notation. If both $a$ and $b$ are edges of two different triangles in $T$, for any non-negative integer $s$, we define the graphs $G_s(a,b)$ and $H_s(a,b)$ as follows. See the following two graphs.

\begin{enumerate}[$(a)$]

  \item $G_s(a,b)$ is isomorphic to $G_{T,\varsigma}$, where $\varsigma$ is the curve in $\mathcal O$ crossing $a$ $s$ times and $b$ $s$ times alternately.

  \item $H_s(a,b)$ is isomorphic to $G_{T,\varsigma}$, where $\varsigma$ is the curve in $\mathcal O$ crossing $a$ $s+1$ times and $b$ $s$ times alternately.

\end{enumerate}

\medskip

\centerline{\begin{tikzpicture}
\draw[-] (0,0) -- (1,0);
\draw[-] (0,0) -- (0,-1);
\draw[-] (0,-1) -- (1,-1);
\draw[-] (1,0) -- (2,0);
\draw[-] (1,-1) -- (2,-1);
\draw[-] (2,0) -- (3,0);
\draw[-] (2,-1) -- (3,-1);
\draw[-] (3,0) -- (3,-1);
\draw[-] (3,0) -- (4,0);
\draw[-] (3,-1) -- (4,-1);
\draw[-] (4,0) -- (4,-1);
\draw[-] (5,0) -- (6,0);
\draw[-] (5,0) -- (5,-1);
\draw[-] (6,0) -- (6,-1);
\draw[-] (5,-1) -- (6,-1);
\draw[-] (6,0) -- (7,0);
\draw[-] (6,-1) -- (7,-1);
\draw[-] (7,0) -- (7,-1);
\draw[dashed] (1,0) -- (2,-1);
\draw[dashed] (0,0) -- (1,-1);
\draw[dashed] (2,0) -- (3,-1);
\draw[dashed] (3,0) -- (4,-1);
\draw[dashed] (5,0) -- (6,-1);
\draw[dashed] (6,0) -- (7,-1);
\node[above] at (0.5,-0.05) {$b$};
\node[above] at (1.5,-0.05) {$a$};
\node[above] at (2.5,-0.05) {$b$};
\node[above] at (3.5,-0.05) {$a$};
\node[above] at (5.5,-0.05) {$b$};
\node[above] at (6.5,-0.05) {$a$};
\node[below] at (0.5,-0.95) {$b$};
\node[below] at (1.5,-0.95) {$a$};
\node[below] at (2.5,-0.95) {$b$};
\node[below] at (3.5,-0.95) {$a$};
\node[below] at (5.5,-0.95) {$b$};
\node[below] at (6.5,-0.95) {$a$};
\node[right] at (1.4,-0.45) {$b$};
\node[right] at (3.4,-0.45) {$b$};
\node[right] at (2.4,-0.5) {$a$};
\node[right] at (6.4,-0.45) {$b$};
\node[right] at (5.4,-0.5) {$a$};
\node[left] at (0.8,-0.5) {$a$};
\draw [-] (1, 0) -- (1,-1);
\draw [-] (2, 0) -- (2,-1);
\draw [fill] (0,0) circle [radius=.05];
\draw [fill] (0,-1) circle [radius=.05];
\draw [fill] (1,0) circle [radius=.05];
\draw [fill] (2,0) circle [radius=.05];
\draw [fill] (1,-1) circle [radius=.05];
\draw [fill] (2,-1) circle [radius=.05];
\draw [fill] (3,0) circle [radius=.05];
\draw [fill] (3,-1) circle [radius=.05];
\draw [fill] (4,0) circle [radius=.05];
\draw [fill] (4,-1) circle [radius=.05];
\draw [fill] (5,0) circle [radius=.05];
\draw [fill] (5,-1) circle [radius=.05];
\draw [fill] (6,0) circle [radius=.05];
\draw [fill] (6,-1) circle [radius=.05];
\draw [fill] (7,0) circle [radius=.05];
\draw [fill] (4.3,-0.5) circle [radius=.02];
\draw [fill] (4.5,-0.5) circle [radius=.02];
\draw [fill] (4.7,-0.5) circle [radius=.02];
\draw [fill] (7,-1) circle [radius=.05];
\end{tikzpicture}}
%
%\medskip

\centerline{\begin{tikzpicture}
\draw[-] (0,0) -- (1,0);
\draw[-] (0,0) -- (0,-1);
\draw[-] (0,-1) -- (1,-1);
\draw[-] (1,0) -- (2,0);
\draw[-] (1,-1) -- (2,-1);
\draw[-] (2,0) -- (3,0);
\draw[-] (2,-1) -- (3,-1);
\draw[-] (3,0) -- (3,-1);
\draw[-] (3,0) -- (4,0);
\draw[-] (3,-1) -- (4,-1);
\draw[-] (4,0) -- (4,-1);
\draw[-] (5,0) -- (6,0);
\draw[-] (5,0) -- (5,-1);
\draw[-] (6,0) -- (6,-1);
\draw[-] (5,-1) -- (6,-1);
\draw[-] (6,0) -- (7,0);
\draw[-] (6,-1) -- (7,-1);
\draw[-] (7,0) -- (7,-1);
\draw[-] (7,0) -- (8,0);
\draw[-] (7,-1) -- (8,-1);
\draw[-] (8,0) -- (8,-1);
\draw[dashed] (1,0) -- (2,-1);
\draw[dashed] (0,0) -- (1,-1);
\draw[dashed] (2,0) -- (3,-1);
\draw[dashed] (3,0) -- (4,-1);
\draw[dashed] (5,0) -- (6,-1);
\draw[dashed] (6,0) -- (7,-1);
\draw[dashed] (7,0) -- (8,-1);
\node[above] at (0.5,-0.05) {$b$};
\node[above] at (1.5,-0.05) {$a$};
\node[above] at (2.5,-0.05) {$b$};
\node[above] at (3.5,-0.05) {$a$};
\node[above] at (5.5,-0.05) {$b$};
\node[above] at (6.5,-0.05) {$a$};
\node[above] at (7.5,-0.05) {$b$};
\node[below] at (0.5,-0.95) {$b$};
\node[below] at (1.5,-0.95) {$a$};
\node[below] at (2.5,-0.95) {$b$};
\node[below] at (3.5,-0.95) {$a$};
\node[below] at (5.5,-0.95) {$b$};
\node[below] at (6.5,-0.95) {$a$};
\node[below] at (7.5,-0.95) {$b$};
\node[right] at (1.4,-0.45) {$b$};
\node[right] at (3.4,-0.45) {$b$};
\node[right] at (2.4,-0.5) {$a$};
\node[right] at (6.4,-0.45) {$b$};
\node[right] at (5.4,-0.5) {$a$};
\node[right] at (7.4,-0.5) {$a$};
\node[left] at (0.8,-0.5) {$a$};
\draw [-] (1, 0) -- (1,-1);
\draw [-] (2, 0) -- (2,-1);
\draw [fill] (0,0) circle [radius=.05];
\draw [fill] (0,-1) circle [radius=.05];
\draw [fill] (1,0) circle [radius=.05];
\draw [fill] (2,0) circle [radius=.05];
\draw [fill] (1,-1) circle [radius=.05];
\draw [fill] (2,-1) circle [radius=.05];
\draw [fill] (3,0) circle [radius=.05];
\draw [fill] (3,-1) circle [radius=.05];
\draw [fill] (4,0) circle [radius=.05];
\draw [fill] (4,-1) circle [radius=.05];
\draw [fill] (5,0) circle [radius=.05];
\draw [fill] (5,-1) circle [radius=.05];
\draw [fill] (6,0) circle [radius=.05];
\draw [fill] (6,-1) circle [radius=.05];
\draw [fill] (7,0) circle [radius=.05];
\draw [fill] (8,0) circle [radius=.05];
\draw [fill] (8,-1) circle [radius=.05];
\draw [fill] (4.3,-0.5) circle [radius=.02];
\draw [fill] (4.5,-0.5) circle [radius=.02];
\draw [fill] (4.7,-0.5) circle [radius=.02];
\draw [fill] (7,-1) circle [radius=.05];
\end{tikzpicture}}

We denote the $i$-th tile of $G_s(a,b)$ by $G_i$. For a given sequence $(\lambda_1,\cdots,\lambda_s)\in \{0,1\}^{s}$, let $\mathcal G_{(\lambda_1,\cdots,\lambda_s)}(a,b)$ be the set of perfect matching of $G_s(a,b)$ containing perfect matchings $P$ which satisfy the following conditions: if $\lambda_i=1$, then the two edges of $G_{2i}$ labeled $a$ are both in $P$; if $\lambda_i=0$, then the two edges of $G_{2i}$ labeled $a$ are both not in $P$. Let $\mathcal G'_{(\lambda_1,\cdots,\lambda_s)}(a,b)$ be the set of perfect matching of $G_s(a,b)$ containing perfect matchings $P$ which satisfy the following conditions: if $\lambda_i=1$, then the two edges of $G_{2i-1}$ labeled $b$ are both in $P$; if $\lambda_i=0$, then the two edges of $G_{2i-1}$ labeled $b$ are both not in $P$. See the figure below for example $\mathcal G_{(1,0)}(a,b)$ and $\mathcal G'_{(1,0)}(a,b)$.

\centerline{\includegraphics{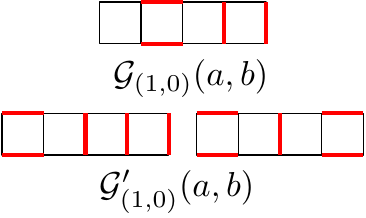}}

We denote the $i$-th tile of $H_s(a,b)$ by $H_i$. For a sequence $(\lambda_1,\cdots,\lambda_s)\in \{0,1\}^{s}$, let $\mathcal H_{(\lambda_1,\cdots,\lambda_s)}(a,b)$ be the set of perfect matching of $H_s(a,b)$ containing all perfect matchings $P$ which satisfy the following conditions: if $\lambda_i=1$, then the two edges of $H_{2i}$ labeled $a$ are both in $P$; if $\lambda_i=0$, then the two edges of $H_{2i}$ labeled $a$ are both not in $P$; for a sequence $(\lambda_1,\cdots,\lambda_{s+1})\in \{0,1\}^{s+1}$. Let $\mathcal H'_{(\lambda_1,\cdots,\lambda_{s+1})}(a,b)$ be the set of perfect matching of $H_s(a,b)$ containing all the perfect matching $P$ which satisfies the following conditions: if $\lambda_i=1$, then the two edges of $H_{2i-1}$ labeled $b$ are both in $P$; if $\lambda_i=0$, then the two edges of $H_{2i-1}$ labeled $b$ are both not in $P$. See the following figure for example $\mathcal H_{(1,0)}(a,b)$ and $\mathcal H'_{(1,0,1)}(a,b)$.

\centerline{\includegraphics{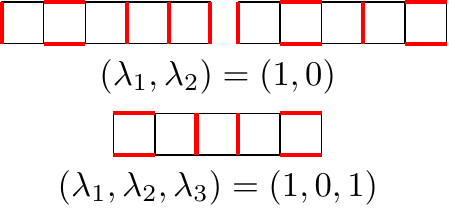}}

The following observations are easy but important.

\medskip

\begin{Lemma}\label{basicdecom}

\begin{enumerate}[$(1)$]

  \item $\mathcal G_{(\lambda_1,\cdots,\lambda_s)}(a,b) = \mathcal P^{a}_{(\lambda_1-1,\lambda_1+\lambda_2-1,\cdots,\lambda_{s-1}+\lambda_s-1,\lambda_s)}(G_s(a,b))$.

  \item $\mathcal G'_{(\lambda_1,\cdots,\lambda_s)}(a,b) = \mathcal P^{b}_{(\lambda_1,\lambda_1+\lambda_2-1,\lambda_2+\lambda_3-1,\cdots,\lambda_{s-1}+\lambda_{s}-1,\lambda_s-1)}(G_s(a,b))$.

  \item \[\begin{array}{ccl} \mathcal P(G_s(a,b))
        & = & \bigsqcup_{(\lambda_1,\cdots,\lambda_s)\in \{0,1\}^{s}}\mathcal G_{(\lambda_1,\cdots,\lambda_s)}(a,b)\\
        & = & \bigsqcup_{(\lambda_1,\cdots,\lambda_s)\in \{0,1\}^{s}}\mathcal G'_{(\lambda_1,\cdots,\lambda_s)}(a,b).
        \end{array}\]

  \item $\mathcal H_{(\lambda_1,\cdots,\lambda_s)}(a,b) = \mathcal P^{a}_{(\lambda_1-1,\lambda_1+\lambda_2-1,\cdots,\lambda_{s-1}+\lambda_s-1,\lambda_s-1)}(H_s(a,b))$.

  \item $\mathcal H'_{(\lambda_1,\cdots,\lambda_{s+1})}(a,b) = \mathcal P^{b}_{(\lambda_1,\lambda_1+\lambda_2-1,\cdots,\lambda_{s}+\lambda_{s+1}-1,\lambda_{s+1})}(H_s(a,b))$.

  \item  \[\begin{array}{ccl} \mathcal P(H_s(a,b))
  & = & \bigsqcup_{(\lambda_1,\cdots,\lambda_s)\in \{0,1\}^{s}}\mathcal H_{(\lambda_1,\cdots,\lambda_s)}(a,b)\\
  & = & \bigsqcup_{(\lambda_1,\cdots,\lambda_{s+1})\in \{0,1\}^{s+1}}\mathcal H'_{(\lambda_1,\cdots,\lambda_{s+1})}(a,b).
  \end{array}\]

\end{enumerate}

\end{Lemma}

\begin{proof}

(1,2,4,5) follow by definition. (3,6) follow immediately by Lemma \ref{inone}.
%We shall only prove (1). It is equivalent to say, for each $i$ satisfies $1\leq i\leq s$ and a perfect matching $P$, the edges of $G_i$ labeled $a$ are either in $P$ or not in $P$ at the same time. Otherwise, choose the maximal $i$ such that one edge of $G_i$ labeled $a$ in $P$ and the other one not in $P$ for some $P$. Clearly $i\neq s$. One can see immediately from the definition of perfect matching, one edge of $G_{i+1}$ labeled $a$ in $P$ and the other one not in $P$ for some $P$. This contracts to the maximality of $i$.
\end{proof}

\medskip

\begin{Lemma}\label{diff}
\begin{enumerate}[$(1)$]

  \item If $P\neq Q\in \mathcal G_{(\lambda_1,\cdots,\lambda_s)}(a,b)$, then there exists $i\in [1,s]$ such that $P$ and $Q$ can twist on the tile $H_{2i-1}$ and $P\cap E(H_{2i-1})\neq Q\cap E(H_{2i-1})$.

  \item If $P'\neq Q'\in \mathcal G'_{(\lambda_1,\cdots,\lambda_{s})}(a,b)$, then there exists $i\in [1,s]$ such that $P'$ and $Q'$ can twist on the tile $H_{2i}$ and $P'\cap E(H_{2i})\neq Q'\cap E(H_{2i})$.

  \item If $P\neq Q\in \mathcal H_{(\lambda_1,\cdots,\lambda_s)}(a,b)$, then there exists $i\in [1,s+1]$ such that $P$ and $Q$ can twist on the tile $H_{2i-1}$ and $P\cap E(H_{2i-1})\neq Q\cap E(H_{2i-1})$.

  \item If $P'\neq Q'\in \mathcal H'_{(\lambda_1,\cdots,\lambda_{s+1})}(a,b)$, then there exists $i\in [1,s]$ such that $P'$ and $Q'$ can twist on the tile $H_{2i}$ and $P'\cap E(H_{2i})\neq Q'\cap E(H_{2i})$.

\end{enumerate}

\end{Lemma}

\begin{proof}

We shall only prove (1) because the remainder can be proved similarly. If we assume $\lambda_0=0$ in convention, according to the definition of $\mathcal G_{(\lambda_1,\cdots,\lambda_s)}(a,b)$, $P$ and $Q$ can twist on tiles $H_{2i-1}$ for $i\in [1,s]$ with $\lambda_{i-1}=\lambda_i=0$, the other edges in $P$ and $Q$ are the same. Therefore, the result follows.
\end{proof}

\medskip

We now turn to the comparison of the perfect matchings of $G_{T,\rho}$ and $G_{T',\rho}$.

\medskip

In case (3), $\rho$ and $\zeta$ have the endpoints. Thus $\rho=\zeta$. Up to a difference of relative orientation, $G_{T,\rho}$ and $G_{T',\rho}$ are the following graphs, respectively,

\centerline{\begin{tikzpicture}
\draw[-] (0,0) -- (1,0);
\draw[-] (0,0) -- (0,-1);
\draw[-] (0,-1) -- (1,-1);
\draw[-] (1,0) -- (2,0);
\draw[-] (1,-1) -- (2,-1);
\draw[-] (2,0) -- (3,0);
\draw[-] (2,-1) -- (3,-1);
\draw[-] (3,0) -- (3,-1);
\draw[-] (3,0) -- (4,0);
\draw[-] (3,-1) -- (4,-1);
\draw[-] (4,0) -- (4,-1);
\draw[-] (5,0) -- (6,0);
\draw[-] (5,0) -- (5,-1);
\draw[-] (6,0) -- (6,-1);
\draw[-] (5,-1) -- (6,-1);
\draw[-] (6,0) -- (7,0);
\draw[-] (6,-1) -- (7,-1);
\draw[-] (7,0) -- (7,-1);
\draw[dashed] (1,0) -- (2,-1);
\draw[dashed] (0,0) -- (1,-1);
\draw[dashed] (2,0) -- (3,-1);
\draw[dashed] (3,0) -- (4,-1);
\draw[dashed] (5,0) -- (6,-1);
\draw[dashed] (6,0) -- (7,-1);
\node[left] at (0.05,-0.5) {$a_4$};
\node[left] at (2.25,-0.5) {$a_4$};
\node[left] at (4.25,-0.5) {$a_4$};
\node[left] at (5.25,-0.5) {$a_4$};
\node[right] at (6.95,-0.5) {$a_4$};
\node[left] at (1.25,-0.5) {$a_2$};
\node[left] at (3.25,-0.5) {$a_2$};
\node[left] at (6.25,-0.5) {$a_2$};
\node[above] at (0.5,-0.05) {$a_1$};
\node[above] at (1.5,-0.05) {$\tau$};
\node[above] at (2.5,-0.05) {$a_1$};
\node[above] at (3.5,-0.05) {$\tau$};
\node[above] at (5.5,-0.05) {$a_1$};
\node[above] at (6.5,-0.05) {$\tau$};
\node[below] at (0.5,-0.95) {$a_1$};
\node[below] at (1.5,-0.95) {$\tau$};
\node[below] at (2.5,-0.95) {$a_1$};
\node[below] at (3.5,-0.95) {$\tau$};
\node[below] at (5.5,-0.95) {$a_1$};
\node[below] at (6.5,-0.95) {$\tau$};
\node[right] at (1.3,-0.45) {$a_1$};
\node[right] at (3.3,-0.45) {$a_1$};
\node[right] at (2.3,-0.5) {$\tau$};
\node[right] at (6.3,-0.45) {$a_1$};
\node[right] at (5.3,-0.5) {$\tau$};
\node[left] at (0.8,-0.5) {$\tau$};
\draw [-] (1, 0) -- (1,-1);
\draw [-] (2, 0) -- (2,-1);
\draw [fill] (0,0) circle [radius=.05];
\draw [fill] (0,-1) circle [radius=.05];
\draw [fill] (1,0) circle [radius=.05];
\draw [fill] (2,0) circle [radius=.05];
\draw [fill] (1,-1) circle [radius=.05];
\draw [fill] (2,-1) circle [radius=.05];
\draw [fill] (3,0) circle [radius=.05];
\draw [fill] (3,-1) circle [radius=.05];
\draw [fill] (4,0) circle [radius=.05];
\draw [fill] (4,-1) circle [radius=.05];
\draw [fill] (5,0) circle [radius=.05];
\draw [fill] (5,-1) circle [radius=.05];
\draw [fill] (6,0) circle [radius=.05];
\draw [fill] (6,-1) circle [radius=.05];
\draw [fill] (7,0) circle [radius=.05];
\draw [fill] (4.3,-0.5) circle [radius=.02];
\draw [fill] (4.5,-0.5) circle [radius=.02];
\draw [fill] (4.7,-0.5) circle [radius=.02];
\draw [fill] (7,-1) circle [radius=.05];
\end{tikzpicture}}

\centerline{\begin{tikzpicture}
\draw[-] (0,0) -- (1,0);
\draw[-] (0,0) -- (0,-1);
\draw[-] (0,-1) -- (1,-1);
\draw[-] (1,0) -- (2,0);
\draw[-] (1,-1) -- (2,-1);
\draw[-] (2,0) -- (3,0);
\draw[-] (2,-1) -- (3,-1);
\draw[-] (3,0) -- (3,-1);
\draw[-] (3,0) -- (4,0);
\draw[-] (3,-1) -- (4,-1);
\draw[-] (4,0) -- (4,-1);
\draw[-] (5,0) -- (6,0);
\draw[-] (5,0) -- (5,-1);
\draw[-] (6,0) -- (6,-1);
\draw[-] (5,-1) -- (6,-1);
\draw[-] (6,0) -- (7,0);
\draw[-] (6,-1) -- (7,-1);
\draw[-] (7,0) -- (7,-1);
\draw[dashed] (1,0) -- (2,-1);
\draw[dashed] (0,0) -- (1,-1);
\draw[dashed] (2,0) -- (3,-1);
\draw[dashed] (3,0) -- (4,-1);
\draw[dashed] (5,0) -- (6,-1);
\draw[dashed] (6,0) -- (7,-1);
\node[left] at (2.25,-0.5) {$a_4$};
\node[left] at (4.25,-0.5) {$a_4$};
\node[left] at (5.25,-0.5) {$a_4$};
\node[right] at (6.95,-0.5) {$a_4$};
\node[left] at (1.25,-0.5) {$a_2$};
\node[left] at (0.05,-0.5) {$a_4$};
\node[left] at (3.25,-0.5) {$a_2$};
\node[left] at (6.25,-0.5) {$a_2$};
\node[above] at (0.5,-0.05) {$\tau'$};
\node[above] at (1.5,-0.05) {$a_1$};
\node[above] at (2.5,-0.05) {$\tau'$};
\node[above] at (3.5,-0.05) {$a_1$};
\node[above] at (5.5,-0.05) {$\tau'$};
\node[above] at (6.5,-0.05) {$a_1$};
\node[below] at (0.5,-0.95) {$\tau'$};
\node[below] at (1.5,-0.95) {$a_1$};
\node[below] at (2.5,-0.95) {$\tau'$};
\node[below] at (3.5,-0.95) {$a_1$};
\node[below] at (5.5,-0.95) {$\tau'$};
\node[below] at (6.5,-0.95) {$a_1$};
\node[right] at (1.3,-0.45) {$\tau'$};
\node[right] at (3.3,-0.45) {$\tau'$};
\node[right] at (2.3,-0.5) {$a_1$};
\node[right] at (6.3,-0.45) {$\tau'$};
\node[right] at (5.3,-0.5) {$a_1$};
\node[left] at (0.8,-0.5) {$a_1$};
\draw [-] (1, 0) -- (1,-1);
\draw [-] (2, 0) -- (2,-1);
\draw [fill] (0,0) circle [radius=.05];
\draw [fill] (0,-1) circle [radius=.05];
\draw [fill] (1,0) circle [radius=.05];
\draw [fill] (2,0) circle [radius=.05];
\draw [fill] (1,-1) circle [radius=.05];
\draw [fill] (2,-1) circle [radius=.05];
\draw [fill] (3,0) circle [radius=.05];
\draw [fill] (3,-1) circle [radius=.05];
\draw [fill] (4,0) circle [radius=.05];
\draw [fill] (4,-1) circle [radius=.05];
\draw [fill] (5,0) circle [radius=.05];
\draw [fill] (5,-1) circle [radius=.05];
\draw [fill] (6,0) circle [radius=.05];
\draw [fill] (6,-1) circle [radius=.05];
\draw [fill] (7,0) circle [radius=.05];
\draw [fill] (4.3,-0.5) circle [radius=.02];
\draw [fill] (4.5,-0.5) circle [radius=.02];
\draw [fill] (4.7,-0.5) circle [radius=.02];
\draw [fill] (7,-1) circle [radius=.05];
\end{tikzpicture}}

We clearly have $n^{\tau}(T,\rho)=n^{\tau'}(T',\rho)$.

\medskip

Since $G_{T,\rho}$ and $G_{T',\rho}$ are isomorphic to $G_{s}(\tau,a_1)$ and $G_{s}(a_1,\tau')$ respectively, by Lemma \ref{basicdecom}, we have $\mathcal P(G_{T,\rho})\cong\bigsqcup_{(\lambda_1,\cdots,\lambda_{s})\in \{0,1\}^{s}}\mathcal G_{(\lambda_1,\cdots,\lambda_{s})}(\tau,a_1)$, and  $\mathcal P(G_{T',\rho})\cong\bigsqcup_{(\lambda_1,\cdots,\lambda_{s})\in \{0,1\}^{s}}\mathcal G'_{(\lambda_1,\cdots,\lambda_{s})}(a_1,\tau')$ as sets .

\medskip

We have the following observation.

\medskip

\begin{Lemma}\label{non-tau-mu0}

Let $P\in \mathcal G_{(\lambda_1,\cdots,\lambda_{s})}(\tau,a_1)$ for some sequence $(\lambda_1,\cdots,\lambda_{s})$. We assume $\lambda_{0}=0$ in convention. Then for any $i\in [1,s]$,

 \begin{enumerate}[$(1)$]

   \item the edge labeled $a_2$ of the $(2i-1)$-th tile is in $P$ but non-$\tau$-mutable if and only if $\lambda_{i-1}=1$ and $\lambda_i=0$;

   \item the edge labeled $a_4$ of the $(2i-1)$-th tile is in $P$ but non-$\tau$-mutable if and only if $\lambda_{i-1}=0$ and $\lambda_{i}=1$;

   \item the edge labeled $a_4$ of the $2s$-th tile is in $P$ but non-$\tau$-mutable if and only if $\lambda_{s}=0$.

 \end{enumerate}

\end{Lemma}

\begin{proof}

We shall only prove (1) because (2) and (3) can be proved similarly.

``Only If Part:"
Since the edge labeled $a_2$ of the $(2i-1)$-th tile is in $P$ but non-$\tau$-mutable, the edges labeled $a_1,a_4$ of the $(2i-1)$-th tile and the edges labeled $\tau$ of the $2i$-th tile are not in $P$. Thus the edges labeled $\tau$ of the $(2i-2)$-th tile are in $P$. Therefore, $\lambda_{i-1}=1$ and $\lambda_i=0$.

``If Part:" Since $\lambda_{i-1}=1$, the edges labeled $\tau$ of the $(2i-2)$-th tile are in $P$, and hence the edges labeled $a_1$ and $a_4$ of the $(2i-1)$-th tile are not in $P$. Since $\lambda_i=0$, the edges labeled $\tau$ of the $2i$-th tile are not in $P$. Therefore, the edge labeled $a_2$ of the $(2i-1)$-th tile is in $P$ but non-$\tau$-mutable.
\end{proof}

\medskip

The following result can be similarly verified.

\medskip

\begin{Lemma}\label{non-tau-mu}

Let $P\in \mathcal G'_{(\lambda_1,\cdots,\lambda_{s})}(a_1,\tau')$ for some sequence $(\lambda_1,\cdots,\lambda_{s})$. We assume $\lambda_{s+1}=0$ in convention. Then for any $i\in [1,s]$,

 \begin{enumerate}[$(1)$]

   \item the edge labeled $a_2$ of the $(2i-1)$-th tile is in $P$ but non-$\tau'$-mutable if and only if $\lambda_{i}=0$ and $\lambda_{i+1}=1$;

   \item the edge labeled $a_4$ of the $2i$-th tile is in $P$ but non-$\tau'$-mutable if and only if $\lambda_{i}=1$ and $\lambda_{i+1}=0$;

   \item the edge labeled $a_4$ of the $1$-th tile is in $P$ but non-$\tau'$-mutable if and only if $\lambda_{1}=1$.

 \end{enumerate}

\end{Lemma}

\medskip

\begin{Lemma}\label{same-num1}

Let $P\in \mathcal G_{(\lambda_1,\cdots,\lambda_{s})}(\tau,a_1)$ and $P'\in\mathcal G'_{(1-\lambda_1,\cdots,1-\lambda_{s})}(a_1,\tau')$. For each arc $\tau \neq a\in T$, the number of the non-$\tau$-mutable edges labeled $a$ in $P$ equals to the number of the non-$\tau'$-mutable edges labeled $a$ in $P'$.

\end{Lemma}

\begin{proof}

By Lemma \ref{inone}, the two numbers both are equal to $0$ if $a\neq a_2,a_4$. If $a=a_2$, by Lemma \ref{non-tau-mu0} and Lemma \ref{non-tau-mu}, for each $i\in [1,s]$, the edge labeled $a_2$ of the $(2i-1)$-th tile of $G_{T,\rho}$ is non-$\tau$-mutable in $P$ if and only if the edge labeled $a_2$ of the $(2i-1)$-th tile of $G_{T',\rho}$ is non-$\tau$-mutable in $P'$. Thus, the two numbers are equal in this case. Similarly, the result holds for $a=a_4$.
\end{proof}

\medskip

$\psi_{\rho}: \mathcal G_{(\lambda_1,\cdots,\lambda_{s})}(\tau,a_1)\rightarrow \mathcal G'_{(1-\lambda_1,\cdots,1-\lambda_{s})}(a_1,\tau')$ clearly gives a partition bijection between $\mathcal P(G_{T,\rho})$ and $\mathcal P(G_{T',\rho})$. Moreover, for any $P\in \mathcal G_{(\lambda_1,\cdots,\lambda_{s})}(\tau,a_1)$ and $P'\in \mathcal G'_{(1-\lambda_1,\cdots,1-\lambda_{s})}(a_1,\tau')$, by Lemma \ref{same-num1} the number of non-$\tau$-mutable edges labeled $a$ in $P$ equals to the number of non-$\tau'$-mutable edges labeled $a$ in $P'$ for any $\tau \neq a\in T$.

\medskip

In this case, as $\rho=\zeta$, $G_{T,\rho}$ and $G_{T',\rho}$ have no first or last gluing edge.

\medskip

In case (4), $\rho$ and $\zeta$ have the same endpoint. We may assume $i=4$. Then up to a difference of relative orientation, $G_{T,\rho}$ and $G_{T',\rho}$ are the following graphs, respectively,

\centerline{\begin{tikzpicture}
\draw[-] (0,-2) -- (1,-2);
\draw[-] (0,-1) -- (0,-2);
\draw[-] (1,-1) -- (1,-2);
\draw[-] (0,0) -- (1,0);
\draw[-] (0,0) -- (0,-1);
\draw[-] (0,-1) -- (1,-1);
\draw[-] (1,0) -- (2,0);
\draw[-] (1,-1) -- (2,-1);
\draw[-] (2,0) -- (3,0);
\draw[-] (2,-1) -- (3,-1);
\draw[-] (3,0) -- (3,-1);
\draw[-] (3,0) -- (4,0);
\draw[-] (3,-1) -- (4,-1);
\draw[-] (4,0) -- (4,-1);
\draw[-] (5,0) -- (6,0);
\draw[-] (5,0) -- (5,-1);
\draw[-] (6,0) -- (6,-1);
\draw[-] (5,-1) -- (6,-1);
\draw[-] (6,0) -- (7,0);
\draw[-] (6,-1) -- (7,-1);
\draw[-] (7,0) -- (7,-1);
\draw[dashed] (1,0) -- (2,-1);
\draw[dashed] (0,0) -- (1,-1);
\draw[dashed] (2,0) -- (3,-1);
\draw[dashed] (3,0) -- (4,-1);
\draw[dashed] (5,0) -- (6,-1);
\draw[dashed] (6,0) -- (7,-1);
\draw[dashed] (0,-1) -- (1,-2);
\node[left] at (0.05,-0.5) {$a_4$};
\node[left] at (2.25,-0.5) {$a_4$};
\node[left] at (4.25,-0.5) {$a_4$};
\node[left] at (5.25,-0.5) {$a_4$};
\node[right] at (6.95,-0.5) {$a_4$};
\node[left] at (1.25,-0.5) {$a_2$};
\node[left] at (3.25,-0.5) {$a_2$};
\node[left] at (6.25,-0.5) {$a_2$};
\node[above] at (0.5,-0.05) {$a_1$};
\node[above] at (1.5,-0.05) {$\tau$};
\node[above] at (2.5,-0.05) {$a_1$};
\node[above] at (3.5,-0.05) {$\tau$};
\node[above] at (5.5,-0.05) {$a_1$};
\node[above] at (6.5,-0.05) {$\tau$};
\node[below] at (0.5,-0.75) {$a_1$};
\node[below] at (1.5,-0.95) {$\tau$};
\node[below] at (2.5,-0.95) {$a_1$};
\node[below] at (3.5,-0.95) {$\tau$};
\node[below] at (5.5,-0.95) {$a_1$};
\node[below] at (6.5,-0.95) {$\tau$};
\node[right] at (1.3,-0.45) {$a_1$};
\node[right] at (3.3,-0.45) {$a_1$};
\node[right] at (2.3,-0.5) {$\tau$};
\node[right] at (0.3,-1.5) {$a_4$};
\node[right] at (6.3,-0.45) {$a_1$};
\node[right] at (5.3,-0.5) {$\tau$};
\node[left] at (0.8,-0.5) {$\tau$};
\node[left] at (0.05,-1.5) {$a_8$};
\node[right] at (0.95,-1.5) {$\tau$};
\node[below] at (0.5,-1.95) {$a_7$};
\draw [-] (1, 0) -- (1,-1);
\draw [-] (2, 0) -- (2,-1);
\draw [fill] (0,0) circle [radius=.05];
\draw [fill] (0,-1) circle [radius=.05];
\draw [fill] (1,0) circle [radius=.05];
\draw [fill] (2,0) circle [radius=.05];
\draw [fill] (1,-1) circle [radius=.05];
\draw [fill] (2,-1) circle [radius=.05];
\draw [fill] (3,0) circle [radius=.05];
\draw [fill] (3,-1) circle [radius=.05];
\draw [fill] (4,0) circle [radius=.05];
\draw [fill] (4,-1) circle [radius=.05];
\draw [fill] (5,0) circle [radius=.05];
\draw [fill] (5,-1) circle [radius=.05];
\draw [fill] (6,0) circle [radius=.05];
\draw [fill] (6,-1) circle [radius=.05];
\draw [fill] (7,0) circle [radius=.05];
\draw [fill] (0,-2) circle [radius=.05];
\draw [fill] (1,-2) circle [radius=.05];
\draw [fill] (4.3,-0.5) circle [radius=.02];
\draw [fill] (4.5,-0.5) circle [radius=.02];
\draw [fill] (4.7,-0.5) circle [radius=.02];
\draw [fill] (7,-1) circle [radius=.05];
\end{tikzpicture}}

\centerline{\begin{tikzpicture}
\draw[-] (0,0) -- (1,0);
\draw[-] (0,0) -- (0,-1);
\draw[-] (0,-1) -- (1,-1);
\draw[-] (1,0) -- (2,0);
\draw[-] (1,-1) -- (2,-1);
\draw[-] (2,0) -- (3,0);
\draw[-] (2,-1) -- (3,-1);
\draw[-] (3,0) -- (3,-1);
\draw[-] (3,0) -- (4,0);
\draw[-] (3,-1) -- (4,-1);
\draw[-] (4,0) -- (4,-1);
\draw[-] (5,0) -- (6,0);
\draw[-] (5,0) -- (5,-1);
\draw[-] (6,0) -- (6,-1);
\draw[-] (5,-1) -- (6,-1);
\draw[-] (6,0) -- (7,0);
\draw[-] (6,-1) -- (7,-1);
\draw[-] (7,0) -- (7,-1);
\draw[-] (0,-2) -- (1,-2);
\draw[-] (1,-2) -- (1,-1);
\draw[-] (0,-2) -- (0,-1);
\draw[dashed] (1,0) -- (2,-1);
\draw[dashed] (0,0) -- (1,-1);
\draw[dashed] (2,0) -- (3,-1);
\draw[dashed] (3,0) -- (4,-1);
\draw[dashed] (5,0) -- (6,-1);
\draw[dashed] (6,0) -- (7,-1);
\draw[dashed] (0,-1) -- (1,-2);
\node[left] at (2.25,-0.5) {$a_4$};
\node[left] at (4.25,-0.5) {$a_4$};
\node[left] at (5.25,-0.5) {$a_4$};
\node[right] at (6.95,-0.5) {$a_4$};
\node[left] at (1.25,-0.5) {$a_2$};
\node[left] at (0.05,-0.5) {$a_4$};
\node[left] at (3.25,-0.5) {$a_2$};
\node[left] at (6.25,-0.5) {$a_2$};
\node[above] at (0.5,-0.05) {$\tau'$};
\node[above] at (1.5,-0.05) {$a_1$};
\node[above] at (2.5,-0.05) {$\tau'$};
\node[above] at (3.5,-0.05) {$a_1$};
\node[above] at (5.5,-0.05) {$\tau'$};
\node[above] at (6.5,-0.05) {$a_1$};
\node[below] at (0.5,-0.75) {$\tau'$};
\node[below] at (1.5,-0.95) {$a_1$};
\node[below] at (2.5,-0.95) {$\tau'$};
\node[below] at (3.5,-0.95) {$a_1$};
\node[below] at (5.5,-0.95) {$\tau'$};
\node[below] at (6.5,-0.95) {$a_1$};
\node[right] at (1.3,-0.45) {$\tau'$};
\node[right] at (3.3,-0.45) {$\tau'$};
\node[right] at (2.3,-0.5) {$a_1$};
\node[right] at (6.3,-0.45) {$\tau'$};
\node[right] at (0.3,-1.55) {$a_4$};
\node[right] at (5.3,-0.5) {$a_1$};
\node[left] at (0.8,-0.5) {$a_1$};
\node[left] at (0.05,-1.5) {$a_8$};
\node[right] at (0.95,-1.5) {$a_1$};
\node[below] at (0.5,-1.95) {$a_7$};
\draw [-] (1, 0) -- (1,-1);
\draw [-] (2, 0) -- (2,-1);
\draw [fill] (0,0) circle [radius=.05];
\draw [fill] (0,-1) circle [radius=.05];
\draw [fill] (1,0) circle [radius=.05];
\draw [fill] (2,0) circle [radius=.05];
\draw [fill] (1,-1) circle [radius=.05];
\draw [fill] (2,-1) circle [radius=.05];
\draw [fill] (3,0) circle [radius=.05];
\draw [fill] (3,-1) circle [radius=.05];
\draw [fill] (4,0) circle [radius=.05];
\draw [fill] (4,-1) circle [radius=.05];
\draw [fill] (5,0) circle [radius=.05];
\draw [fill] (5,-1) circle [radius=.05];
\draw [fill] (6,0) circle [radius=.05];
\draw [fill] (6,-1) circle [radius=.05];
\draw [fill] (7,0) circle [radius=.05];
\draw [fill] (1,-2) circle [radius=.05];
\draw [fill] (0,-2) circle [radius=.05];
\draw [fill] (4.3,-0.5) circle [radius=.02];
\draw [fill] (4.5,-0.5) circle [radius=.02];
\draw [fill] (4.7,-0.5) circle [radius=.02];
\draw [fill] (7,-1) circle [radius=.05];
\end{tikzpicture}}

We clearly have $n^{\tau}(T,\rho)=n^{\tau'}(T',\rho)$.

\medskip

Since $G_{T,\rho}$ can be obtained by gluing a tile with a graph which is isomorphic to $G_s(\tau,a_1)$, we have $\mathcal P(G_{T,\rho})=\mathcal P_1\sqcup \mathcal P_2$, where $\mathcal P_1$ contains all $P$ which contains the left lower edge labeled $a_7$, $\mathcal P_2$ contains all $P$ which contains the left lower edge labeled $a_8$. By Lemma \ref{basicdecom}, as sets, we have $$\mathcal P_1\cong\textstyle\bigsqcup_{(\lambda_1,\cdots,\lambda_{s})\in \{0,1\}^{s}}\mathcal G_{(\lambda_1,\cdots,\lambda_{s})}(\tau,a_1),$$ $$\mathcal P_2\cong\textstyle\bigsqcup_{(\lambda_2,\cdots,\lambda_{s})\in \{0,1\}^{s-1}}\mathcal G_{(\lambda_2,\cdots,\lambda_{s})}(\tau,a_1).$$

\medskip

Similarly, $\mathcal P(G_{T',\rho})=\mathcal P'_1\sqcup \mathcal P'_2$, where $\mathcal P'_1$ contains all $P$ which contains the left lower edge labeled $a_7$, $\mathcal P'_2$ contains all $P$ which contains the left lower edge labeled $a_8$. As sets, we have $$\mathcal P'_1\cong\textstyle\bigsqcup_{(\lambda_1,\cdots,\lambda_{s})\in \{0,1\}^{s}}\mathcal G'_{(\lambda_1,\cdots,\lambda_{s})}(a_1,\tau'),$$ $$\mathcal P'_2\cong\textstyle\bigsqcup_{(\lambda_2,\cdots,\lambda_{s})\in \{0,1\}^{s-1}}\mathcal G'_{(\lambda_2,\cdots,\lambda_{s})}(a_1,\tau').$$

\medskip

As case (3), under the above isomorphisms, $\psi_{\rho}: \mathcal G_{\lambda}(\tau,a_1)\rightarrow \mathcal G'_{1-\lambda}(a_1,\tau')$ gives a partition bijection between $\mathcal P(G_{T,\rho})$ and $\mathcal P(G_{T',\rho})$ for $\lambda=(\lambda_1,\cdots,\lambda_{s})$ or $(\lambda_2,\cdots,\lambda_{s})$.
For any $P\in \mathcal G_{\lambda}(\tau,a_1)$ and $P'\in \mathcal G'_{1-\lambda}(a_1,\tau')$, by Lemma \ref{same-num1} the number of non-$\tau$-mutable edges labeled $a$ in $P$ equals to the number of non-$\tau'$-mutable edges labeled $a$ in $P'$ for any $\tau \neq a\in T$. Moreover, the first gluing edge of $G_{T,\rho}$ is in $P$ if and only if the first gluing edge of $G_{T',\rho}$ in $P'$. Since $\rho$ and $\zeta$ have the same endpoint, $G_{T,\rho}$ and $G_{T',\rho}$ do not have last gluing edge.

\medskip

In case (5) we may assume $i=4$. Up to a difference of relative orientation, $G_{T,\rho}$ and $G_{T',\rho}$ are the following graphs, respectively,

\centerline{\begin{tikzpicture}
\draw[-] (0,-1) -- (0,-2);
\draw[-] (1,-1) -- (1,-2);
\draw[-] (0,-2) -- (1,-2);
\draw[-] (0,0) -- (1,0);
\draw[-] (0,0) -- (0,-1);
\draw[-] (0,-1) -- (1,-1);
\draw[-] (1,0) -- (2,0);
\draw[-] (1,-1) -- (2,-1);
\draw[-] (2,0) -- (3,0);
\draw[-] (2,-1) -- (3,-1);
\draw[-] (3,0) -- (3,-1);
\draw[-] (3,0) -- (4,0);
\draw[-] (3,-1) -- (4,-1);
\draw[-] (4,0) -- (4,-1);
\draw[-] (5,0) -- (6,0);
\draw[-] (5,0) -- (5,-1);
\draw[-] (6,0) -- (6,-1);
\draw[-] (5,-1) -- (6,-1);
\draw[-] (6,0) -- (7,0);
\draw[-] (6,-1) -- (7,-1);
\draw[-] (7,0) -- (7,-1);
\draw[-] (6,1) -- (7,1);
\draw[-] (6,0) -- (6,1);
\draw[-] (7,0) -- (7,1);
\draw[dashed] (1,0) -- (2,-1);
\draw[dashed] (0,0) -- (1,-1);
\draw[dashed] (0,-1) -- (1,-2);
\draw[dashed] (2,0) -- (3,-1);
\draw[dashed] (3,0) -- (4,-1);
\draw[dashed] (5,0) -- (6,-1);
\draw[dashed] (6,0) -- (7,-1);
\draw[dashed] (6,1) -- (7,0);
\node[left] at (0.05,-0.5) {$a_4$};
\node[left] at (2.25,-0.5) {$a_4$};
\node[left] at (4.25,-0.5) {$a_4$};
\node[left] at (5.25,-0.5) {$a_4$};
\node[right] at (6.95,-0.5) {$a_4$};
\node[right] at (6.95,0.5) {$a_7$};
\node[left] at (1.25,-0.5) {$a_2$};
\node[left] at (3.25,-0.5) {$a_2$};
\node[left] at (6.25,-0.5) {$a_2$};
\node[above] at (0.5,-0.05) {$a_1$};
\node[above] at (1.5,-0.05) {$\tau$};
\node[above] at (2.5,-0.05) {$a_1$};
\node[above] at (3.5,-0.05) {$\tau$};
\node[above] at (5.5,-0.05) {$a_1$};
\node[above] at (6.5,-0.25) {$\tau$};
\node[below] at (0.5,-0.8) {$a_1$};
\node[below] at (1.5,-0.95) {$\tau$};
\node[below] at (2.5,-0.95) {$a_1$};
\node[below] at (3.5,-0.95) {$\tau$};
\node[below] at (0.5,-1.95) {$a_7$};
\node[below] at (5.5,-0.95) {$a_1$};
\node[below] at (6.5,-0.95) {$\tau$};
\node[right] at (1.3,-0.45) {$a_1$};
\node[right] at (3.3,-0.45) {$a_1$};
\node[right] at (2.3,-0.5) {$\tau$};
\node[right] at (0.95,-1.5) {$\tau$};
\node[right] at (6.3,0.55) {$a_4$};
\node[right] at (6.3,-0.45) {$a_1$};
\node[right] at (5.3,-0.5) {$\tau$};
\node[left] at (0.8,-0.5) {$\tau$};
\node[left] at (0.05,-1.5) {$a_8$};
\node[left] at (0.8,-1.5) {$a_4$};
\node[left] at (6.05,0.5) {$a_1$};
\node[above] at (6.5,0.95) {$a_8$};
\draw [-] (1, 0) -- (1,-1);
\draw [-] (2, 0) -- (2,-1);
\draw [fill] (0,0) circle [radius=.05];
\draw [fill] (0,-1) circle [radius=.05];
\draw [fill] (0,-2) circle [radius=.05];
\draw [fill] (1,-2) circle [radius=.05];
\draw [fill] (1,0) circle [radius=.05];
\draw [fill] (2,0) circle [radius=.05];
\draw [fill] (1,-1) circle [radius=.05];
\draw [fill] (2,-1) circle [radius=.05];
\draw [fill] (3,0) circle [radius=.05];
\draw [fill] (3,-1) circle [radius=.05];
\draw [fill] (4,0) circle [radius=.05];
\draw [fill] (4,-1) circle [radius=.05];
\draw [fill] (5,0) circle [radius=.05];
\draw [fill] (5,-1) circle [radius=.05];
\draw [fill] (6,0) circle [radius=.05];
\draw [fill] (6,-1) circle [radius=.05];
\draw [fill] (7,0) circle [radius=.05];
\draw [fill] (7,1) circle [radius=.05];
\draw [fill] (6,1) circle [radius=.05];
\draw [fill] (4.3,-0.5) circle [radius=.02];
\draw [fill] (4.5,-0.5) circle [radius=.02];
\draw [fill] (4.7,-0.5) circle [radius=.02];
\draw [fill] (7,-1) circle [radius=.05];
\end{tikzpicture}}

\centerline{\begin{tikzpicture}
\draw[-] (0,-1) -- (0,-2);
\draw[-] (1,-1) -- (1,-2);
\draw[-] (0,-2) -- (1,-2);
\draw[-] (0,0) -- (1,0);
\draw[-] (0,0) -- (0,-1);
\draw[-] (0,-1) -- (1,-1);
\draw[-] (1,0) -- (2,0);
\draw[-] (1,-1) -- (2,-1);
\draw[-] (2,0) -- (3,0);
\draw[-] (2,-1) -- (3,-1);
\draw[-] (3,0) -- (3,-1);
\draw[-] (3,0) -- (4,0);
\draw[-] (3,-1) -- (4,-1);
\draw[-] (4,0) -- (4,-1);
\draw[-] (5,0) -- (6,0);
\draw[-] (5,0) -- (5,-1);
\draw[-] (6,0) -- (6,-1);
\draw[-] (5,-1) -- (6,-1);
\draw[-] (6,0) -- (7,0);
\draw[-] (6,-1) -- (7,-1);
\draw[-] (7,0) -- (7,-1);
\draw[-] (6,1) -- (7,1);
\draw[-] (6,0) -- (6,1);
\draw[-] (7,0) -- (7,1);
\draw[dashed] (1,0) -- (2,-1);
\draw[dashed] (0,0) -- (1,-1);
\draw[dashed] (0,-1) -- (1,-2);
\draw[dashed] (2,0) -- (3,-1);
\draw[dashed] (3,0) -- (4,-1);
\draw[dashed] (5,0) -- (6,-1);
\draw[dashed] (6,0) -- (7,-1);
\draw[dashed] (6,1) -- (7,0);
\node[left] at (0.05,-0.5) {$a_4$};
\node[left] at (2.25,-0.5) {$a_4$};
\node[left] at (4.25,-0.5) {$a_4$};
\node[left] at (5.25,-0.5) {$a_4$};
\node[right] at (6.95,-0.5) {$a_4$};
\node[right] at (6.95,0.5) {$a_7$};
\node[left] at (1.25,-0.5) {$a_2$};
\node[left] at (3.25,-0.5) {$a_2$};
\node[left] at (6.25,-0.5) {$a_2$};
\node[above] at (0.5,-0.05) {$\tau'$};
\node[above] at (1.5,-0.05) {$a_1$};
\node[above] at (2.5,-0.05) {$\tau'$};
\node[above] at (3.5,-0.05) {$a_1$};
\node[above] at (5.5,-0.05) {$\tau'$};
\node[above] at (6.5,-0.25) {$a_1$};
\node[below] at (0.5,-0.8) {$\tau'$};
\node[below] at (1.5,-0.95) {$a_1$};
\node[below] at (2.5,-0.95) {$\tau'$};
\node[below] at (3.5,-0.95) {$a_1$};
\node[below] at (0.5,-1.95) {$a_7$};
\node[below] at (5.5,-0.95) {$\tau'$};
\node[below] at (6.5,-0.95) {$a_1$};
\node[right] at (1.3,-0.45) {$\tau'$};
\node[right] at (3.3,-0.45) {$\tau'$};
\node[right] at (2.3,-0.5) {$a_1$};
\node[right] at (0.95,-1.5) {$a_1$};
\node[right] at (6.3,0.55) {$a_4$};
\node[right] at (6.3,-0.45) {$\tau'$};
\node[right] at (5.3,-0.5) {$a_1$};
\node[left] at (0.8,-0.5) {$a_1$};
\node[left] at (0.05,-1.5) {$a_8$};
\node[left] at (0.8,-1.5) {$a_4$};
\node[left] at (6.05,0.5) {$\tau'$};
\node[above] at (6.5,0.95) {$a_8$};
\draw [-] (1, 0) -- (1,-1);
\draw [-] (2, 0) -- (2,-1);
\draw [fill] (0,0) circle [radius=.05];
\draw [fill] (0,-1) circle [radius=.05];
\draw [fill] (0,-2) circle [radius=.05];
\draw [fill] (1,-2) circle [radius=.05];
\draw [fill] (1,0) circle [radius=.05];
\draw [fill] (2,0) circle [radius=.05];
\draw [fill] (1,-1) circle [radius=.05];
\draw [fill] (2,-1) circle [radius=.05];
\draw [fill] (3,0) circle [radius=.05];
\draw [fill] (3,-1) circle [radius=.05];
\draw [fill] (4,0) circle [radius=.05];
\draw [fill] (4,-1) circle [radius=.05];
\draw [fill] (5,0) circle [radius=.05];
\draw [fill] (5,-1) circle [radius=.05];
\draw [fill] (6,0) circle [radius=.05];
\draw [fill] (6,-1) circle [radius=.05];
\draw [fill] (7,0) circle [radius=.05];
\draw [fill] (7,1) circle [radius=.05];
\draw [fill] (6,1) circle [radius=.05];
\draw [fill] (4.3,-0.5) circle [radius=.02];
\draw [fill] (4.5,-0.5) circle [radius=.02];
\draw [fill] (4.7,-0.5) circle [radius=.02];
\draw [fill] (7,-1) circle [radius=.05];
\end{tikzpicture}}

We clearly have $n^{\tau}(T,\rho)=n^{\tau'}(T',\rho)$.

\medskip

 As $G_{T,\rho}$ can be obtained by gluing two tiles with a graph which is isomorphic to $G_s(\tau,a_1)$, we have $\mathcal P(G_{T,\rho})=\mathcal P_1\sqcup \mathcal P_2\sqcup \mathcal P_3\sqcup \mathcal P_4$, where $\mathcal P_1$ contains all $P$ which contains the lower left edge labeled $a_7$ and the upper right edge labeled $a_8$, $\mathcal P_2$ contains all $P$ which contains the lower left edge labeled $a_7$ and the upper right edge labeled $a_7$, $\mathcal P_3$ contains all $P$ which contains the lower left edge labeled $a_8$ and the upper right edge labeled $a_8$, $\mathcal P_4$ contains all $P$ which contains the lower left edge labeled $a_8$ and the upper right edge labeled $a_7$. By Lemma \ref{basicdecom}, as sets, we have
$$\mathcal P_1\cong\textstyle\bigsqcup_{(\lambda_1,\cdots,\lambda_{s})\in \{0,1\}^{s}}\mathcal G_{(\lambda_1,\cdots,\lambda_{s})}(\tau,a_1),$$
$$\mathcal P_2\cong\textstyle\bigsqcup_{(\lambda_1,\cdots,\lambda_{s-1})\in \{0,1\}^{s-1}}\mathcal G_{(\lambda_1,\cdots,\lambda_{s-1})}(\tau,a_1),$$
$$\mathcal P_3\cong\textstyle\bigsqcup_{(\lambda_2,\cdots,\lambda_{s})\in \{0,1\}^{s-1}}\mathcal G_{(\lambda_2,\cdots,\lambda_{s})}(\tau,a_1),$$
$$\mathcal P_4\cong\textstyle\bigsqcup_{(\lambda_2,\cdots,\lambda_{s-1})\in \{0,1\}^{s-2}}\mathcal G_{(\lambda_2,\cdots,\lambda_{s-1})}(\tau,a_1).$$
\medskip

Similarly, $\mathcal P(G_{T',\rho})=\mathcal P'_1\sqcup \mathcal P'_2\sqcup \mathcal P'_3\sqcup \mathcal P'_4$, where $\mathcal P'_1$ contains all $P$ which contains the lower left edge labeled $a_7$ and the upper right edge labeled $a_8$, $\mathcal P'_2$ contains all $P$ which contains the lower left edge labeled $a_7$ and the upper right edge labeled $a_7$, $\mathcal P'_3$ contains all $P$ which contains the lower left edge labeled $a_8$ and the upper right edge labeled $a_8$, $\mathcal P'_4$ contains all $P$ which contains the lower left edge labeled $a_8$ and the upper right edge labeled $a_7$. As sets, we have $$\mathcal P'_1\cong\textstyle\bigsqcup_{(\lambda_1,\cdots,\lambda_{s})\in \{0,1\}^{s}}\mathcal G'_{(\lambda_1,\cdots,\lambda_{s})}(a_1,\tau'),$$
$$\mathcal P'_2\cong\textstyle\bigsqcup_{(\lambda_1,\cdots,\lambda_{s-1})\in \{0,1\}^{s-1}}\mathcal G'_{(\lambda_1,\cdots,\lambda_{s-1})}(a_1,\tau'),$$
$$\mathcal P'_3\cong\textstyle\bigsqcup_{(\lambda_2,\cdots,\lambda_{s})\in \{0,1\}^{s-1}}\mathcal G'_{(\lambda_2,\cdots,\lambda_{s})}(a_1,\tau'),$$
$$\mathcal P'_4\cong\textstyle\bigsqcup_{(\lambda_2,\cdots,\lambda_{s-1})\in \{0,1\}^{s-2}}\mathcal G'_{(\lambda_2,\cdots,\lambda_{s-1})}(a_1,\tau').$$

\medskip

As cases (3,4), under the above isomorphisms, $\psi_{\rho}: \mathcal G_{\lambda}(\tau,a_1)\rightarrow \mathcal G'_{1-\lambda}(a_1,\tau')$ gives a partition bijection between $\mathcal P(G_{T,\rho})$ and $\mathcal P(G_{T',\rho})$ for $\lambda=(\lambda_1,\cdots,\lambda_{s}),(\lambda_1,\cdots,\lambda_{s-1})$, $(\lambda_2,\cdots,\lambda_{s})$ or $(\lambda_2,\cdots,\lambda_{s-1})$. For any $P\in \mathcal G_{\lambda}(\tau,a_1)$ and $P'\in \mathcal G'_{1-\lambda}(a_1,\tau')$, by Lemma \ref{same-num1} the number of non-$\tau$-mutable edges labeled $a$ in $P$ equals to the number of non-$\tau'$-mutable edges labeled $a$ in $P'$ for any $\tau \neq a\in T$. Moreover, the first/last gluing edge of $G_{T,\rho}$ is in $P$ if and only if the first/last gluing edge of $G_{T',\rho}$ in $P'$.

\medskip

In case (6), $\rho=\zeta$. Up to a difference of relative orientation, $G_{T,\rho}$ and $G_{T',\rho}$ are the following graphs, respectively,

\centerline{\begin{tikzpicture}
\draw[-] (0,0) -- (1,0);
\draw[-] (0,0) -- (0,-1);
\draw[-] (0,-1) -- (1,-1);
\draw[-] (1,0) -- (2,0);
\draw[-] (1,-1) -- (2,-1);
\draw[-] (2,0) -- (3,0);
\draw[-] (2,-1) -- (3,-1);
\draw[-] (3,0) -- (3,-1);
\draw[-] (3,0) -- (4,0);
\draw[-] (3,-1) -- (4,-1);
\draw[-] (4,0) -- (4,-1);
\draw[-] (5,0) -- (6,0);
\draw[-] (5,0) -- (5,-1);
\draw[-] (6,0) -- (6,-1);
\draw[-] (5,-1) -- (6,-1);
\draw[-] (6,0) -- (7,0);
\draw[-] (6,-1) -- (7,-1);
\draw[-] (7,0) -- (7,-1);
\draw[-] (7,0) -- (8,0);
\draw[-] (7,-1) -- (8,-1);
\draw[-] (8,0) -- (8,-1);
\draw[dashed] (1,0) -- (2,-1);
\draw[dashed] (0,0) -- (1,-1);
\draw[dashed] (2,0) -- (3,-1);
\draw[dashed] (3,0) -- (4,-1);
\draw[dashed] (5,0) -- (6,-1);
\draw[dashed] (6,0) -- (7,-1);
\draw[dashed] (7,0) -- (8,-1);
\node[left] at (0.05,-0.5) {$a_2$};
\node[left] at (2.25,-0.5) {$a_2$};
\node[left] at (4.25,-0.5) {$a_2$};
\node[left] at (5.25,-0.5) {$a_2$};
\node[left] at (7.25,-0.5) {$a_2$};
\node[left] at (1.25,-0.5) {$a_4$};
\node[left] at (3.25,-0.5) {$a_4$};
\node[left] at (6.25,-0.5) {$a_4$};
\node[right] at (7.95,-0.5) {$a_4$};
\node[above] at (0.5,-0.05) {$\tau$};
\node[above] at (1.5,-0.05) {$a_1$};
\node[above] at (2.5,-0.05) {$\tau$};
\node[above] at (3.5,-0.05) {$a_1$};
\node[above] at (5.5,-0.05) {$\tau$};
\node[above] at (6.5,-0.05) {$a_1$};
\node[above] at (7.5,-0.05) {$\tau$};
\node[below] at (0.5,-0.95) {$\tau$};
\node[below] at (1.5,-0.95) {$a_1$};
\node[below] at (2.5,-0.95) {$\tau$};
\node[below] at (3.5,-0.95) {$a_1$};
\node[below] at (5.5,-0.95) {$\tau$};
\node[below] at (6.5,-0.95) {$a_1$};
\node[below] at (7.5,-0.95) {$\tau$};
\node[right] at (1.3,-0.45) {$\tau$};
\node[right] at (3.3,-0.45) {$\tau$};
\node[right] at (2.3,-0.5) {$a_1$};
\node[right] at (6.3,-0.45) {$\tau$};
\node[right] at (5.3,-0.5) {$a_1$};
\node[right] at (7.3,-0.5) {$a_1$};
\node[left] at (0.8,-0.5) {$a_1$};
\draw [-] (1, 0) -- (1,-1);
\draw [-] (2, 0) -- (2,-1);
\draw [fill] (0,0) circle [radius=.05];
\draw [fill] (0,-1) circle [radius=.05];
\draw [fill] (1,0) circle [radius=.05];
\draw [fill] (2,0) circle [radius=.05];
\draw [fill] (1,-1) circle [radius=.05];
\draw [fill] (2,-1) circle [radius=.05];
\draw [fill] (3,0) circle [radius=.05];
\draw [fill] (3,-1) circle [radius=.05];
\draw [fill] (4,0) circle [radius=.05];
\draw [fill] (4,-1) circle [radius=.05];
\draw [fill] (5,0) circle [radius=.05];
\draw [fill] (5,-1) circle [radius=.05];
\draw [fill] (6,0) circle [radius=.05];
\draw [fill] (6,-1) circle [radius=.05];
\draw [fill] (7,0) circle [radius=.05];
\draw [fill] (8,0) circle [radius=.05];
\draw [fill] (8,-1) circle [radius=.05];
\draw [fill] (4.3,-0.5) circle [radius=.02];
\draw [fill] (4.5,-0.5) circle [radius=.02];
\draw [fill] (4.7,-0.5) circle [radius=.02];
\draw [fill] (7,-1) circle [radius=.05];
\end{tikzpicture}}

\medskip

\centerline{\begin{tikzpicture}
\draw[-] (-1,0) -- (0,0);
\draw[-] (-1,-1) -- (0,-1);
\draw[-] (-1,0) -- (-1,-1);
\draw[-] (0,0) -- (1,0);
\draw[-] (0,0) -- (0,-1);
\draw[-] (0,-1) -- (1,-1);
\draw[-] (1,0) -- (2,0);
\draw[-] (1,-1) -- (2,-1);
\draw[-] (2,0) -- (3,0);
\draw[-] (2,-1) -- (3,-1);
\draw[-] (3,0) -- (3,-1);
\draw[-] (3,0) -- (4,0);
\draw[-] (3,-1) -- (4,-1);
\draw[-] (4,0) -- (4,-1);
\draw[-] (5,0) -- (6,0);
\draw[-] (5,0) -- (5,-1);
\draw[-] (6,0) -- (6,-1);
\draw[-] (5,-1) -- (6,-1);
\draw[-] (6,0) -- (7,0);
\draw[-] (6,-1) -- (7,-1);
\draw[-] (7,0) -- (7,-1);
\draw[-] (7,0) -- (8,0);
\draw[-] (7,-1) -- (8,-1);
\draw[-] (8,0) -- (8,-1);
\draw[-] (8,0) -- (9,0);
\draw[-] (8,-1) -- (9,-1);
\draw[-] (9,0) -- (9,-1);
\draw[dashed] (-1,0) -- (0,-1);
\draw[dashed] (1,0) -- (2,-1);
\draw[dashed] (0,0) -- (1,-1);
\draw[dashed] (2,0) -- (3,-1);
\draw[dashed] (3,0) -- (4,-1);
\draw[dashed] (5,0) -- (6,-1);
\draw[dashed] (6,0) -- (7,-1);
\draw[dashed] (7,0) -- (8,-1);
\draw[dashed] (8,0) -- (9,-1);
\node[left] at (-0.95,-0.5) {$a_2$};
\node[left] at (2.25,-0.5) {$a_4$};
\node[left] at (4.25,-0.5) {$a_4$};
\node[left] at (5.25,-0.5) {$a_4$};
\node[left] at (7.25,-0.5) {$a_4$};
\node[left] at (1.25,-0.5) {$a_2$};
\node[left] at (0.25,-0.5) {$a_4$};
\node[left] at (3.25,-0.5) {$a_2$};
\node[left] at (6.25,-0.5) {$a_2$};
\node[left] at (8.25,-0.5) {$a_2$};
\node[right] at (8.95,-0.5) {$a_4$};
\node[above] at (0.5,-0.05) {$\tau'$};
\node[above] at (1.5,-0.05) {$a_1$};
\node[above] at (-0.5,-0.05) {$a_1$};
\node[above] at (2.5,-0.05) {$\tau'$};
\node[above] at (3.5,-0.05) {$a_1$};
\node[above] at (5.5,-0.05) {$\tau'$};
\node[above] at (6.5,-0.05) {$a_1$};
\node[above] at (8.5,-0.05) {$a_1$};
\node[above] at (7.5,-0.05) {$\tau'$};
\node[below] at (0.5,-0.95) {$\tau'$};
\node[below] at (1.5,-0.95) {$a_1$};
\node[below] at (-0.5,-0.95) {$a_1$};
\node[below] at (2.5,-0.95) {$\tau'$};
\node[below] at (3.5,-0.95) {$a_1$};
\node[below] at (5.5,-0.95) {$\tau'$};
\node[below] at (6.5,-0.95) {$a_1$};
\node[below] at (8.5,-0.95) {$a_1$};
\node[below] at (7.5,-0.95) {$\tau'$};
\node[right] at (1.3,-0.45) {$\tau'$};
\node[right] at (-0.7,-0.45) {$\tau'$};
\node[right] at (3.3,-0.45) {$\tau'$};
\node[right] at (2.3,-0.5) {$a_1$};
\node[right] at (6.3,-0.45) {$\tau'$};
\node[right] at (8.3,-0.45) {$\tau'$};
\node[right] at (5.3,-0.5) {$a_1$};
\node[right] at (7.3,-0.5) {$a_1$};
\node[left] at (0.8,-0.5) {$a_1$};
\draw [-] (1, 0) -- (1,-1);
\draw [-] (2, 0) -- (2,-1);
\draw [fill] (0,0) circle [radius=.05];
\draw [fill] (0,-1) circle [radius=.05];
\draw [fill] (-1,0) circle [radius=.05];
\draw [fill] (-1,-1) circle [radius=.05];
\draw [fill] (1,0) circle [radius=.05];
\draw [fill] (2,0) circle [radius=.05];
\draw [fill] (1,-1) circle [radius=.05];
\draw [fill] (2,-1) circle [radius=.05];
\draw [fill] (3,0) circle [radius=.05];
\draw [fill] (3,-1) circle [radius=.05];
\draw [fill] (4,0) circle [radius=.05];
\draw [fill] (4,-1) circle [radius=.05];
\draw [fill] (5,0) circle [radius=.05];
\draw [fill] (5,-1) circle [radius=.05];
\draw [fill] (6,0) circle [radius=.05];
\draw [fill] (6,-1) circle [radius=.05];
\draw [fill] (7,0) circle [radius=.05];
\draw [fill] (8,0) circle [radius=.05];
\draw [fill] (8,-1) circle [radius=.05];
\draw [fill] (9,0) circle [radius=.05];
\draw [fill] (9,-1) circle [radius=.05];
\draw [fill] (4.3,-0.5) circle [radius=.02];
\draw [fill] (4.5,-0.5) circle [radius=.02];
\draw [fill] (4.7,-0.5) circle [radius=.02];
\draw [fill] (7,-1) circle [radius=.05];
\end{tikzpicture}}

We clearly have $n^{\tau}(T,\rho)=n^{\tau'}(T',\rho)$.

\medskip

Since $G_{T,\rho}$ is isomorphic to $H_{s}(a_1,\tau)$ and $G_{T',\rho}$ is isomorphic to $H_{s+1}(\tau',a_1)$. By Lemma \ref{basicdecom}, as sets, we have $\mathcal P(G_{T,\rho})\cong\bigsqcup_{(\lambda_1,\cdots,\lambda_{s+1})\in \{0,1\}^{s+1}}\mathcal H'_{(\lambda_1,\cdots,\lambda_{s+1})}(a_1,\tau)$, and  $\mathcal P(G_{T',\rho})\cong\bigsqcup_{(\lambda_1,\cdots,\lambda_{s+1})\in \{0,1\}^{s+1}}\mathcal H_{(\lambda_1,\cdots,\lambda_{s+1})}(\tau',a_1)$.

\medskip

We have the following observation.

\medskip

\begin{Lemma}\label{non-tau-mu1}

Let $P\in \mathcal H'_{(\lambda_1,\cdots,\lambda_{s+1})}(a_1,\tau)$ for some sequence $(\lambda_1,\cdots,\lambda_{s+1})$. We assume $\lambda_{0}=1$ in convention. Then for any $i\in [1,s+1]$,

 \begin{enumerate}[$(1)$]

   \item the edge labeled $a_2$ of the $(2i-1)$-th tile is in $P$ but non-$\tau$-mutable if and only if $\lambda_{i-1}=1$ and $\lambda_i=0$;

   \item the edge labeled $a_4$ of the $(2i-1)$-th tile is in $P$ but non-$\tau$-mutable if and only if $\lambda_{i}=0$ and $\lambda_{i+1}=1$.

 \end{enumerate}

\end{Lemma}

\begin{proof}

By symmetry, we shall only prove (1).

``Only If Part:"
Since the edge labeled $a_2$ of the $(2i-1)$-th tile is in $P$ but non-$\tau$-mutable, the edges labeled $a_1,a_4$ of the $(2i-2)$-th tile and the edges labeled $\tau$ of the $(2i-1)$-th tile are not in $P$. Thus the edges labeled $\tau$ of the $(2i-3)$-th tile are in $P$. Therefore, $\lambda_{i-1}=1$ and $\lambda_i=0$.

``If Part:" Since $\lambda_{i-1}=1$, the edges labeled $\tau$ of the $(2i-3)$-th tile are in $P$, and hence the edges labeled $a_1$ and $a_4$ of the $(2i-2)$-th tile are not in $P$. Since $\lambda_i=0$, the edges labeled $\tau$ of the $(2i-1)$-th tile are not in $P$. Therefore, the edge labeled $a_2$ of the $(2i-1)$-th tile is in $P$ but non-$\tau$-mutable.
\end{proof}

\medskip

The following result can be verified similarly.

\medskip

\begin{Lemma}\label{non-tau-mu2}

Let $P\in \mathcal H_{(\lambda_1,\cdots,\lambda_{s+1})}(\tau',a_1)$ for some sequence $(\lambda_1,\cdots,\lambda_{s+1})$. We assume $\lambda_{0}=\lambda_{s+2}=0$ in convention. Then for any $i\in [1,s+2]$,

 \begin{enumerate}[$(1)$]

   \item the edge labeled $a_2$ of the $(2i-1)$-th tile is in $P$ but non-$\tau'$-mutable if and only if $\lambda_{i-1}=0$ and $\lambda_i=1$;

   \item the edge labeled $a_4$ of the $(2i-1)$-th tile is in $P$ but non-$\tau'$-mutable if and only if $\lambda_{i}=1$ and $\lambda_{i+1}=0$.

 \end{enumerate}

\end{Lemma}

\medskip

\begin{Lemma}\label{same-num}

Let $P\in \mathcal H'_{(\lambda_1,\cdots,\lambda_{s+1})}(a_1,\tau)$ and $P'\in\mathcal H_{(1-\lambda_1,\cdots,1-\lambda_{s+1})}(\tau',a_1)$. For each edge $\tau \neq a\in T$, the number of the non-$\tau$-mutable edges labeled $a$ in $P$ equals to the number of the non-$\tau'$-mutable edges labeled $a$ in $P'$.

\end{Lemma}

\begin{proof}

By Lemma \ref{inone}, if $a\neq a_2,a_4$, then the two numbers both are equal to $0$. If $a=a_2$, by Lemma \ref{non-tau-mu1} and Lemma \ref{non-tau-mu2}, for each $i\in [1,s+1]$, the edge labeled $a_2$ of the $(2i-1)$-th tile of $G_{T,\rho}$ is non-$\tau$-mutable in $P$ if and only if the edge labeled $a_2$ of the $(2i-1)$-th tile of $G_{T',\rho}$ is non-$\tau$-mutable in $P'$. By Lemma \ref{non-tau-mu2}, the edge labeled $a_2$ of the $(2(s+2)-1)$-th tile of $G_{T,\rho}$ is not a non-$\tau$-mutable edge in $P$. Therefore, the result holds for $a=a_2$. Similarly, the result holds for $a=a_4$.
\end{proof}

\medskip

$\psi_{\rho}: \mathcal H'_{(\lambda_1,\cdots,\lambda_{s+1})}(a_1,\tau)\rightarrow \mathcal H_{(1-\lambda_1,\cdots,1-\lambda_{s+1})}(\tau',a_1)$ clearly gives a partition bijection between $\mathcal P(G_{T,\rho})$ and $\mathcal P(G_{T',\rho})$. Moreover, for each $P\in \mathcal H'_{(\lambda_1,\cdots,\lambda_{s+1})}(a_1,\tau)$ and $P'\in \mathcal H_{(1-\lambda_1,\cdots,1-\lambda_{s+1})}(\tau',a_1)$, by Lemma \ref{same-num} the number of non-$\tau$-mutable edges labeled $a$ in $P$ equals to the number of non-$\tau'$-mutable edges labeled $a$ in $P'$ for any $\tau \neq a\in T$.

\medskip

In this case, since $\rho=\zeta$, $G_{T,\rho}$ and $G_{T',\rho}$ have no first or last gluing edge.

\medskip

In case (7), $\rho$ and $\zeta$ have the same starting point. We may assume $i=4$. Then up to a difference of relative orientation, $G_{T,\rho}$ and $G_{T',\rho}$ are the following graphs, respectively,

\centerline{\begin{tikzpicture}
\draw[-] (0,0) -- (1,0);
\draw[-] (0,0) -- (0,-1);
\draw[-] (0,-1) -- (1,-1);
\draw[-] (1,0) -- (2,0);
\draw[-] (1,-1) -- (2,-1);
\draw[-] (2,0) -- (3,0);
\draw[-] (2,-1) -- (3,-1);
\draw[-] (3,0) -- (3,-1);
\draw[-] (3,0) -- (4,0);
\draw[-] (3,-1) -- (4,-1);
\draw[-] (4,0) -- (4,-1);
\draw[-] (5,0) -- (6,0);
\draw[-] (5,0) -- (5,-1);
\draw[-] (6,0) -- (6,-1);
\draw[-] (5,-1) -- (6,-1);
\draw[-] (6,0) -- (7,0);
\draw[-] (6,-1) -- (7,-1);
\draw[-] (7,0) -- (7,-1);
\draw[-] (7,0) -- (8,0);
\draw[-] (7,-1) -- (8,-1);
\draw[-] (8,0) -- (8,-1);
\draw[-] (7,1) -- (8,1);
\draw[-] (7,0) -- (7,1);
\draw[-] (8,0) -- (8,1);
\draw[dashed] (1,0) -- (2,-1);
\draw[dashed] (0,0) -- (1,-1);
\draw[dashed] (2,0) -- (3,-1);
\draw[dashed] (3,0) -- (4,-1);
\draw[dashed] (5,0) -- (6,-1);
\draw[dashed] (6,0) -- (7,-1);
\draw[dashed] (7,0) -- (8,-1);
\draw[dashed] (7,1) -- (8,0);
\node[left] at (0.05,-0.5) {$a_2$};
\node[left] at (2.25,-0.5) {$a_2$};
\node[left] at (4.25,-0.5) {$a_2$};
\node[left] at (5.25,-0.5) {$a_2$};
\node[left] at (7.25,-0.5) {$a_2$};
\node[left] at (1.25,-0.5) {$a_4$};
\node[left] at (3.25,-0.5) {$a_4$};
\node[left] at (6.25,-0.5) {$a_4$};
\node[right] at (7.95,-0.5) {$a_4$};
\node[above] at (0.5,-0.05) {$\tau$};
\node[above] at (1.5,-0.05) {$a_1$};
\node[above] at (2.5,-0.05) {$\tau$};
\node[above] at (3.5,-0.05) {$a_1$};
\node[above] at (5.5,-0.05) {$\tau$};
\node[above] at (6.5,-0.05) {$a_1$};
\node[above] at (7.5,-0.2) {$\tau$};
\node[below] at (0.5,-0.95) {$\tau$};
\node[below] at (1.5,-0.95) {$a_1$};
\node[below] at (2.5,-0.95) {$\tau$};
\node[below] at (3.5,-0.95) {$a_1$};
\node[below] at (5.5,-0.95) {$\tau$};
\node[below] at (6.5,-0.95) {$a_1$};
\node[below] at (7.5,-0.95) {$\tau$};
\node[right] at (1.3,-0.45) {$\tau$};
\node[right] at (3.3,-0.45) {$\tau$};
\node[right] at (2.3,-0.5) {$a_1$};
\node[right] at (7.3,0.5) {$a_4$};
\node[right] at (6.3,-0.45) {$\tau$};
\node[right] at (5.3,-0.5) {$a_1$};
\node[right] at (7.3,-0.5) {$a_1$};
\node[left] at (0.8,-0.5) {$a_1$};
\node[left] at (7.05,0.5) {$a_1$};
\node[right] at (7.95,0.5) {$a_8$};
\node[above] at (7.5,0.95) {$a_7$};
\draw [-] (1, 0) -- (1,-1);
\draw [-] (2, 0) -- (2,-1);
\draw [fill] (0,0) circle [radius=.05];
\draw [fill] (0,-1) circle [radius=.05];
\draw [fill] (1,0) circle [radius=.05];
\draw [fill] (2,0) circle [radius=.05];
\draw [fill] (1,-1) circle [radius=.05];
\draw [fill] (2,-1) circle [radius=.05];
\draw [fill] (3,0) circle [radius=.05];
\draw [fill] (3,-1) circle [radius=.05];
\draw [fill] (4,0) circle [radius=.05];
\draw [fill] (4,-1) circle [radius=.05];
\draw [fill] (5,0) circle [radius=.05];
\draw [fill] (5,-1) circle [radius=.05];
\draw [fill] (6,0) circle [radius=.05];
\draw [fill] (6,-1) circle [radius=.05];
\draw [fill] (7,0) circle [radius=.05];
\draw [fill] (8,0) circle [radius=.05];
\draw [fill] (8,-1) circle [radius=.05];
\draw [fill] (7,1) circle [radius=.05];
\draw [fill] (8,1) circle [radius=.05];
\draw [fill] (4.3,-0.5) circle [radius=.02];
\draw [fill] (4.5,-0.5) circle [radius=.02];
\draw [fill] (4.7,-0.5) circle [radius=.02];
\draw [fill] (7,-1) circle [radius=.05];
\end{tikzpicture}}

\centerline{\begin{tikzpicture}
\draw[-] (-1,0) -- (0,0);
\draw[-] (-1,-1) -- (0,-1);
\draw[-] (-1,0) -- (-1,-1);
\draw[-] (0,0) -- (1,0);
\draw[-] (0,0) -- (0,-1);
\draw[-] (0,-1) -- (1,-1);
\draw[-] (1,0) -- (2,0);
\draw[-] (1,-1) -- (2,-1);
\draw[-] (2,0) -- (3,0);
\draw[-] (2,-1) -- (3,-1);
\draw[-] (3,0) -- (3,-1);
\draw[-] (3,0) -- (4,0);
\draw[-] (3,-1) -- (4,-1);
\draw[-] (4,0) -- (4,-1);
\draw[-] (5,0) -- (6,0);
\draw[-] (5,0) -- (5,-1);
\draw[-] (6,0) -- (6,-1);
\draw[-] (5,-1) -- (6,-1);
\draw[-] (6,0) -- (7,0);
\draw[-] (6,-1) -- (7,-1);
\draw[-] (7,0) -- (7,-1);
\draw[-] (7,0) -- (8,0);
\draw[-] (7,-1) -- (8,-1);
\draw[-] (8,0) -- (8,-1);
\draw[-] (8,0) -- (9,0);
\draw[-] (8,-1) -- (9,-1);
\draw[-] (9,0) -- (9,-1);
\draw[-] (8,1) -- (9,1);
\draw[-] (9,0) -- (9,1);
\draw[-] (8,0) -- (8,1);
\draw[dashed] (-1,0) -- (0,-1);
\draw[dashed] (1,0) -- (2,-1);
\draw[dashed] (0,0) -- (1,-1);
\draw[dashed] (2,0) -- (3,-1);
\draw[dashed] (3,0) -- (4,-1);
\draw[dashed] (5,0) -- (6,-1);
\draw[dashed] (6,0) -- (7,-1);
\draw[dashed] (7,0) -- (8,-1);
\draw[dashed] (8,0) -- (9,-1);
\draw[dashed] (8,1) -- (9,0);
\node[left] at (-0.95,-0.5) {$a_2$};
\node[left] at (2.25,-0.5) {$a_4$};
\node[left] at (4.25,-0.5) {$a_4$};
\node[left] at (5.25,-0.5) {$a_4$};
\node[left] at (7.25,-0.5) {$a_4$};
\node[left] at (1.25,-0.5) {$a_2$};
\node[left] at (0.25,-0.5) {$a_4$};
\node[left] at (3.25,-0.5) {$a_2$};
\node[left] at (6.25,-0.5) {$a_2$};
\node[left] at (8.25,-0.5) {$a_2$};
\node[right] at (8.95,-0.5) {$a_4$};
\node[above] at (0.5,-0.05) {$\tau'$};
\node[above] at (1.5,-0.05) {$a_1$};
\node[above] at (-0.5,-0.05) {$a_1$};
\node[above] at (2.5,-0.05) {$\tau'$};
\node[above] at (3.5,-0.05) {$a_1$};
\node[above] at (5.5,-0.05) {$\tau'$};
\node[above] at (6.5,-0.05) {$a_1$};
\node[above] at (8.5,-0.2) {$a_1$};
\node[above] at (7.5,-0.05) {$\tau'$};
\node[below] at (0.5,-0.95) {$\tau'$};
\node[below] at (1.5,-0.95) {$a_1$};
\node[below] at (-0.5,-0.95) {$a_1$};
\node[below] at (2.5,-0.95) {$\tau'$};
\node[below] at (3.5,-0.95) {$a_1$};
\node[below] at (5.5,-0.95) {$\tau'$};
\node[below] at (6.5,-0.95) {$a_1$};
\node[below] at (8.5,-0.95) {$a_1$};
\node[below] at (7.5,-0.95) {$\tau'$};
\node[right] at (1.3,-0.45) {$\tau'$};
\node[right] at (-0.7,-0.45) {$\tau'$};
\node[right] at (3.3,-0.45) {$\tau'$};
\node[right] at (2.3,-0.5) {$a_1$};
\node[right] at (6.3,-0.45) {$\tau'$};
\node[right] at (8.3,-0.45) {$\tau'$};
\node[right] at (8.3,0.55) {$a_4$};
\node[right] at (5.3,-0.5) {$a_1$};
\node[right] at (7.3,-0.5) {$a_1$};
\node[left] at (0.8,-0.5) {$a_1$};
\node[left] at (8.05,0.5) {$\tau'$};
\node[right] at (8.95,0.5) {$a_8$};
\node[above] at (8.5,0.95) {$a_7$};
\draw [-] (1, 0) -- (1,-1);
\draw [-] (2, 0) -- (2,-1);
\draw [fill] (0,0) circle [radius=.05];
\draw [fill] (0,-1) circle [radius=.05];
\draw [fill] (-1,0) circle [radius=.05];
\draw [fill] (-1,-1) circle [radius=.05];
\draw [fill] (1,0) circle [radius=.05];
\draw [fill] (2,0) circle [radius=.05];
\draw [fill] (1,-1) circle [radius=.05];
\draw [fill] (2,-1) circle [radius=.05];
\draw [fill] (3,0) circle [radius=.05];
\draw [fill] (3,-1) circle [radius=.05];
\draw [fill] (4,0) circle [radius=.05];
\draw [fill] (4,-1) circle [radius=.05];
\draw [fill] (5,0) circle [radius=.05];
\draw [fill] (5,-1) circle [radius=.05];
\draw [fill] (6,0) circle [radius=.05];
\draw [fill] (6,-1) circle [radius=.05];
\draw [fill] (7,0) circle [radius=.05];
\draw [fill] (8,0) circle [radius=.05];
\draw [fill] (8,-1) circle [radius=.05];
\draw [fill] (9,0) circle [radius=.05];
\draw [fill] (9,-1) circle [radius=.05];
\draw [fill] (9,1) circle [radius=.05];
\draw [fill] (8,1) circle [radius=.05];
\draw [fill] (4.3,-0.5) circle [radius=.02];
\draw [fill] (4.5,-0.5) circle [radius=.02];
\draw [fill] (4.7,-0.5) circle [radius=.02];
\draw [fill] (7,-1) circle [radius=.05];
\end{tikzpicture}}

We clearly have $n^{\tau}(T,\rho)=n^{\tau'}(T',\rho)$.

\medskip

As $G_{T,\rho}$ can be obtained by gluing a tile with a graph which is isomorphic to $H_s(a_1,\tau)$, we have $\mathcal P(G_{T,\rho})=\mathcal P_1\sqcup \mathcal P_2$, where $\mathcal P_1$ contains all $P$ which contains the upper right edge labeled $a_7$, $\mathcal P_2$ contains all $P$ which contains the upper right edge labeled $a_8$. By Lemma \ref{basicdecom}, as sets, we have $$\mathcal P_1\cong\textstyle\bigsqcup_{(\lambda_1,\cdots,\lambda_{s+1})\in \{0,1\}^{s+1}}\mathcal H'_{(\lambda_1,\cdots,\lambda_{s+1})}(a_1,\tau),$$ $$\mathcal P_2\cong\textstyle\bigsqcup_{(\lambda_1,\cdots,\lambda_{s})\in \{0,1\}^{s}}\mathcal H'_{(\lambda_1,\cdots,\lambda_{s})}(a_1,\tau).$$

\medskip

Similarly, $\mathcal P(G_{T',\rho})=\mathcal P'_1\sqcup \mathcal P'_2$, where $\mathcal P'_1$ contains all $P$ which contains the upper right edge labeled $a_7$, $\mathcal P'_2$ contains all $P$ which contains the upper right edge labeled $a_8$. As sets, we have $$\mathcal P'_1\cong\textstyle\bigsqcup_{(\lambda_1,\cdots,\lambda_{s+1})\in \{0,1\}^{s+1}}\mathcal H_{(\lambda_1,\cdots,\lambda_{s+1})}(\tau',a_1),$$ $$\mathcal P'_2\cong\textstyle\bigsqcup_{(\lambda_1,\cdots,\lambda_{s})\in \{0,1\}^{s}}\mathcal H_{(\lambda_1,\cdots,\lambda_{s})}(\tau',a_1).$$

\medskip

As case (6), under the above isomorphisms, $\psi_{\rho}: \mathcal H'_{\lambda}(a_1,\tau)\rightarrow \mathcal H_{1-\lambda}(\tau',a_1)$ for $\lambda=(\lambda_1,\cdots,\lambda_{s+1})$ or $(\lambda_1,\cdots,\lambda_{s})$ gives a partition bijection between $\mathcal P(G_{T,\rho})$ and $\mathcal P(G_{T',\rho})$. For each $P\in \mathcal H'_{\lambda}(a_1,\tau)$ and $P'\in \mathcal H_{1-\lambda}(\tau',a_1)$, by Lemma \ref{same-num} the number of non-$\tau$-mutable edges labeled $a$ in $P$ equals to the number of non-$\tau'$-mutable edges labeled $a$ in $P'$ for each $\tau\neq a\in T$. Moreover, the last gluing edge of $G_{T,\rho}$ is in $P$ if and only if the last gluing edge of $G_{T',\rho}$ in $P'$. Since $\rho$ and $\zeta$ have the same starting point, $G_{T,\rho}$ and $G_{T',\rho}$ do not have first gluing edge.

\medskip

In case (8), up to a difference of relative orientation, $G_{T,\rho}$ and $G_{T',\rho}$ are the following graphs, respectively,

\centerline{\begin{tikzpicture}
\draw[-] (0,-1) -- (0,-2);
\draw[-] (1,-1) -- (1,-2);
\draw[-] (0,-2) -- (1,-2);
\draw[-] (0,0) -- (1,0);
\draw[-] (0,0) -- (0,-1);
\draw[-] (0,-1) -- (1,-1);
\draw[-] (1,0) -- (2,0);
\draw[-] (1,-1) -- (2,-1);
\draw[-] (2,0) -- (3,0);
\draw[-] (2,-1) -- (3,-1);
\draw[-] (3,0) -- (3,-1);
\draw[-] (3,0) -- (4,0);
\draw[-] (3,-1) -- (4,-1);
\draw[-] (4,0) -- (4,-1);
\draw[-] (5,0) -- (6,0);
\draw[-] (5,0) -- (5,-1);
\draw[-] (6,0) -- (6,-1);
\draw[-] (5,-1) -- (6,-1);
\draw[-] (6,0) -- (7,0);
\draw[-] (6,-1) -- (7,-1);
\draw[-] (7,0) -- (7,-1);
\draw[-] (7,0) -- (8,0);
\draw[-] (7,-1) -- (8,-1);
\draw[-] (8,0) -- (8,-1);
\draw[-] (7,1) -- (8,1);
\draw[-] (7,0) -- (7,1);
\draw[-] (8,0) -- (8,1);
\draw[dashed] (1,0) -- (2,-1);
\draw[dashed] (0,0) -- (1,-1);
\draw[dashed] (0,-1) -- (1,-2);
\draw[dashed] (2,0) -- (3,-1);
\draw[dashed] (3,0) -- (4,-1);
\draw[dashed] (5,0) -- (6,-1);
\draw[dashed] (6,0) -- (7,-1);
\draw[dashed] (7,0) -- (8,-1);
\draw[dashed] (7,1) -- (8,0);
\node[left] at (0.05,-0.5) {$a_2$};
\node[left] at (2.25,-0.5) {$a_2$};
\node[left] at (4.25,-0.5) {$a_2$};
\node[left] at (5.25,-0.5) {$a_2$};
\node[left] at (7.25,-0.5) {$a_2$};
\node[left] at (1.25,-0.5) {$a_4$};
\node[left] at (3.25,-0.5) {$a_4$};
\node[left] at (6.25,-0.5) {$a_4$};
\node[right] at (7.95,-0.5) {$a_4$};
\node[above] at (0.5,-0.05) {$\tau$};
\node[above] at (1.5,-0.05) {$a_1$};
\node[above] at (2.5,-0.05) {$\tau$};
\node[above] at (3.5,-0.05) {$a_1$};
\node[above] at (5.5,-0.05) {$\tau$};
\node[above] at (6.5,-0.05) {$a_1$};
\node[above] at (7.5,-0.2) {$\tau$};
\node[below] at (0.5,-0.8) {$\tau$};
\node[below] at (1.5,-0.95) {$a_1$};
\node[below] at (2.5,-0.95) {$\tau$};
\node[below] at (3.5,-0.95) {$a_1$};
\node[below] at (0.5,-1.95) {$a_6$};
\node[below] at (5.5,-0.95) {$\tau$};
\node[below] at (6.5,-0.95) {$a_1$};
\node[below] at (7.5,-0.95) {$\tau$};
\node[right] at (1.3,-0.45) {$\tau$};
\node[right] at (3.3,-0.45) {$\tau$};
\node[right] at (2.3,-0.5) {$a_1$};
\node[right] at (7.3,0.5) {$a_4$};
\node[right] at (0.95,-1.5) {$a_1$};
\node[right] at (6.3,-0.45) {$\tau$};
\node[right] at (5.3,-0.5) {$a_1$};
\node[right] at (7.3,-0.5) {$a_1$};
\node[left] at (0.8,-0.5) {$a_1$};
\node[left] at (0.05,-1.5) {$a_5$};
\node[left] at (0.8,-1.5) {$a_2$};
\node[left] at (7.05,0.5) {$a_1$};
\node[right] at (7.95,0.5) {$a_8$};
\node[above] at (7.5,0.95) {$a_7$};
\draw [-] (1, 0) -- (1,-1);
\draw [-] (2, 0) -- (2,-1);
\draw [fill] (0,0) circle [radius=.05];
\draw [fill] (0,-1) circle [radius=.05];
\draw [fill] (0,-2) circle [radius=.05];
\draw [fill] (1,-2) circle [radius=.05];
\draw [fill] (1,0) circle [radius=.05];
\draw [fill] (2,0) circle [radius=.05];
\draw [fill] (1,-1) circle [radius=.05];
\draw [fill] (2,-1) circle [radius=.05];
\draw [fill] (3,0) circle [radius=.05];
\draw [fill] (3,-1) circle [radius=.05];
\draw [fill] (4,0) circle [radius=.05];
\draw [fill] (4,-1) circle [radius=.05];
\draw [fill] (5,0) circle [radius=.05];
\draw [fill] (5,-1) circle [radius=.05];
\draw [fill] (6,0) circle [radius=.05];
\draw [fill] (6,-1) circle [radius=.05];
\draw [fill] (7,0) circle [radius=.05];
\draw [fill] (8,0) circle [radius=.05];
\draw [fill] (8,-1) circle [radius=.05];
\draw [fill] (7,1) circle [radius=.05];
\draw [fill] (8,1) circle [radius=.05];
\draw [fill] (4.3,-0.5) circle [radius=.02];
\draw [fill] (4.5,-0.5) circle [radius=.02];
\draw [fill] (4.7,-0.5) circle [radius=.02];
\draw [fill] (7,-1) circle [radius=.05];
\end{tikzpicture}}

\centerline{\begin{tikzpicture}
\draw[-] (-1,0) -- (0,0);
\draw[-] (-1,-1) -- (0,-1);
\draw[-] (-1,-2) -- (0,-2);
\draw[-] (-1,-1) -- (-1,-2);
\draw[-] (0,-1) -- (0,-2);
\draw[-] (-1,0) -- (-1,-1);
\draw[-] (0,0) -- (1,0);
\draw[-] (0,0) -- (0,-1);
\draw[-] (0,-1) -- (1,-1);
\draw[-] (1,0) -- (2,0);
\draw[-] (1,-1) -- (2,-1);
\draw[-] (2,0) -- (3,0);
\draw[-] (2,-1) -- (3,-1);
\draw[-] (3,0) -- (3,-1);
\draw[-] (3,0) -- (4,0);
\draw[-] (3,-1) -- (4,-1);
\draw[-] (4,0) -- (4,-1);
\draw[-] (5,0) -- (6,0);
\draw[-] (5,0) -- (5,-1);
\draw[-] (6,0) -- (6,-1);
\draw[-] (5,-1) -- (6,-1);
\draw[-] (6,0) -- (7,0);
\draw[-] (6,-1) -- (7,-1);
\draw[-] (7,0) -- (7,-1);
\draw[-] (7,0) -- (8,0);
\draw[-] (7,-1) -- (8,-1);
\draw[-] (8,0) -- (8,-1);
\draw[-] (8,0) -- (9,0);
\draw[-] (8,-1) -- (9,-1);
\draw[-] (9,0) -- (9,-1);
\draw[-] (8,1) -- (9,1);
\draw[-] (9,0) -- (9,1);
\draw[-] (8,0) -- (8,1);
\draw[dashed] (-1,0) -- (0,-1);
\draw[dashed] (-1,-1) -- (0,-2);
\draw[dashed] (1,0) -- (2,-1);
\draw[dashed] (0,0) -- (1,-1);
\draw[dashed] (2,0) -- (3,-1);
\draw[dashed] (3,0) -- (4,-1);
\draw[dashed] (5,0) -- (6,-1);
\draw[dashed] (6,0) -- (7,-1);
\draw[dashed] (7,0) -- (8,-1);
\draw[dashed] (8,0) -- (9,-1);
\draw[dashed] (8,1) -- (9,0);
\node[left] at (-0.95,-0.5) {$a_2$};
\node[left] at (2.25,-0.5) {$a_4$};
\node[left] at (4.25,-0.5) {$a_4$};
\node[left] at (5.25,-0.5) {$a_4$};
\node[left] at (7.25,-0.5) {$a_4$};
\node[left] at (1.25,-0.5) {$a_2$};
\node[left] at (0.25,-0.5) {$a_4$};
\node[left] at (3.25,-0.5) {$a_2$};
\node[left] at (6.25,-0.5) {$a_2$};
\node[left] at (8.25,-0.5) {$a_2$};
\node[right] at (8.95,-0.5) {$a_4$};
\node[above] at (0.5,-0.05) {$\tau'$};
\node[above] at (1.5,-0.05) {$a_1$};
\node[above] at (-0.5,-0.05) {$a_1$};
\node[above] at (2.5,-0.05) {$\tau'$};
\node[above] at (3.5,-0.05) {$a_1$};
\node[above] at (5.5,-0.05) {$\tau'$};
\node[above] at (6.5,-0.05) {$a_1$};
\node[above] at (8.5,-0.2) {$a_1$};
\node[above] at (7.5,-0.05) {$\tau'$};
\node[below] at (0.5,-0.95) {$\tau'$};
\node[below] at (1.5,-0.95) {$a_1$};
\node[below] at (-0.5,-0.8) {$a_1$};
\node[below] at (2.5,-0.95) {$\tau'$};
\node[below] at (3.5,-0.95) {$a_1$};
\node[below] at (5.5,-0.95) {$\tau'$};
\node[below] at (6.5,-0.95) {$a_1$};
\node[below] at (8.5,-0.95) {$a_1$};
\node[below] at (7.5,-0.95) {$\tau'$};
\node[below] at (-0.5,-1.95) {$a_6$};
\node[right] at (1.3,-0.45) {$\tau'$};
\node[right] at (-0.05,-1.5) {$\tau'$};
\node[right] at (-0.7,-0.45) {$\tau'$};
\node[right] at (-0.8,-1.45) {$a_2$};
\node[right] at (3.3,-0.45) {$\tau'$};
\node[right] at (2.3,-0.5) {$a_1$};
\node[right] at (6.3,-0.45) {$\tau'$};
\node[right] at (8.3,-0.45) {$\tau'$};
\node[right] at (8.3,0.55) {$a_4$};
\node[right] at (5.3,-0.5) {$a_1$};
\node[right] at (7.3,-0.5) {$a_1$};
\node[left] at (0.8,-0.5) {$a_1$};
\node[left] at (-0.95,-1.5) {$a_5$};
\node[left] at (8.05,0.5) {$\tau'$};
\node[right] at (8.95,0.5) {$a_8$};
\node[above] at (8.5,0.95) {$a_7$};
\draw [-] (1, 0) -- (1,-1);
\draw [-] (2, 0) -- (2,-1);
\draw [fill] (0,0) circle [radius=.05];
\draw [fill] (0,-1) circle [radius=.05];
\draw [fill] (0,-2) circle [radius=.05];
\draw [fill] (-1,-2) circle [radius=.05];
\draw [fill] (-1,0) circle [radius=.05];
\draw [fill] (-1,-1) circle [radius=.05];
\draw [fill] (1,0) circle [radius=.05];
\draw [fill] (2,0) circle [radius=.05];
\draw [fill] (1,-1) circle [radius=.05];
\draw [fill] (2,-1) circle [radius=.05];
\draw [fill] (3,0) circle [radius=.05];
\draw [fill] (3,-1) circle [radius=.05];
\draw [fill] (4,0) circle [radius=.05];
\draw [fill] (4,-1) circle [radius=.05];
\draw [fill] (5,0) circle [radius=.05];
\draw [fill] (5,-1) circle [radius=.05];
\draw [fill] (6,0) circle [radius=.05];
\draw [fill] (6,-1) circle [radius=.05];
\draw [fill] (7,0) circle [radius=.05];
\draw [fill] (8,0) circle [radius=.05];
\draw [fill] (8,-1) circle [radius=.05];
\draw [fill] (9,0) circle [radius=.05];
\draw [fill] (9,-1) circle [radius=.05];
\draw [fill] (9,1) circle [radius=.05];
\draw [fill] (8,1) circle [radius=.05];
\draw [fill] (4.3,-0.5) circle [radius=.02];
\draw [fill] (4.5,-0.5) circle [radius=.02];
\draw [fill] (4.7,-0.5) circle [radius=.02];
\draw [fill] (7,-1) circle [radius=.05];
\end{tikzpicture}}

We clearly have $n^{\tau}(T,\rho)=n^{\tau'}(T',\rho)$.

\medskip

Since $G_{T,\rho}$ can be obtained  by gluing two tiles with a graph which is isomorphic to $H_s(a_1,\tau)$, we have $\mathcal P(G_{T,\rho})=\mathcal P_1\sqcup \mathcal P_2\sqcup \mathcal P_3\sqcup \mathcal P_4$, where $\mathcal P_1$ contains all $P$ which contains the lower left edge labeled $a_6$ and the upper right edge labeled $a_7$, $\mathcal P_2$ contains all $P$ which contains the lower left edge labeled $a_6$ and the upper right edge labeled $a_8$, $\mathcal P_3$ contains all $P$ which contains the lower left edge labeled $a_5$ and the upper right edge labeled $a_7$, $\mathcal P_4$ contains all $P$ which contains the lower left edge labeled $a_5$ and the upper right edge labeled $a_8$. By Lemma \ref{basicdecom}, as sets, we have
$$\mathcal P_1\cong\textstyle\bigsqcup_{(\lambda_1,\cdots,\lambda_{s+1})\in \{0,1\}^{s+1}}\mathcal H'_{(\lambda_1,\cdots,\lambda_{s+1})}(a_1,\tau),$$ $$\mathcal P_2\cong\textstyle\bigsqcup_{(\lambda_1,\cdots,\lambda_{s})\in \{0,1\}^{s}}\mathcal H'_{(\lambda_1,\cdots,\lambda_{s})}(a_1,\tau),$$ $$\mathcal P_3\cong\textstyle\bigsqcup_{(\lambda_2,\cdots,\lambda_{s+1})\in \{0,1\}^{s}}\mathcal H'_{(\lambda_1,\cdots,\lambda_{s+1})}(a_1,\tau),$$ $$\mathcal P_4\cong\textstyle\bigsqcup_{(\lambda_2,\cdots,\lambda_{s})\in \{0,1\}^{s-1}}\mathcal H'_{(\lambda_2,\cdots,\lambda_{s})}(a_1,\tau).$$

Similarly, $\mathcal P(G_{T',\rho})=\mathcal P'_1\sqcup \mathcal P'_2\sqcup \mathcal P'_3\sqcup \mathcal P'_4$, where $\mathcal P'_1$ contains all $P$ which contains the lower left edge labeled $a_6$ and the upper right edge labeled $a_7$, $\mathcal P'_2$ contains all $P$ which contains the lower left edge labeled $a_6$ and the upper right edge labeled $a_8$, $\mathcal P'_3$ contains all $P$ which contains the lower left edge labeled $a_5$ and the upper right edge labeled $a_7$, $\mathcal P'_4$ contains all $P$ which contains the lower left edge labeled $a_5$ and the upper right edge labeled $a_8$. As sets, we have $$\mathcal P'_1\cong\textstyle\bigsqcup_{(\lambda_1,\cdots,\lambda_{s+1})\in \{0,1\}^{s+1}}\mathcal H_{(\lambda_1,\cdots,\lambda_{s+1})}(\tau',a_1),$$ $$\mathcal P'_2\cong\textstyle\bigsqcup_{(\lambda_1,\cdots,\lambda_s)\in \{0,1\}^{s}}\mathcal H_{(\lambda_1,\cdots,\lambda_s)}(\tau',a_1),$$ $$\mathcal P'_3\cong\textstyle\bigsqcup_{(\lambda_2,\cdots,\lambda_{s+1})\in \{0,1\}^{s}}\mathcal H_{(\lambda_2,\cdots,\lambda_{s+1})}(\tau',a_1),$$ $$\mathcal P'_4\cong\textstyle\bigsqcup_{(\lambda_2,\cdots,\lambda_{s})\in \{0,1\}^{s-1}}\mathcal H_{(\lambda_2,\cdots,\lambda_{s})}(\tau',a_1).$$

\medskip

As cases (6,7), under the above isomorphisms, $\psi_{\rho}: \mathcal H'_{\lambda}(a_1,\tau)\rightarrow \mathcal H_{1-\lambda}(\tau',a_1)$ for $\lambda=(\lambda_1,\cdots,\lambda_{s+1})$, $(\lambda_1,\cdots,\lambda_{s})$, $(\lambda_2,\cdots,\lambda_{s+1})$ or $(\lambda_2,\cdots,\lambda_{s})$ gives a partition bijection between $\mathcal P(G_{T,\rho})$ and $\mathcal P(G_{T',\rho})$. For each $P\in \mathcal H'_{\lambda}(a_1,\tau)$ and $P'\in \mathcal H_{1-\lambda}(\tau',a_1)$, by Lemma \ref{same-num} the number of non-$\tau$-mutable edges labeled $a$ in $P$ equals to the number of non-$\tau'$-mutable edges labeled $a$ in $P'$ for each $\tau\neq a\in T$. Moreover, the first/last gluing edge of $G_{T,\rho}$ is in $P$ if and only if the first/last gluing edge of $G_{T',\rho}$ in $P'$.

\medskip

\subsection{Partition bijection between $\mathcal P(G_{T,\rho})$ and $\mathcal P(G_{T',\rho})$} Herein, we summarize the discussion of Subsections \ref{local1} and \ref{local2}. We recall that to give a partition bijection $\psi_{\rho}$ from $\mathcal P(G_{T,\rho})$ to $\mathcal P(G_{T',\rho})$ is equivalent to associate each $P\in \mathcal P(G_{T,\rho})$ with a non-empty subset $\psi_{\rho}(P)\subset \mathcal P(G_{T',\rho})$ such that $\psi_{\rho}(P_1)=\psi_{\rho}(P_2)$ if $\psi_{\rho}(P_1)\cap \psi_{\rho}(P_2)\neq \emptyset$ for any $P_1,P_2\in \mathcal P(G_{T,\rho})$ and $\cup_{P\in \mathcal P(G_{T,\rho})}\psi_{\rho}(P)=\mathcal P(G_{T',\rho})$.

\medskip

\begin{Lemma}\label{partition3}

To give a partition bijection $\psi_{\rho}$ from $\mathcal P(G_{T,\rho})$ to $\mathcal P(G_{T',\rho})$ and its inverse $\psi'_{\rho}$ is equivalent to associate  each $P\in \mathcal P(G_{T,\rho})$ with a non-empty subset $\psi_{\rho}(P)\subset \mathcal P(G_{T',\rho})$ and each $P'\in \mathcal P(G_{T',\rho})$ with a non-empty subset $\psi'_{\rho}(P')\subset \mathcal P(G_{T,\rho})$, which satisfy the following

\begin{enumerate}[$(1)$]

  \item $\psi_{\rho}(P_1)=\psi_{\rho}(P_2)$ if $\psi_{\rho}(P_1)\cap \psi_{\rho}(P_2)\neq \emptyset$ for any $P_1,P_2\in \mathcal P(G_{T,\rho})$.

  \item $\psi'_{\rho}(P'_1)=\psi'_{\rho}(P'_2)$ if $\psi'_{\rho}(P'_1)\cap \psi'_{\rho}(P'_2)\neq \emptyset$ for any $P'_1,P'_2\in \mathcal P(G_{T',\rho})$.

  \item $P'\in \psi_{\rho}(P)$ if and only if $P\in \psi'_{\rho}(P')$.

\end{enumerate}

\end{Lemma}

\begin{proof}

(1) and (3) imply that $\psi_{\rho}$ is a partition bijection. (2) and (3) imply that $\psi'_{\rho}$ is a partition bijection. Furthermore, (3) shows $\psi_{\rho}$ and $\psi'_{\rho}$ are inverse to each other.
\end{proof}

\medskip

\begin{Proposition}\label{compare-loc}

With the same notation as above, let $T$ be a triangulation of $\mathcal O$ and $\tau$ be an arc in $T$. We assume that $\tau'\neq\zeta\notin T$ is an arc in $\mathcal O$. Let $\rho$ be a subcurve of $\zeta$ satisfying the assumptions at the beginning of this section. Then

\begin{enumerate}[$(1)$]

  \item $n^{\tau}(T,\rho)=n^{\tau'}(T',\rho)$.

  \item There is a partition bijection $\psi_{\rho}$ from $\mathcal P(G_{T},\rho)$ to $\mathcal P(G_{T'},\rho)$ and its inverse $\psi'_{\rho}$, which satisfy for any $P\in \mathcal P(G_{T},\rho)$ and $P'\in \psi_{\rho}(P)$,

      \begin{enumerate}[$(a)$]

        \item $P\in \mathcal P^{\tau}_\nu(G_{T,\rho})$ for some $\nu$ if and only if $P'\in \mathcal P^{\tau'}_{-\nu}(G_{T',\rho})$;

        \item the number of non-$\tau$-mutable edges labeled $a$ in $P$ equals to the number of non-$\tau'$-mutable edges labeled $a$ in $P'$ for each $\tau\neq a\in T$;

        \item the first/last gluing edge of $G_{T,\rho}$ is in $P$ if and only if the first/last gluing edge of $G_{T',\rho}$ is in $P'$.

     \end{enumerate}

\end{enumerate}

\end{Proposition}

\begin{proof}

We already proved the proposition in case $\mathcal O_0$ contains no orbifold points. It can be seen that Lemmas \ref{formglobal}, \ref{diff2} hold for the unpunctured surface case.

Let $\widetilde{\mathcal O}_0$ be the lift of $\mathcal O_0$ given in Subsection \ref{local1} and let $\widetilde \rho$ be a curve in $\widetilde{\mathcal O}_0$ lift by $\rho$. Using the strategy at the beginning of this section, we divide $\widetilde \rho$ into subcurves $\widetilde \rho_1,\cdots, \widetilde \rho_s$ for some $s$. There are two possibilities: (I) $\widetilde{\mathcal O}_0$ contains no orbifold point, and (II) $\widetilde{\mathcal O}_0$ contains orbifold points.

(I) $\widetilde{\mathcal O}_0$ contains no orbifold point.

Given an arc $\widetilde\tau\in \widetilde T_0$ with $\chi(\widetilde\tau)=\tau$, for each $i\in [1,s]$, there is a partition bijection $\psi_{\widetilde \rho_i}$ from $\mathcal P(G_{\widetilde T_0,\widetilde \rho_i})$ to $\mathcal P(G_{\mu_{\widetilde\tau}\widetilde T_0,\widetilde \rho_i})$. For each $\widetilde P\in \mathcal P(G_{\widetilde T_0,\widetilde \rho})$, by Proposition \ref{decompose}, we may write $\widetilde P$ as $(\widetilde P_i)$ for $\widetilde P_i\in \mathcal P(G_{\widetilde T_0,\widetilde \rho_i})$. According to Lemma \ref{formglobal} in the unpunctured surface case, $(\widetilde P'_i)$ is a perfect matching in $\mathcal P(G_{\mu_{\widetilde\tau}\widetilde T_0,\widetilde \rho})$ for $\widetilde P'_i\in \psi_{\widetilde \rho_i}(\widetilde P_i)$. Thus, we can associate $\widetilde P$ with a subset
$$\psi_{\widetilde \tau,\widetilde \rho}(\widetilde P)=\{(\widetilde P'_i)\mid \widetilde P'_i\in \psi_{\widetilde \rho_i}(\widetilde P_i)\}\subset \mathcal P(G_{\mu_{\widetilde\tau}\widetilde T_0,\widetilde \rho}).$$
For any $\widetilde P=(\widetilde P_i), \widetilde Q=(\widetilde Q_i)\in \mathcal P(G_{\widetilde T_0,\widetilde \rho})$, if $\psi_{\widetilde \tau,\widetilde \rho}(\widetilde P)\cap \psi_{\widetilde \tau,\widetilde \rho}(\widetilde Q)\neq \emptyset$, then $\psi_{\widetilde \tau,\widetilde \rho_i}(\widetilde P_i)\cap \psi_{\widetilde \tau,\widetilde \rho_i}(\widetilde Q_i)\neq \emptyset$ for each $i$. Therefore, $\psi_{\widetilde \tau,\widetilde \rho_i}(\widetilde P_i)= \psi_{\widetilde \tau,\widetilde \rho_i}(\widetilde Q_i)$ for each $i$. Consequently, $\psi_{\widetilde \tau,\widetilde \rho}(\widetilde P)=\psi_{\widetilde \tau,\widetilde \rho}(\widetilde Q)$.

Dually, we can associate $\widetilde P'\in \mathcal P(G_{\mu_{\widetilde\tau}\widetilde T_0,\widetilde \rho})$ with a subset
$$\psi'_{\widetilde \tau,\widetilde \rho}(\widetilde P')=\{(\widetilde P_i)\mid \widetilde P_i\in \psi'_{\widetilde \rho_i}(\widetilde P'_i)\}\subset \mathcal P(G_{\widetilde T_0,\widetilde \rho}).$$
Moreover, $\psi'_{\widetilde \tau,\widetilde \rho}(\widetilde P')=\psi'_{\widetilde \tau,\widetilde \rho}(\widetilde Q')$ if $\psi'_{\widetilde \tau,\widetilde \rho}(\widetilde P')\cap \psi'_{\widetilde \tau,\widetilde \rho}(\widetilde Q')\neq \emptyset$ for any $\widetilde P',\widetilde Q'\in \mathcal P(G_{\mu_{\widetilde\tau}\widetilde T_0,\widetilde \rho})$.

Since $\psi_{\widetilde \rho_i}$ and $\psi'_{\widetilde \rho_i}$ are partition bijections inverse to each other, using Lemma \ref{partition3}, we have $\widetilde P\in \psi'_{\widetilde \tau,\widetilde \rho}(\widetilde P')$ if and only if $\widetilde P'\in \psi_{\widetilde \tau,\widetilde \rho}(\widetilde P)$. Therefore, by Lemma \ref{partition3}, $\psi_{\widetilde \tau,\widetilde \rho}$ and $\psi'_{\widetilde \tau,\widetilde \rho}$ are partition bijections inverse to each other. Moreover, for any $i$, since conditions (b) and (c) hold for $\psi_{\widetilde \tau,\widetilde\rho_i}$ and $\psi'_{\widetilde \tau,\widetilde\rho_i}$, we have the conditions (b) and (c) hold for $\psi_{\widetilde \tau,\widetilde\rho}$ and $\psi'_{\widetilde \tau,\widetilde\rho}$.

We recall that there are at most two arcs $\widetilde{\tau}$ in $\widetilde{\mathcal O}_0$ such that $\chi(\widetilde{\tau})=\tau$.

If such $\widetilde \tau$ uniquely exists, then $\mu_{\widetilde\tau}\widetilde T_0$ is the triangulation of $\widetilde{\mathcal O}_0$ lift by $\mu_{\tau}T_0$. Therefore, the partition maps $\psi_{\widetilde \tau,\widetilde\rho}$ and $\psi'_{\widetilde \tau,\widetilde\rho}$ induce the required partition maps $\psi_{\rho}$ and $\psi'_{\rho}$. Since $n^{\widetilde\tau}(\widetilde T_0,\widetilde \rho_i)=n^{\widetilde \tau'}(\mu_{\widetilde\tau}\widetilde T_0,\widetilde\rho_i)$ for each $i$, we obtain $n^{\tau}(T,\rho)=n^{\tau'}(T',\rho)$.

If there are $\widetilde{\tau}_1$ and $\widetilde{\tau}_2$ such that $\chi(\widetilde{\tau}_j)=\tau$ for $j=1,2$, then  $\mu_{\widetilde\tau_2}\mu_{\widetilde\tau_1}\widetilde T_0$ is the triangulation of $\widetilde{\mathcal O}_0$ lift by $\mu_{\tau}T_0$. Thus, we have partition bijections $\psi_{\widetilde \tau_1,\widetilde\rho}$ and $\psi'_{\widetilde \tau_1,\widetilde\rho}$ between $\mathcal P(G_{\widetilde T_0,\widetilde \rho})$ and $\mathcal P(G_{\mu_{\widetilde\tau_1}\widetilde T_0,\widetilde \rho})$, and partition bijections $\psi_{\widetilde \tau_2,\widetilde\rho}$ and $\psi'_{\widetilde \tau_2,\widetilde\rho}$ between $\mathcal P(G_{\mu_{\widetilde\tau_1}\widetilde T_0,\widetilde \rho})$ and $\mathcal P(G_{\mu_{\widetilde \tau_2}\mu_{\widetilde\tau_1}\widetilde T_0,\widetilde \rho})$. For each $\widetilde P\in \mathcal P(G_{\widetilde T_0,\widetilde \rho})$, we associate with
$$\psi_{\widetilde \rho}(\widetilde P)=\{\widetilde P''\mid \widetilde P''\in \psi_{\widetilde \tau_2,\widetilde\rho}(\widetilde P')\text{ for some }\widetilde P'\in \psi_{\widetilde \tau_1,\widetilde\rho}(\widetilde P)\}\subset \mathcal P(G_{\mu_{\widetilde \tau_2}\mu_{\widetilde\tau_1}\widetilde T_0,\widetilde \rho}).$$
Dually, for each $\widetilde P''\in \mathcal P(G_{\mu_{\widetilde \tau_2}\mu_{\widetilde\tau_1}\widetilde T_0,\widetilde \rho})$, we associate with
$$\psi'_{\widetilde \rho}(\widetilde P'')=\{\widetilde P\mid \widetilde P\in \psi'_{\widetilde \tau_1,\widetilde\rho}(\widetilde P')\text{ for some }\widetilde P'\in \psi'_{\widetilde \tau_2,\widetilde\rho}(\widetilde P'')\}\subset \mathcal P(G_{\widetilde T_0,\widetilde \rho}).$$

We claim that $\psi_{\widetilde \rho}$ and $\psi'_{\widetilde \rho}$ are partition bijections inverse to each other.

If $\psi'_{\widetilde \rho}(\widetilde P'')\cap \psi'_{\widetilde \rho}(\widetilde Q'')\neq \emptyset$ for some $\widetilde P'',\widetilde Q''\in \mathcal P(G_{\mu_{\widetilde \tau_2}\mu_{\widetilde\tau_1}\widetilde T_0,\widetilde \rho})$, we may assume $\widetilde P\in \psi'_{\widetilde \rho}(\widetilde P'')\cap \psi'_{\widetilde \rho}(\widetilde Q'')$, then $\widetilde P\in \psi'_{\widetilde \tau_1,\widetilde\rho}(\widetilde P')\cap \psi'_{\widetilde \tau_1,\widetilde\rho}(\widetilde Q')$ for some $\widetilde P',\widetilde Q'$ which satisfy $\widetilde P'\in \psi'_{\widetilde \tau_2,\widetilde\rho}(\widetilde P'')$ and $\widetilde Q'\in \psi'_{\widetilde \tau_2,\widetilde\rho}(\widetilde Q'')$. By Lemma \ref{diff2} in the unpunctured surface case, $\widetilde Q'$ can be obtained from $\widetilde P'$ by twists on some tiles with diagonals labeled $\widetilde \tau_1$, we may assume $\widetilde Q'=\prod_{i\in I}\mu_{p_i}\widetilde P'$ for some $I$ and $p_i$. Then $\widetilde P'=\prod_{i\in I}\mu_{p_i}\widetilde Q'$. We assume that $p_i$ corresponds to the tile $G(p_i)$ of $G_{\mu_{\widetilde\tau_1}\widetilde T_0,\widetilde \rho}$. Then the diagonals of $G(p_i),i\in I$ are labeled $\widetilde \tau_1$.

For each $\widetilde R\in \psi'_{\widetilde \rho}(\widetilde P'')$, $\widetilde R\in \psi'_{\widetilde \tau_1,\widetilde\rho}(\widetilde R')$ for some $\widetilde R'$ which satisfies $\widetilde R'\in \psi'_{\widetilde \tau_2,\widetilde\rho}(\widetilde P'')$. By Lemma \ref{diff2} in the unpunctured surface case, $\widetilde R'$ can be obtained from $\widetilde P'$ by twists on some tiles with diagonals labeled $\widetilde \tau_2$, we may assume $\widetilde R'=\prod_{j\in J}\mu_{p_j}\widetilde P'$ for some $J$ and $p_j$. Then $\widetilde P'=\prod_{j\in J}\mu_{p_j}\widetilde R'$. We assume $p_j$ corresponds to the tile $G(p_j)$ of $G_{\mu_{\widetilde\tau_1}\widetilde T_0,\widetilde \rho}$. Then the diagonals of $G(p_j),j\in J$ are labeled $\widetilde \tau_2$. Let $\widetilde O'=\prod_{i\in I}\mu_{p_i}\widetilde R'$. As the diagonals of $G(p_i),i\in I$ are labeled $\widetilde \tau_1$, by Lemma \ref{diff2} in the unpunctured surface case, $\widetilde R\in \psi'_{\widetilde \tau_1,\widetilde\rho}(\widetilde R')=\psi'_{\widetilde \tau_1,\widetilde\rho}(\widetilde O')$. Since $\widetilde \tau_1$ and $\widetilde \tau_2$ are not in the same triangle in $\widetilde T_0$, $\prod_{i\in I}\mu_{p_i}$ commutes with $\prod_{j\in J}\mu_{p_j}$. Thus, $\widetilde O'=\prod_{j\in J}\mu_{p_j}\widetilde Q'$. As the diagonals of $G(p_j),j\in J$ are labeled $\widetilde \tau_2$ and $\widetilde Q'\in \psi'_{\widetilde \tau_2,\widetilde\rho}(\widetilde Q'')$, by Lemma \ref{diff2} in the unpunctured surface case, $\widetilde O'\in \psi'_{\widetilde \tau_2,\widetilde\rho}(\widetilde Q'')$. Thus, $\widetilde R\in \psi'_{\widetilde \rho}(\widetilde Q'')$ and hence $\psi'_{\widetilde \rho}(\widetilde P'')\subset \psi'_{\widetilde \rho}(\widetilde Q'')$. It can be proved similarly the inverse inclusion holds. Therefore, $\psi'_{\widetilde \rho}(\widetilde P'')=\psi'_{\widetilde \rho}(\widetilde Q'')$. It can be proved similarly that $\psi_{\widetilde \rho}(\widetilde P)=\psi_{\widetilde \rho}(\widetilde Q)$ if $\psi_{\widetilde \rho}(\widetilde P)\cap \psi_{\widetilde \rho}(\widetilde Q)\neq \emptyset$. Moreover, since $\psi_{\widetilde \tau_1,\widetilde\rho}$ and $\psi'_{\widetilde \tau_1,\widetilde\rho}$ are partition bijections inverse to each other and $\psi_{\widetilde \tau_2,\widetilde\rho}$ and $\psi'_{\widetilde \tau_2,\widetilde\rho}$ are partition bijections inverse to each other, $\widetilde P''\in \psi_{\widetilde \rho}(\widetilde P)$ if and only if $\widetilde P\in \psi'_{\widetilde \rho}(\widetilde P'')$. Thus, the claim follows by Lemma \ref{partition3}.

Therefore, the partition maps $\psi_{\widetilde\rho}$ and $\psi'_{\widetilde\rho}$ induce the required partition maps $\psi_{\rho}$ and $\psi_{\rho}$. Since $n^{\widetilde\tau_1}(\widetilde T_0,\widetilde \rho)=n^{\mu_{\widetilde \tau_1}\widetilde \tau_1}(\mu_{\widetilde\tau_1}\widetilde T_0,\widetilde\rho)$ and $n^{\widetilde\tau_2}(\mu_{\widetilde\tau_1}\widetilde T_0,\widetilde \rho)=n^{\mu_{\widetilde \tau_2}\widetilde \tau_2}(\mu_{\widetilde\tau_2}\mu_{\widetilde\tau_1}\widetilde T_0,\widetilde\rho)$, we have $n^{\tau}(T,\rho)=n^{\tau'}(T',\rho)$.

The above discussion implies that the result holds for all the surfaces with $\mathcal O_0$ admits at most one orbifold point.

(II) $\widetilde{\mathcal O}_0$ contains orbifold points.

By the construction in Subsection \ref{local1}, $\widetilde{\mathcal O}_0$ contains at most one orbifold point. In this case, it can be proved similarly by the unfolding of $\widetilde{\mathcal O}_0$.
\end{proof}

\medskip

\begin{Lemma}\label{formglobal}

Let $P=(P_i)\in \mathcal P(G_{T,\zeta})$. For any $(P'_i)$ such that $P'_i\in \psi_{\zeta_i}(P_i), i\in [1,k]$, $(P'_i)$ forms a perfect matching of $G_{T',\zeta}$. Moreover, if $u_j\notin P_j\cap P_{j+1}$ for some $j\in [1,k-1]$, then $u'_j\notin P'_j\cap P'_{j+1}$, where $u_j$ are the gluing edges.

\end{Lemma}

\begin{proof}

For any $j\in [1,k-1]$, since $(P_i)\in \mathcal P(G_{T,\zeta})$, $u_j\in P_j\cup P_{j+1}$ by Proposition \ref{decompose}. Hence $u'_j\in P'_j\cup P'_{j+1}$ by Proposition \ref{compare-loc} (2.c). Therefore, by the dual version of Proposition \ref{decompose}, $(P'_i)$ forms a perfect matching of $G_{T',\zeta}$. By Proposition \ref{compare-loc} (2.c), if $u_j\notin P_j\cap P_{j+1}$ for some $j\in [1,k-1]$, then $u'_j\notin P'_j\cap P'_{j+1}$.
\end{proof}

\medskip

The following observation would help us to prove Theorem \ref{partition bi}.

\medskip

\begin{Lemma}\label{diff2}

With the same notation as above,

\begin{enumerate}[$(1)$]

\item given $P\in \mathcal P^{\tau}_{\nu}(G_{T,\rho})$ for some $\nu$, for any $P',Q'\in \mathcal \psi_{\rho}(P)$, $Q'$ can be obtained from $P'$ by twists on some tiles with diagonals labeled $\tau'$. Moreover, $|\psi_{\rho}(P)|=2^{\sum_{\nu_l=1}\nu_l}$. More precisely, for any $\nu_i=1$ and $P'\in \psi_{\rho}(P)$, $P'$ can twist on the tile corresponding the $i$-th $\tau'$-equivalence class.

\item given $P'\in \mathcal P^{\tau}_{\nu}(G_{T',\rho})$ for some $\nu$, for any $P,Q\in \mathcal \psi'_{\rho}(P')$, $Q$ can be obtained from $P$ by twists on some tiles with diagonals labeled $\tau$. Moreover, $|\psi'_{\rho}(P')|=2^{\sum_{\nu_l=1}\nu_l}$. More precisely, for any $\nu_i=1$ and $P\in \psi'_{\rho}(P')$, $P$ can twist on the tile corresponding the $i$-th $\tau$-equivalence class.

\end{enumerate}

\end{Lemma}

\begin{proof}

We shall only prove (1) by dually. From the proof of Proposition \ref{compare-loc}, it suffices to prove the case that $\mathcal O_0$ contains no orbifold point. When $\mathcal O_0$ contains no orbifold point, the result follows immediately by the constructions of partition bijections in Subsection \ref{local2}.
\end{proof}

\medskip

The following observation are useful to prove Theorem \ref{mainthm}.

\medskip

\begin{Lemma}\label{maxtomax}
With the same notation as above,
\begin{enumerate}[$(1)$]

  \item $P_{+}(G_{T',\rho})\subset \psi_{\rho}(P_{+}(G_{T,\rho}))$,\;\;  $P_{-}(G_{T',\rho})\subset \psi_{\rho}(P_{-}(G_{T,\rho}))$.

  \item $P_{+}(G_{T,\rho})\subset \psi'_{\rho}(P_{+}(G_{T',\rho}))$,\;\;  $P_{-}(G_{T,\rho})\subset \psi'_{\rho}(P_{-}(G_{T',\rho}))$.

\end{enumerate}

\end{Lemma}

\begin{proof}

We shall only prove (1) because (2) can be proved dually. Since $G_{T,\rho}\cong G_{\widetilde T_0,\widetilde \rho}, G_{T',\rho}\cong G_{\widetilde T'_0,\widetilde \rho}$ as graphs and the orientation of $\widetilde{\mathcal O}_0$ is compatible with the orientation of $\mathcal O_0$ under the map $\chi$, it suffices to prove the result holds for the case $\mathcal O_0$ contains no orbifold points. In such case, the result follows immediately by the construction of the partition bijections in Subsection \ref{local2}.
\end{proof}

\medskip

\begin{Lemma}\label{com}

With the same notation as above,

\begin{enumerate}[$(1)$]

  \item we assume that $P\in \mathcal P(G_{T,\rho})$ can twist on a tile $G(p)$ with diagonal labeled $a=a_q,q=1,2,3,4$. If all $\tau$-mutable edges pairs in $P$ are labeled $a_{q-1},a_{q+1}$, then there exists $P'\in \psi_{\rho}(P)$ so that the following hold

       \begin{enumerate}[$(a)$]

         \item $P'$ can twist on $G'(p)$, where $G'(p)$ is the tile of $G_{T',\rho}$ determined by $p$.

         \item $\mu_pP'\in \psi_{\rho}(\mu_{p}P)$.

         \item All $\tau'$-mutable edges pairs in $P'$ are labeled $a_{q-1}, a_{q+1}$.

         \item $n^+_p(a,P)-m^+_p(a,\rho)=n^+_p(a,P')-m^+_p(a,\rho)$ and $n^-_p(a,P)-m^-_p(a,\rho)=n^-_p(a,P')-m^-_p(a,\rho)$.

       \end{enumerate}

  \item we assume that $P'\in \mathcal P(G_{T',\rho})$ can twist on a tile $G'(p)$ with diagonal labeled $a=a_q,q=1,2,3,4$. If all $\tau'$-mutable edges pairs in $P'$ are labeled $a_{q-1},a_{q+1}$, then there exists $P\in \psi'_{\rho}(P')$ so that the following hold

      \begin{enumerate}[$(a)$]

         \item $P$ can twist on $G(p)$, where $G(p)$ is the tile of $G_{T,\rho}$ determined by $p$.

         \item $\mu_pP\in \psi'_{\rho}(\mu_{p}P')$.

         \item All $\tau$-mutable edges pairs in $P$ are labeled $a_{q-1}, a_{q+1}$.

         \item $n^+_p(a,P)-m^+_p(a,\rho)=n^+_p(a,P')-m^+_p(a,\rho)$ and $n^-_p(a,P)-m^-_p(a,\rho)=n^-_p(a,P')-m^-_p(a,\rho)$.

       \end{enumerate}

\end{enumerate}

\end{Lemma}

\begin{proof}

We shall only prove (1) because (2) can be proved dually.

We first consider the case that $\mathcal O$ is an unpunctured surface.

When $a_1\neq a_3, a_2\neq a_4$, by checking the partition bijection between $\mathcal P(G_{T,\rho})$ and $\mathcal P(G_{T',\rho})$, it can be seen that the result holds.

When $a_1= a_3$ or $a_2= a_4$, we may assume $a_1=a_3$. The cases (1,2,9) can be verified similarly as the case $a_1\neq a_3,a_2\neq a_4$.

In case (6), $a=a_1$ and $P\in \mathcal H'_{(\lambda_1,\cdots,\lambda_{s+1})}(a_1,\tau)$ for some $\lambda\in \{0,1\}^{s+1}$. Thus, $\psi_{\rho}(P)=\mathcal H_{(1-\lambda_1,\cdots,1-\lambda_{s+1})}(\tau',a_1)$. We assume that $G(p)$ is the $(2i-1)$-th tile. Then $\mu_{p}P\in \mathcal H'_{(\lambda_1,\cdots,\lambda_{i-1},1-\lambda_i,\lambda_{i+1},\cdots,\lambda_{s+1})}(a_1,\tau)$. For any $P'\in\psi_{\rho}(P)$, $P'$ can twist on $G'(p)$, the $2i$-th tile of $G_{T',\rho}$, and $\mu_pP'\in \mathcal H_{(1-\lambda_1,\cdots,\lambda_i,\cdots,1-\lambda_{s+1})}(\tau',a_1)=\psi_{\rho}(\mu_pP)$. Let $P'\in \psi_{\rho}(P)$ so that all the $\tau'$-mutation edges pairs are labeled $a_2,a_4$. Then $n^+_p(a,P)-m^+_p(a,\rho)=n^+_p(a,P')-m^+_p(a,\rho)=s+1-i$ and $n^-_p(a,P)-m^-_p(a,\rho)=n^-_p(a,P')-m^-_p(a,\rho)=i-1$.

In case (7), if $a=a_1$, then the result can be proved similarly as case (6). If $a\neq a_1$, then $a=a_4$ and $G(p)$ is the last tile of $G_{T,\rho}$.

If the upper right edge labeled $a_8$ of $G(p)$ is in $P$, then $P\in \mathcal H'_{(\lambda_1,\cdots,\lambda_{s})}(a_1,\tau)$ for some $\lambda\in \{0,1\}^{s}$. Thus $\psi_{\rho}(P)=\mathcal H_{(1-\lambda_1,\cdots,1-\lambda_{s})}(\tau',a_1)$. We have the upper right edge labeled $a_8$ of $G'(p)$ is in $P'$ for any $P'\in \psi_{\rho}(P)$, where $G'(p)$ is the last tile of $G_{T',\rho}$. Thus, $P'$ can twist on the tile $G'(p)$ for any $P'\in \psi_{\rho}(P)$. Let $P'\in \psi_{\rho}(P)$ such that all the $\tau'$-mutation edges pairs are labeled $a_1$. Since $\mu_pP\in \mathcal H'_{(\lambda_1,\cdots,\lambda_{s},1)}(a_1,\tau)$ and $\mu_pP'\in \mathcal H_{(1-\lambda_1,\cdots,1-\lambda_{s},1)}(\tau',a_1)$, $\mu_pP'\in \psi_{\rho}(\mu_pP)$.

If the upper right edge labeled $a_7$ of $G(p)$ is in $P$, then the upper right edge labeled $a_8$ of $G(p)$ is in $\mu_pP$. Change the roles of $P$ and $\mu_pP$, this case can be proved as the above case.

In both cases, by Lemma \ref{non-tau-mu1} and Lemma \ref{non-tau-mu2}, the edge labeled $a_4$ of the $(2i-1)$-th tile of $G_{T,\rho}$ is non-$\tau$-mutable in $P$ if and only if the edge labeled $a_4$ of the $(2i+1)$-th tile of $G_{T',\rho}$ is non-$\tau'$-mutable in $P'$, we have $n^+_p(a,P)-m^+_p(a,\rho)=n^+_p(a,P')-m^+_p(a,\rho)$ and $n^-_p(a,P)-m^-_p(a,\rho)=n^-_p(a,P')-m^-_p(a,\rho)$, note that the edge of the first tile of $G_{T',\rho}$ labeled $a_4$ can not be a non-$\tau'$-mutable edge.

Therefore, the result holds in case (7).

For case (8), the result can be proved similarly to cases (6) and (7). Therefore, the result follows in cases (4,5,6). Dually, the result can be proved in cases (10,11,12).

The cases (3,4,5) can be proved similarly.

Dually, the result can be verified if we change the roles of $\tau$ and $\tau'$. Therefore, the result follows in case $\mathcal O$ contains no orbifold point.

We now consider the case that $\mathcal O$ contains orbifold points. From the proof of Proposition \ref{compare-loc}, $\psi_{\rho}$ and $\psi'_{\rho}$ are induced by two partition bijections $\widetilde\psi_{\widetilde\rho}$ and $\widetilde\psi'_{\widetilde\rho}$, where $\widetilde\rho$ is an arc in an unpunctured surface. Therefore, the result follows by the corresponding result for $\widetilde\psi_{\widetilde\rho}$ and $\widetilde\psi'_{\widetilde\rho}$.

The proof is complete.
\end{proof}

\medskip

\section{Proof of the Theorem \ref{partition bi}}\label{main1}

In this section, the proof of Theorem \ref{partition bi} will be given. We follow the idea in \cite{H,H1}. More precisely, let $\gamma$ be an oriented arc, $\zeta$ be the corresponding arc and $\zeta_i,i\in [1,k]$ be the subcurves of $\zeta$ constructed at the beginning of Section \ref{compare}. Then the partition bijection between $\mathcal P(G_{T,\gamma})$ and $\mathcal P(G_{T',\gamma})$ is constructed by the partition bijections between $\mathcal P(G_{T,\zeta_i})$ and $\mathcal P(G_{T',\zeta_i})$.

\medskip

\subsection{Proof of the Theorem \ref{partition bi} (1)}\label{main11}

By Proposition \ref{decompose}, we would write a perfect matching $P\in \mathcal P(G_{T,\gamma})$ as $(P_i)$ in the sequel. We arrange the order of $\tau$-equivalence classes in $G_{T,\zeta_i}$ according to the orientation of $\zeta$. By Lemma \ref{localdec}, for any perfect matching $P$ of $G_{T,\gamma}$, $P\in \mathcal P^{\tau}_{\nu}(G_{T,\gamma})$ for some sequence $\nu$. By the truncation of $\zeta$, $\nu$ can be truncated as $\nu^1,\cdots,\nu^k$ so that $P_i\in \mathcal P^{\tau}_{\nu^i}(G_{T,\zeta_i})$ for any $i\in [1,k]$. Similarly, we would write $P'\in \mathcal P(G_{T',\gamma})$ as $(P'_i)$ with $P'_i\in \mathcal P(G_{T',\zeta_i})$.

\medskip

For a sequence $\nu=(\nu_l)\in S_{T,\zeta}$, as in \cite{H,H1}, we define $\nu$-pairs as follows. Firstly, choose a pair $(\nu_s,\nu_t), s<t$ satisfies (1) $\nu_s\nu_t=-1$; (2) $\nu_l=0$ for all $s<l<t$; (3) $\nu_l\nu_s\geq 0$ for all $l<s$, and we call $\{s,t\}$ a \emph{$\nu$-pair}. Then delete $\nu_s,\nu_t$ from $(\nu_l)$, we get a sequence $(\nu'_l)$ and we do the same step on $(\nu'_l)$ as on $(\nu_l)$. Particularly, the set of $\nu$-pairs is empty if $\{s,t\}$ does not exist in the first step. It is easy to see the set of $\nu$-pairs equals to the set of $(-\nu)$-pairs.

\medskip

Let $P\in \mathcal P^{\tau}_{\nu}(G_{T,\gamma})$ for some $\nu$. If $\nu_s=-1$, then the $s$-th $\tau$-equivalence class in $G_{T,\gamma}$ corresponds to a diagonal of a tile, denote by $G(p_{l_s})$, of $G_{T,\gamma}$. We similarly denote by $G'(p_{l_t})$ the tile of $G_{T',\gamma}$ which corresponds to the $t$-th $\tau'$-equivalence class. Let $\pi(P)$ be the subset of $\mathcal P(G_{T',\gamma})$ containing all $P'=(P'_i)$ such that $P'_i\in \psi_{\zeta_i}(P_i)$ and the edges in $P'\cap E(G'(p_{l_t}))$ have the same labels with the edges in $P\cap E(G(p_{l_s}))$ for any $\nu$-pair $\{s,t\}$ with $\nu_t=1$.

\medskip

Dually, change the roles of $\tau$ and $\tau'$, for any $P'\in \mathcal P^{\tau'}_{\nu}(G_{T',\gamma})$, let $\pi'(P')$ be the subset of $\mathcal P(G_{T,\gamma})$ containing all $P=(P_i)$ such that $P_i\in \psi'_{\zeta_i}(P'_i)$ and the edges in $P\cap E(G(p_{l_s}))$ have the same labels with the edges in $P'\cap E(G'(p_{l_t}))$ for any $\nu$-pair $\{s,t\}$ with $\nu_s=1$.

\medskip

We have the following characterization of $\pi(P)$ for $P\in \mathcal P^{\tau}_{\nu}(G_{T,\gamma})$.

\medskip

\begin{Proposition}\label{2^n}

With the same notation as above, we assume $\tau'\neq \gamma \notin T$. For any $P\in \mathcal P^{\tau}_{\nu}(G_{T,\gamma})$, we have $|\pi(P)|=2^{max\{0, \sum\nu_l\}}$. More precisely, for each $s\in [1,\sum\nu_l]$ not in a $\nu$-pair with $\nu_s=1$ and $P'\in \mathcal P(G_{T',\gamma})$, $P'$ can twist on the tile $G'_{l_{s}}$. Moreover, any two perfect matchings in $\pi(P)$ can obtained from each other by such kind of twists.

\end{Proposition}

\begin{proof}

It follows by the construction of $\pi(P)$ and Lemma \ref{diff2}.
\end{proof}

\medskip

Given $\nu$ with $\sum\nu_l\geq 0$, let $d=\sum\nu_l$. Let $s_1,\cdots,s_d$ be the indices in order such that $\nu_{s_i}=1$ and $s_i$ is not in a $\nu$-pair for $i\in [1,d]$. Given $P\in \mathcal P_{\nu}^{\tau}(G_{T,\gamma})$, by Proposition \ref{2^n}, for each $\lambda\in \{0,1\}^{\sum\nu_l}$, there uniquely exists $P'(\lambda)\in \pi(P)$ such that the edges of $G'(p_{l_{s_i}})$ labeled $a_1,a_3$ are in $P'(\lambda)$ when $\lambda_i=1$ and the edges of $G'(p_{l_{s_i}})$ labeled $a_2,a_4$ are in $P'(\lambda)$ when $\lambda_i=0$ for $i\in [1,\sum\nu_l]$. Dually, given $P'\in \mathcal P^{\tau'}_{\nu}(G_{T',\gamma})$, for each $\lambda\in \{0,1\}^d$, there uniquely exists $P(\lambda)\in \pi'(P')$ such that the edges of $G(p_{l_{s_i}})$ labeled $a_1,a_3$ are in $P(\lambda)$ when $\lambda_i=1$ and the edges of $G(p_{l_{s_i}})$ labeled $a_2,a_4$ are in $P(\lambda)$ when $\lambda_i=0$ for $i\in [1,\sum\nu_l]$.

\medskip

We are now ready to prove Theorem \ref{partition bi} (1).

\medskip

{\bf Proof of Theorem \ref{partition bi} (1):} For any $P=(P_i), Q=(Q_i)\in \mathcal P(G_{T,\gamma})$ with $\pi(P)\cap \pi(Q)\neq \emptyset$, we have $\psi_{\zeta_i}(P_i)\cap \psi_{\zeta_i}(Q_i)\neq \emptyset$ for any $i$ and $P\cap E(G(p_{l_s}))=Q\cap E(G(p_{l_s}))$ for all $\nu$-pair $\{s,t\}$ with $\nu_s=-1$. Since $\psi_{\zeta_i}$ are partition bijections, $\psi_{\zeta_i}(P_i)=\psi_{\zeta_i}(Q_i)$. Therefore, by the construction of $\pi$, we have $\pi(P)=\pi(Q)$.

Dually, we have $\pi'(P')=\pi'(Q')$ if $\pi'(P')\cap \pi'(Q')\neq \emptyset$ for any $P',Q'\in \mathcal P(G_{T',\gamma})$.

Given $P\in \mathcal P(G_{T,\gamma})$ and $P'\in \pi(P)$, since the set of $\nu$-pairs is the same as the set of $(-\nu)$-pairs and $\psi_{\zeta_i}$ and $\psi'_{\zeta_i}$ are partition bijections inverse to each other, we have $P\in \pi'(P')$. Dually, we have $P'\in \pi(P)$ for any $P'\in \mathcal P(G_{T',\gamma})$ and $P\in \mathcal P(G_{T,\gamma})$ with $P\in \pi'(P')$.

By Lemma \ref{partition3} and Proposition \ref{2^n}, the result follows.\ \ \ \ \ \ \ \;\;\; \ \ \ \ \ \ \ \ \ \ \ \ \ \ \ \ \ \ \ \ \ \ \ \  $\square$

\medskip

\subsection{Proof of the Theorem \ref{partition bi} (2)}\label{main12}

Herein, we prove Theorem \ref{partition bi} (2). We follow the strategy in \cite{H1}.

\begin{Lemma}\label{basic1}

\begin{enumerate}[$(1)$]

  \item Let $Q\in \mathcal P_{\nu}^{\tau}(G_{T,\gamma})$ for $\sum\nu_l\geq 0$ and $Q'\in \pi(Q)$ corresponding to some $\lambda\in \{0,1\}^{\sum\nu_l}$ as in Proposition \ref{2^n}.

  \begin{enumerate}[$(a)$]

    \item If $\lambda=(0,0,\cdots,0)$, then $\bigoplus_{P'\in \pi(Q)}y^{T'}(P')=y^{T'}(Q')\cdot(1\oplus y^{T'}_{\tau'})^{\sum\nu_l}$ in $\mathbb P$.

    \item If $\lambda=(1,1,\cdots,1)$, then $\bigoplus_{P'\in \pi(Q)}y^{T'}(P')=y^{T'}(Q')\cdot(1\oplus y^{T}_{\tau})^{\sum\nu_l}$ in $\mathbb P$.

  \end{enumerate}

  \item Let $Q'\in \mathcal P_{\nu}^{\tau'}(G_{T',\gamma})$ for $\sum\nu_l\geq 0$ and $Q\in \pi'(Q')$ corresponding to some $\lambda\in \{0,1\}^{\sum\nu_l}$ as in Proposition \ref{2^n}.

  \begin{enumerate}[$(a)$]

    \item If $\lambda=(0,0,\cdots,0)$, then $\bigoplus_{P\in \pi'(Q')}y^{T}(P)=y^{T}(Q)\cdot(1\oplus y^{T'}_{\tau'})^{\sum\nu_l}$ in $\mathbb P$.

    \item If $\lambda=(1,1,\cdots,1)$, then $\bigoplus_{P\in \pi'(Q')}y^{T}(P)=y^{T}(Q)\cdot(1\oplus y^{T}_{\tau})^{\sum\nu_l}$ in $\mathbb P$.

  \end{enumerate}

\end{enumerate}

\end{Lemma}

\begin{proof}

The proof is the same as that of \cite[Lemma 7.7]{H1}.
\end{proof}

\medskip

\begin{Lemma}\label{non-a-F}

With the same notation as above, we assume $\tau'\neq \gamma\notin T$. Suppose that $P\in \mathcal P^{\tau}_{\nu}(G_{T,\gamma})$ can twist on a tile $G(p)$ with diagonal labeled $\varsigma$. We assume $S_1,S_2\in \mathfrak P$ so that $Q\in S_1$ and $\mu_{p}Q\in S_2$. If $\varsigma\neq a_1,a_2,a_3,a_4,\tau$, then in $\mathbb P$,
$$\frac{\bigoplus_{P\in S_1}y^{T}(P)}{\bigoplus_{P'\in \pi(S_1)}y^{T'}(P')}=\frac{\bigoplus_{P\in S_2}y^{T}(P)}{\bigoplus_{P'\in \pi(S_2)}y^{T'}(P')}.$$

\end{Lemma}

\begin{proof}

The proof is the same as that of \cite[Lemma 7.8]{H1}.
\end{proof}

\medskip

\begin{Lemma}\label{tau-F}

With the same notation as above, we assume $\tau'\neq \gamma\notin T$. Suppose that $P\in \mathcal P^{\tau}_{\nu}(G_{T,\gamma})$ can twist on a tile $G(p_{l_t})$ with diagonal labeled $\tau$. We assume $S_1,S_2\in \mathfrak P$ so that $Q\in S_1$ and $\mu_{p_{l_t}}Q\in S_2$. Then in $\mathbb P$,
$$\frac{\bigoplus_{P\in S_1}y^{T}(P)}{\bigoplus_{P'\in \pi(S_1)}y^{T'}(P')}=\frac{\bigoplus_{P\in S_2}y^{T}(P)}{\bigoplus_{P'\in \pi(S_2)}y^{T'}(P')}.$$

\end{Lemma}

\begin{proof}

The proof is the same as the proof of \cite[Lemma 7.9]{H1}.
\end{proof}

\medskip

\begin{Lemma}\label{special-a-F}

With the same notation as above, we assume $\tau'\neq \gamma\notin T$. Let $Q$ be a perfect matching of $G_{T,\gamma}$ that can twist on a tile $G(p)$ with diagonal labeled $a=a_q$ for $q=1,2,3,4$. We assume $S_1,S_2\in \mathfrak P$ so that $Q\in S_1$ and $\mu_{p}Q\in S_2$. If all the $\tau$-mutable edges pairs in $Q$ are labeled $a_{q-1},a_{q+1}$, then in $\mathbb P$,
$$\frac{\bigoplus_{P\in S_1}y^{T}(P)}{\bigoplus_{P'\in \pi(S_1)}y^{T'}(P')}=\frac{\bigoplus_{P\in S_2}y^{T}(P)}{\bigoplus_{P'\in \pi(S_2)}y^{T'}(P')}.$$

\end{Lemma}

\begin{proof}

We may assume that the edge labeled $\tau$ of $G(p)$ is in $Q$. We write $Q$ as $(Q_i)$ and assume that $G(p)$ is a tile of $G_{T,\zeta_j}$ for some $j$. Thus $\mu_pQ=(\mu_pQ_i)$ with $\mu_pQ_i=Q_i$ for $i\neq j$. By Lemma \ref{com}, there exists $Q'_j\in \psi_{\zeta_j}(Q_j)$ satisfies (a) $Q'_j$ can twist on $G'(p)$, (b) $\mu_pQ'_j\in \psi_{\zeta_j}(\mu_pQ_j)$, (c) all $\tau'$-mutable edges pairs in $Q'_j$ are labeled $a_{q-1},a_{q+1}$. For each $i\neq j$, since all the $\tau$-mutable edges pairs in $Q$ are labeled $a_{q-1},a_{q+1}$, by Lemma \ref{diff2}, we can choose $Q_i'\in \psi_{\zeta_i}(Q_i)$ such that all $\tau'$-mutation edges pairs in $Q'_i$ are labeled $a_{q-1},a_{q+1}$. Since all $\tau$-mutable edges pairs in $Q$ are labeled $a_{q-1},a_{q+1}$, all $\tau$-mutable edges pairs in $\mu_pQ$ are labeled $a_{q-1},a_{q+1}$. By the constructions of $\pi(Q)$ and $\pi(\mu_pQ)$, we have $Q'=(Q'_i)\in \pi(Q)$ and $Q''=(\mu_pQ'_i)\in \pi(\mu_pQ)$, here $\mu_pQ'_i=Q'_i$ for $i\neq j$. Clearly, $\mu_{p}Q'=Q''$. We assume $\mu_pQ\in \mathcal P^{\tau}_{\nu'}(G_{T,\gamma})$. It is clear $\sum\nu_l>\sum\nu'_l$ and $|b^{T}_{\tau a_q}|=\sum\nu_l-\sum\nu'_l$.

(I) In case $\sum\nu_l\geq 0$ and $\sum\nu'_l\geq 0$, we have $S_1=\{Q\}$ and $S_2=\{\mu_{p}Q\}$.

If $a=a_1,a_3$, since all mutation edges pairs in $Q'$ and $Q''$ are labeled $a_2,a_4$, $Q'$ and $Q''$ correspond to $(0,0,\cdots,0)\in \{0,1\}^{\sum\nu_l}$ and $(0,0,\cdots,0)\in \{0,1\}^{\sum\nu'_l}$, respectively. By Lemma \ref{basic1},
$$\textstyle\bigoplus_{P'\in \pi(S_1)}y^{T'}(P')=y^{T'}(Q')\cdot(1\oplus y_{\tau'}^{T'})^{\sum\nu_l},$$
and
$$\textstyle\bigoplus_{P'\in \pi(S_2)}y^{T'}(P')=y^{T'}(Q'')\cdot(1\oplus y_{\tau'}^{T'})^{\sum\nu'_l}.$$
By Lemma \ref{brick1}, $\frac{y^{T'}(Q')}{y^{T'}(Q'')}=y^{T'}_{a_q}$ and $\frac{y^{T}(Q)}{y^{T}(\mu_p Q)}=y^{T}_{a_q}$. Since $y^{T'}_{\tau'}=(y^{T}_{\tau})^{-1}$ and
$$y^{T'}_{a_q}=y^{T}_{a_q}\cdot (y^T_{\tau})^{\sum\nu_l-\sum\nu'_l}\cdot(1\oplus y^{T}_{\tau})^{\sum\nu'_l-\sum\nu_l},$$
$y_{a_q}^{T}=y_{a_q}^{T'}\cdot(1\oplus y^{T'}_{\tau'})^{\sum\nu_l-\sum\nu'_l}$.
Thus, the result follows in this case.

If $a=a_2,a_4$, since all mutation edges pairs in $Q'$ and $Q''$ are labeled $a_1,a_3$, $Q'$ and $Q''$ correspond to $(1,1,\cdots,1)\in \{0,1\}^{\sum\nu_l}$ and $(1,1,\cdots,1)\in \{0,1\}^{\sum\nu'_l}$, respectively. By Lemma \ref{basic1},
$$\textstyle\bigoplus_{P'\in \pi(S_1)}y^{T'}(P')=y^{T'}(Q')\cdot(1\oplus y_{\tau}^{T})^{\sum\nu_l},$$
and
$$\textstyle\bigoplus_{P'\in \pi(S_2)}y^{T'}(P')=y^{T'}(Q'')\cdot(1\oplus y_{\tau}^{T})^{\sum\nu'_l}.$$
By Lemma \ref{brick1}, $\frac{y^{T'}(Q'')}{y^{T'}(Q')}=y^{T'}_{a_q}$ and $\frac{y^{T}(\mu_p Q)}{y^{T}(Q)}=y^{T}_{a_q}$. Because of $y^{T'}_{a_q}=y^{T}_{a_q}\cdot(1\oplus y^{T}_{\tau})^{\sum\nu_l-\sum\nu'_l},$ the result follows in this case.

(II) In case $\sum\nu_l\geq 0$ and $\sum\nu'_l< 0$, we have $S_1=\{Q\}$ and $\pi(S_2)=\{Q''\}$.

If $q=1,3$, since all mutation edges pairs in $Q'$ and $\mu_pQ$ are labeled $a_2,a_4$, $Q'$ and $\mu_pQ$ correspond to $(0,0,\cdots,0)\in \{0,1\}^{\sum\nu_l}$ and $(0,0,\cdots,0)\in \{0,1\}^{-\sum\nu'_l}$, respectively. By Lemma \ref{basic1},
$$\textstyle\bigoplus_{P'\in \pi(S_1)}y^{T'}(P')=y^{T'}(Q')\cdot(1\oplus y_{\tau'}^{T'})^{\sum\nu_l},$$
and
$$\textstyle\bigoplus_{P\in S_2}y^{T}(P)=y^{T}(\mu_pQ)\cdot(1\oplus y_{\tau'}^{T'})^{-\sum\nu'_l}.$$
By Lemma \ref{brick1}, $\frac{y^{T'}(Q')}{y^{T'}(Q'')}=y^{T'}_{a_q}$ and $\frac{y^{T}(Q)}{y^{T}(\mu_p Q)}=y^{T}_{a_q}$. Since $y^{T'}_{\tau'}=(y^{T}_{\tau})^{-1}$ and
$$y^{T'}_{a_q}=y^{T}_{a_q} \cdot(y^T_{\tau})^{\sum\nu_l-\sum\nu'_l}\cdot(1\oplus y^{T}_{\tau})^{\sum\nu'_l-\sum\nu_l},$$
$y_{a_q}^{T}=y_{a_q}^{T'}\cdot(1\oplus y^{T'}_{\tau'})^{\sum\nu_l-\sum\nu'_l}$, and hence
$$\frac{y^{T}(\mu_pQ)\cdot(1\oplus y_{\tau'}^{T'})^{-\sum\nu'_l}}{y^{T}(Q)}=\frac{y^{T'}(Q'')}{y^{T'}(Q')\cdot(1\oplus y^{T'}_{\tau'})^{\sum\nu_l}}.$$ Thus, the result follows in this case.

If $q=2,4$, since all mutation edges pairs in $Q'$ and $\mu_pQ$ are labeled $a_1,a_3$, $Q'$ and $\mu_pQ$ correspond to $(1,1,\cdots,1)\in \{0,1\}^{\sum\nu_l}$ and $(1,1,\cdots,1)\in \{0,1\}^{-\sum\nu'_l}$, respectively. By Lemma \ref{basic1},
$$\textstyle\bigoplus_{P'\in \pi(S_1)}y^{T'}(P')=y^{T'}(Q')\cdot(1\oplus y_{\tau}^{T})^{\sum\nu_l},$$
and
$$\textstyle\bigoplus_{P\in S_2}y^{T}(P)=y^{T}(\mu_pQ)\cdot(1\oplus y_{\tau}^{T})^{-\sum\nu'_l}.$$
By Lemma \ref{brick1}, $\frac{y^{T'}(Q'')}{y^{T'}(Q')}=y^{T'}_{a_q}$ and $\frac{y^{T}(\mu_p Q)}{y^{T}(Q)}=y^{T}_{a_q}$. Because of $y^{T'}_{a_q}=y^{T}_{a_q}\cdot(1\oplus y^{T}_{\tau})^{\sum\nu_l-\sum\nu'_l},$ the result follows in this case.

(III) In case $\sum\nu_l< 0$ and $\sum\nu'_l< 0$, change the roles of $T$ and $T'$, the result follows by the similar discussion in the case $\sum\nu_l,\sum\nu'_l\geq 0$.
\end{proof}

\medskip

\begin{Lemma}\label{all-a-F}

With the same notation as above, we assume that $\tau'\neq \gamma\notin T$. Let $Q$ be a perfect matching of $G_{T,\gamma}$ that can twist on a tile $G(p)$ with diagonal labeled $a_q$ for $q=1,2,3,4$. We assume that $S_1,S_2\in \mathfrak P$ so that $Q\in S_1$ and $\mu_{p}Q\in S_2$. Then in $\mathbb P$, we have
$$\frac{\bigoplus_{P\in S_1}y^{T}(P)}{\bigoplus_{P'\in \pi(S_1)}y^{T'}(P')}=\frac{\bigoplus_{P\in S_2}y^{T}(P)}{\bigoplus_{P'\in \pi(S_2)}y^{T'}(P')}.$$

\end{Lemma}

\begin{proof}

After $Q$ twists on the tiles $G$ which satisfy the diagonals labeled $\tau$ and the edges labeled $a_{q},a_{q+2}$ are in $Q$, we obtain a perfect matching $R$ which satisfies the conditions of Lemma \ref{special-a-F}. Then the result follows by Lemmas \ref{special-a-F} and \ref{tau-F}.
\end{proof}

\medskip

In summary, we have the following proposition.

\medskip

\begin{Proposition}\label{all-F}

With the same notation as above, we assume that $\tau'\neq \gamma\notin T$. Let $Q$ be a perfect matching of $G_{T,\gamma}$ that can twist on a tile $G(p)$. We assume that $S_1,S_2\in \mathfrak P$ so that $Q\in S_1$ and $\mu_{p}Q\in S_2$. Then in $\mathbb P$, we have
$$\frac{\bigoplus_{P\in S_1}y^{T}(P)}{\bigoplus_{P'\in \pi(S_1)}y^{T'}(P')}=\frac{\bigoplus_{P\in S_2}y^{T}(P)}{\bigoplus_{P'\in \pi(S_2)}y^{T'}(P')}.$$

\end{Proposition}

\begin{proof}

It follows immediately by Lemma \ref{non-a-F}, Lemma \ref{tau-F} and Lemma \ref{all-a-F}.
\end{proof}

\medskip

As a corollary of Proposition \ref{all-F}, the following theorem follows.

\medskip

\begin{Theorem}\label{F}

With the same notation as above, we assume that $\tau'\neq \gamma\notin T$.

\begin{enumerate}[$(1)$]

  \item For any $S\in \mathfrak P$, we have
  $$\frac{\bigoplus_{P\in S}y^{T}(P)}{\bigoplus_{P'\in \pi(S)}y^{T'}(P')}=\frac{\bigoplus_{P\in \mathcal P(G_{T,\gamma})}y^{T}(P)}{\bigoplus_{P'\in \mathcal P(G_{T',\gamma})}y^{T'}(P')}.$$

  \item For any $S'\in \mathfrak P'$, we have
  $$\frac{\bigoplus_{P\in \pi'(S')}y^{T}(P)}{\bigoplus_{P'\in S'}y^{T'}(P')}=\frac{\bigoplus_{P\in \mathcal P(G_{T,\gamma})}y^{T}(P)}{\bigoplus_{P'\in \mathcal P(G_{T',\gamma})}y^{T'}(P')}.$$

\end{enumerate}

\end{Theorem}

\begin{proof}

We shall only prove (1) because (2) can be proved dually. By Lemma \ref{transitive} and Proposition \ref{all-F}, $\frac{\bigoplus_{P\in S_1}y^{T}(P)}{\bigoplus_{P'\in \pi(S_1)}y^{T'}(P')}=\frac{\bigoplus_{P\in S_2}y^{T}(P)}{\bigoplus_{P'\in \pi(S_2)}y^{T'}(P')}$ for any $S_1,S_2\in \mathfrak P$. Therefore, the result follows.
\end{proof}

\medskip

{\bf Proof of Theorem \ref{partition bi} (2):}
By Propositions \ref{compare-loc} (2.b), \ref{2^n}, Lemma \ref{brick1} and Theorem \ref{F}, the proof is the same as that of Theorem 4.6 (2) in \cite{H1}. \;\ \ \ \  $\square$

\medskip

\section{Proof of Theorem \ref{mainthm}}\label{main2}

\subsection{Valuation maps on \textbf{$\mathcal P(G_{T,\gamma})$}}

Herein, the valuation maps $v_{+},v_{-}: \mathcal P(G_{T,\gamma})\rightarrow \mathbb Z$ constructed in \cite{H,H1} will be generalized to the orbifold case. We will prove $v_{+}=v_{-}$, which is the required $v$ in Theorem \ref{mainthm}. Throughout this section let $\mathcal O$ be an unpunctured orbifold and $T$ be a triangulation. Let $\gamma$ be an oriented arc in $\mathcal O$ and $\zeta$ be the corresponding ordinary arc which crosses $T$ with points $p_1,\cdots, p_d$ in order. We assume that $p_1,\cdots,p_d$ belong to the arcs $\tau_{i_1},\cdots,\tau_{i_d}$, respectively in $T$.

\medskip

\begin{Lemma}\label{invo}\cite[Lemma 6.1]{H1}
If $P\in \mathcal P(G_{T,\gamma})$ can twist on $G(p_s)$, then $\mu_{p_s}P$ can twist on $G(p_s)$, and $\mu_{p_s}\mu_{p_s}P=P$.

\end{Lemma}

\medskip

\begin{Lemma}\label{commuta}

If $P\in \mathcal P(G_{T,\gamma})$ can twist on $G(p_s)$ and $G(p_t)$ for some $s,t$ so that $|s-t|>1$, then
\begin{enumerate}[$(1)$]

  \item $\mu_{p_s}P$ and $\mu_{p_t}P$ can twist on $G(p_t)$ and $G(p_s)$, respectively, and $$\mu_{p_s}\mu_{p_t}P=\mu_{p_t}\mu_{p_s}P.$$

  \item $\Omega(p_s,P)+\Omega(p_t,\mu_{p_s}P)=\Omega(p_t,P)+\Omega(p_s,\mu_{p_t}P)$.

\end{enumerate}

\end{Lemma}

\begin{proof}

(1)~ This is Lemma 6.2 (1) in \cite{H}.

(2)~ As $\Omega(p,Q)=-\Omega(p,\mu_{p}Q)$ for any $Q$ which can twist on $G(p)$, we may assume that the edges of $G(p_s)$ and $G(p_t)$ which are labeled $a_{2_s},a_{4_s}$ and $a_{2_t},a_{4_t}$, respectively, in Definition \ref{omega} belong to $P$.

In case $\tau_{i_s}$ and $\tau_{i_t}$ are not in a same triangle in $T$ or $\tau_{i_s}=\tau_{i_t}$, we have $n^{\pm}_{p_s}(\tau_{i_s},P)=n^{\pm}_{p_s}(\tau_{i_s},\mu_{p_t}P)$, and $n^{\pm}_{p_t}(\tau_{i_t},P)=n^{\pm}_{p_t}(\tau_{i_t},\mu_{p_s}P)$. Thus,
$$
\begin{array}{rcl} && \Omega(p_s,P)+\Omega(p_t,\mu_{p_s}P) \vspace{2pt}  \\

& = & [n^{+}_{p_s}(\tau_{i_s},P)-m^{+}_{p_s}(\tau_{i_s},\gamma)-n^{-}_{p_s}(\tau_{i_s},P)+m^{-}_{p_s}(\tau_{i_s},\gamma)]d^T(\tau_{i_s}) \vspace{2pt}  \\

& + & [n^{+}_{p_t}(\tau_{i_t},\mu_{p_s}P)-m^{+}_{p_t}(\tau_{i_t},\gamma)-n^{-}_{p_t}(\tau_{i_t},\mu_{p_s}P)+m^{-}_{p_t}(\tau_{i_t},\gamma)]d^T(\tau_{i_t}) \vspace{2pt}  \\

& = & [n^{+}_{p_s}(\tau_{i_s},\mu_{p_t}P)-m^{+}_{p_s}(\tau_{i_s},\gamma)-n^{-}_{p_s}(\tau_{i_s},\mu_{p_t}P)+m^{-}_{p_s}(\tau_{i_s},\gamma)]d^T(\tau_{i_s}) \vspace{2pt}  \\

& + & [n^{+}_{p_t}(\tau_{i_t},P)-m^{+}_{p_t}(\tau_{i_t},\gamma)-n^{-}_{p_t}(\tau_{i_t},P)+m^{-}_{p_t}(\tau_{i_t},\gamma)]d^T(\tau_{i_t}) \vspace{2pt}  \\

& = & \Omega(p_s,\mu_{p_t}P)+\Omega(p_t,P).\vspace{2pt}\\
\end{array}$$

As $\Omega(p,Q)=-\Omega(p,\mu_{p}Q)$ for any $Q$ which can twist on $G(p)$, the following equalities are equivalence
$$\Omega(p_s,P)+\Omega(p_t,\mu_{p_s}P)=\Omega(p_t,P)+\Omega(p_s,\mu_{p_t}P),$$ $$\Omega(p_s,\mu_{p_s}P)+\Omega(p_t,P)=\Omega(p_t,\mu_{p_s}P)+\Omega(p_s,\mu_{p_s}\mu_{p_t}P),$$
$$\Omega(p_s,\mu_{p_t}P)+\Omega(p_t,\mu_{p_s}\mu_{p_t}P)=\Omega(p_t,\mu_{p_t}P)+\Omega(p_s,P).$$

In case $\tau_{i_s}$ and $\tau_{i_t}$ are in a same triangle in $T$, we therefore shall assume that the edges labeled $\tau_{i_t}$ and $\tau_{i_s}$ of $G(p_s)$ and $G(p_t)$, respectively are in $P$, because otherwise $P$ can twists on some tiles to make the assumption holds. Without loss of generality, assume that $s<t$ and $\tau_{i_t}$ is counterclockwise to $\tau_{i_s}$ in the triangle of $T$. Then we have $n_{p_s}^{+}(\tau_{i_s},\mu_{p_t}P)=n_{p_s}^{+}(\tau_{i_s},P)-|b^T_{i_si_t}|$, $n_{p_s}^{-}(\tau_{i_s},\mu_{p_t}P)=n_{p_s}^{-}(\tau_{i_s},P)$ and $n_{p_t}^{+}(\tau_{i_t},\mu_{p_s}P)=n_{p_t}^{+}(\tau_{i_t},P)$, $n_{p_t}^{-}(\tau_{i_t},\mu_{p_s}P)=n_{p_t}^{-}(\tau_{i_t},P)-|b^T_{i_ti_s}|$. As $B^T$ is skew-symmetrizable, $d(\tau_{i_s})b^T_{i_si_t}=-d(\tau_{i_t})b^T_{i_ti_s}$.
Thus, $$
\begin{array}{rcl} && \Omega(p_s,P)+\Omega(p_t,\mu_{p_s}P) \vspace{2pt}  \\

& = & [n^{+}_{p_s}(\tau_{i_s},P)-m^{+}_{p_s}(\tau_{i_s},\gamma)-n^{-}_{p_s}(\tau_{i_s},P)+m^{-}_{p_s}(\tau_{i_s},\gamma)]d^T(\tau_{i_s}) \vspace{2pt}  \\

& - & [n^{+}_{p_t}(\tau_{i_t},\mu_{p_s}P)-m^{+}_{p_t}(\tau_{i_t},\gamma)-n^{-}_{p_t}(\tau_{i_t},\mu_{p_s}P)+m^{-}_{p_t}(\tau_{i_t},\gamma)]d^T(\tau_{i_t}) \vspace{2pt}  \\

& = & [n^{+}_{p_s}(\tau_{i_s},\mu_{p_t}P)+|b^T_{i_si_t}|-m^{+}_{p_s}(\tau_{i_s},\gamma)-n^{-}_{p_s}(\tau_{i_s},\mu_{p_t}P)+m^{-}_{p_s}(\tau_{i_s},\gamma)]d^T(\tau_{i_s}) \vspace{2pt}  \\

& - & [n^{+}_{p_t}(\tau_{i_t},P)-m^{+}_{p_t}(\tau_{i_t},\gamma)-n^{-}_{p_t}(\tau_{i_t},P)+|b^T_{i_ti_s}|+m^{-}_{p_t}(\tau_{i_t},\gamma)]d^T(\tau_{i_t}) \vspace{2pt}  \\

& = & \Omega(p_s,\mu_{p_t}P)+\Omega(p_t,P).\vspace{2pt}\\
\end{array}$$
The proof is complete.
\end{proof}

\medskip

\begin{Lemma}\label{reduction}\cite[Lemma 6.4]{H} Let $P$ be a perfect matching of $G_{T,\gamma}$. Suppose that $(p_{i_1},\cdots,p_{i_r})$ is a sequence so that $\mu_{p_{i_{t-1}}}\cdots\mu_{p_{i_1}}$ can twist on $G(p_{i_t})$ for $1\leq t\leq r$ and $\mu_{p_{i_r}}\cdots\mu_{p_{i_1}}P=P$. Then

\begin{enumerate}[$(1)$]

  \item there exists $2\leq t\leq r$ such that $p_{i_t}=p_{i_1}$.

  \item There exist $t<t'$ such that $p_{i_{t'}}=p_{i_{t}}$ and $p_{i_{s'}}\neq p_{i_s}$ for any $s'\neq s$ with $t<s,s'<t'$; in this case, $|i_s-i_t|> 1$ for any $s$ satisfies $t<s<t'$.

\end{enumerate}

\end{Lemma}

\medskip

\begin{Theorem}\label{valumap}

Let $(\mathcal O,M)$ be an unpunctured orbifold, $\gamma$ be an oriented arc in $\mathcal O$ and $T$ be an indexed triangulation of $\mathcal O$. Then

\begin{enumerate}[$(1)$]

\item there uniquely exists a \emph{maximal valuation map} $v_{+}:\mathcal P(G_{T,\gamma})\rightarrow \mathbb Z$ satisfying

\begin{enumerate}[$(a)$]

  \item (initial condition) $v_{+}(P_{+}(G_{T,\gamma}))=0$.

  \item (iterated relation) If $P\in \mathcal P(G_{T,\gamma}))$ can twist on $G(p)$, then $$v_{+}(P)-v_{+}(\mu_pP)=\Omega(p,P).$$

\end{enumerate}

\item There uniquely exists a \emph{minimal valuation map} $v_{-}:\mathcal P(G_{T,\gamma})\rightarrow \mathbb Z$ satisfying

\begin{enumerate}[$(a)$]

  \item (initial condition) $v_{-}(P_{-}(G_{T,\gamma}))=0$.

  \item (iterated relation) If $P\in \mathcal P(G_{T,\gamma}))$ can twist on $G(p)$, then $$v_{-}(P)-v_{-}(\mu_pP)=\Omega(p,P).$$

\end{enumerate}

\end{enumerate}

\end{Theorem}

\begin{proof}

By Lemma \ref{commuta} and Lemma \ref{reduction}, the proof is similar to that of \cite[Theorem 6.5]{H} and \cite[Theorem 8.4]{H1}.
\end{proof}

\medskip

\subsection{$v_{+}=v_{-}$}

Herein, we prove $v_+=v_-$ and Theorem \ref{mainthm}. By definition, $X^{T'}_{\tau'}=(X^T)^{-e_{\tau}+(b^{T}_{\tau})_{+}}+(X^T)^{-e_{\tau}+(b^{T}_{\tau})_{-}}$. We denote $(X^T)^{-e_\tau+(b^{T}_{\tau})_{\pm}}$ by $\prod_{\pm}$ for convenience. Let $d$ be a non-negative integer. For any sequence $\lambda=(\lambda_1,\cdots,\lambda_d)\in \{0,1\}^d$, denote by $n(\lambda)$ the integer so that
$$Z_1Z_2\cdots Z_d=q^{n(\lambda)/2}(X^T)^{-de_{\tau}+(\sum\lambda_i)(b^{T}_{\tau})_{-}+(d-\sum\lambda_i)(b^{T}_{\tau})_{+}},$$
where
  \[\begin{array}{ccl} Z_i &=&

         \left\{\begin{array}{ll}

             \prod_{-}, &\mbox{if $\lambda_i=1$}, \\

             \prod_{+}, &\mbox{if $\lambda_i=0$}.

         \end{array}\right.

 \end{array}\]
Clearly, $n(\lambda)=0$ if all $\lambda_i=1$ or all $\lambda_i=0$. Using the notation of $n(\lambda)$, we have
$$(X^{T'}_{\tau'})^d=\textstyle\sum_{\lambda}q^{n(\lambda)/2}(X^T)^{-de_{\tau}+(\sum\lambda_i)(b^{T}_{\tau})_{-}+(d-\sum\lambda_i)(b^{T}_{\tau})_{+}}.$$

\medskip

The following lemma generalizes \cite[Lemma 8.5]{H1} from the surface to orbifold.

\medskip

\begin{Lemma}\label{iterate}

With the same notation as above, for any $i\in [1,d]$, if $\lambda$ and $\lambda'$ satisfy $\lambda_j=\lambda'_j$ for $j\neq i$ and $\lambda_i=1$, $\lambda'_i=0$, then $$n(\lambda)-n(\lambda')=(d-2i+1)d^T(\tau),$$
where $d^T(\tau)$ is the integer determined by the compatibility of $(\widetilde B^T,\Lambda^T)$.

\end{Lemma}

\begin{proof}

We assume $Z_1Z_2\cdots Z_{i-1}=q^{d_1/2}(X^T)^{\vec{a}}$ and $Z_{i+1}\cdots Z_{d}=q^{d_2/2}(X^T)^{\vec{b}}$ for some $d_1,d_2\in \mathbb Z$ and $\vec{a},\vec{b}\in \mathbb Z^m$, where
  \[\begin{array}{ccl} Z_j &=&

         \left\{\begin{array}{ll}

             \prod_{-}, &\mbox{if $\lambda_j=1$}, \\

             \prod_{+}, &\mbox{if $\lambda_j=0$}.

         \end{array}\right.

 \end{array}\]
Thus, we have $\vec{a}=-(i-1)e_{\tau}+(\sum_{j<i}\lambda_j)(b^{T}_{\tau})_{-}+(i-1-\sum_{j<i}\lambda_j)(b^{T}_{\tau})_{+},$ and
$\vec{b}=-(d-i)e_{\tau}+(\sum_{j>i}\lambda_j)(b^{T}_{\tau})_{-}+(d-i-\sum_{j>i}\lambda_j)(b^{T}_{\tau})_{+}$.

According to the definitions of $n(\lambda)$ and $n(\lambda')$, we obtain
$$q^{d_1/2}(X^T)^{\vec{a}}\cdot \textstyle\prod_{-}\cdot q^{d_2/2}(X^T)^{\vec{b}}=q^{n(\lambda)/2}(X^T)^{-de_{\tau}+(\sum\lambda_i)(b^{T}_{\tau})_{-}+(d-\sum\lambda_i)(b^{T}_{\tau})_{+}},$$
$$q^{d_1/2}(X^T)^{\vec{a}}\cdot\textstyle\prod_{+}\cdot q^{d_2/2}(X^T)^{\vec{b}}=q^{n(\lambda')/2}(X^T)^{-de_{\tau}+(\sum\lambda'_i)(b^{T}_{\tau})_{-}+(d-\sum\lambda'_i)(b^{T}_{\tau})_{+}}.$$

Thus,
$$n(\lambda)=d_1+d_2+\Lambda^T(\vec{a}, \vec{b})+\Lambda^T(\vec{a}, (b^{T}_{\tau})_{-})+\Lambda^T((b^{T}_{\tau})_{-},\vec{b}),$$
and
$$n(\lambda')=d_1+d_2+\Lambda^T(\vec{a}, \vec{b})+\Lambda^{T}(\vec{a},(b^{T}_{\tau})_{+})+\Lambda^T((b^{T}_{\tau})_{+},\vec{b}).$$
Consequently,
$$
\begin{array}{rcl} n(\lambda)-n(\lambda') & = & \Lambda^T(\vec{a}-\vec{b}, (b^{T}_{\tau})_{-}-(b^{T}_{\tau})_{+}) \vspace{2pt} \\

& = &  \Lambda^T(\vec{a}-\vec{b}, -b^{T}_{\tau})\vspace{2pt} \\

& = &  \Lambda^T(b^{T}_{\tau}, \vec{a}-\vec{b})\vspace{2pt} \\

& = &  [d-i-(i-1)]d^T(\tau) = (d-2i+1)d^T(\tau),
\end{array}$$
where the last second equality follows by the compatibility of $(\widetilde B^T,\Lambda^T)$.
\end{proof}

\medskip

We denote by $v'_{\pm}$ the maximal and minimal valuation maps on $\mathcal P(G_{T',\gamma})$, which exist by Theorem \ref{valumap}. For $P'\in \mathcal P^{\tau'}_{\nu}(G_{T',\gamma})$ such that $d=\sum \nu_l\geq 0$, let $P(\lambda)\in \pi'(P')$ corresponding to $\lambda\in \{0,1\}^d$ as in Proposition \ref{2^n}.

\medskip

\begin{Lemma}\cite[Lemma 6.9]{H}\label{induction}
With the same notation as above, we assume $P'\in \mathcal P^{\tau'}_{\nu}(G_{T',\gamma})$ so that $\sum \nu_l\geq 0$. Then

\begin{enumerate}[$(1)$]

  \item the following are equivalent.

  \begin{enumerate}[$(a)$]

  \item $q^{v'_{+}(P')/2}X^{T'}(P')=\sum_{P\in \pi'(P')}q^{v_{+}(P)/2}X^{T}(P)$;

  \item $v'_{+}(P')=v_{+}(P(0,0,\cdots,0))$;

  \item $v'_{+}(P')=v_{+}(P(1,1,\cdots,1))$.

  \end{enumerate}

  \item The following are equivalent.

  \begin{enumerate}[$(a)$]

  \item $q^{v'_{-}(P')/2}X^{T'}(P')=\sum_{P\in \pi'(P')}q^{v_{-}(P)/2}X^{T}(P)$;

  \item $v'_{-}(P')=v_{-}(P(0,0,\cdots,0))$;

  \item $v'_{-}(P')=v_{-}(P(1,1,\cdots,1))$.

\end{enumerate}

\end{enumerate}

\end{Lemma}

\begin{proof}

Using Lemma \ref{iterate}, the proof is the same as that of Lemma 6.9 in \cite{H}.
\end{proof}

\begin{Lemma}\label{minmax}

With the same notation as above, let $P=P_{\pm}(G_{T,\gamma})\in \mathcal P^{\tau}_{\nu}(G_{T,\gamma})$ for some $\nu$. If $\nu_i<0$ for some $i$, then $\nu_j\leq 0$ for any $j$.

\end{Lemma}

\begin{proof}

Without loss of generality, we may assume $P=P_{+}(G_{T,\gamma})$. We assume $\nu_j>0$ for some $j$. Then the $j$-th $\tau$-equivalence class is of type (I,IV). Since $\nu_i<0$ and the edges in $P$ are boundary edges, the $i$-th $\tau$-equivalence class is of type (I,III).

We denote the endpoints of the two triangles containing $\tau$ by $o_1,o_2,o_3,o_4$. As shown in the following graph. It should be noted that $o_1=o_3,o_2=o_4$ if $\tau$ is a pending arc.

\centerline{\begin{tikzpicture}
\draw[-] (1,0) -- (2,0);
\draw[-] (1,-1) -- (2,-1);
\draw[-] (1,0) -- (2,-1);
\node[right] at (1.4,-0.45) {$\tau$};
\node[above] at (1.5,-0.05) {$a_2$};
\node[left] at (1,-0.5) {$a_1$};
\node[left] at (1,0) {$o_1$};
\node[right] at (2,0) {$o_2$};
\node[left] at (1,-1) {$o_3$};
\node[right] at (2,-1) {$o_4$};
\node[below] at (1.5,-0.95) {$a_4$};
\node[right] at (2,-0.5) {$a_3$};
\draw [-] (1, 0) -- (1,-1);
\draw [-] (2, 0) -- (2,-1);
\draw [fill] (1,0) circle [radius=.05];
\draw [fill] (2,0) circle [radius=.05];
\draw [fill] (1,-1) circle [radius=.05];
\draw [fill] (2,-1) circle [radius=.05];
\end{tikzpicture}}

In case the $j$-th $\tau$-equivalence class is of type (I), then $\zeta$ crosses $a_2,\tau,a_4$ consequently by Lemma \ref{max-min}. If the $i$-th $\tau$-equivalence class is of type (I), then $\zeta$ crosses $a_1,\tau,a_3$ consequently. If $i$-th $\tau$-equivalence class is of type (III), then $\zeta$ starts from $o_3$ and crosses $\tau,a_3$ or starts from $o_2$ and $\tau,a_1$ consequently. In both cases, $\zeta$ crosses itself.

In case the $j$-th $\tau$-equivalence class is of type (IV), then the $j$-th $\tau$ equivalence class contains an edge of a tile $G$ with diagonal labeled $a_2$ or $a_4$ by Lemma \ref{max-min}. We may assume that the diagonal is labeled $a_2$. Moreover, since the edge of $G$ labeled $\tau$ is a boundary edge, $\zeta$ starts from $o_4$ and crosses $a_2$. If the $i$-th $\tau$-equivalence classes is of type (I), then $\zeta$ crosses $a_1,\tau,a_3$ consequently. If $i$-th $\tau$-equivalence class is of type (III), then $\zeta$ starts from $o_3$ and crosses $\tau,a_3$ or starts from $o_2$ and $\tau,a_1$ consequently. In both cases, $\zeta$ crosses itself.

Therefore, the result follows.
\end{proof}

\medskip

\begin{Lemma}\label{maxtomax1}

With the same notation as above, if $\tau'\neq \gamma\notin T$, then

\begin{enumerate}[$(1)$]

  \item $P_{\pm}(G_{T',\gamma})\subset \pi(P_{\pm}(G_{T,\gamma}))$.

  \item $P_{\pm}(G_{T,\gamma})\subset \pi'(P_{\pm}(G_{T',\gamma}))$

\end{enumerate}

\end{Lemma}

\begin{proof}

Using Lemmas \ref{max-min}, \ref{maxtomax}, \ref{formglobal} and \ref{minmax}, the proof is similar to that of \cite[Lemma 8.7]{H1}.
\end{proof}

\medskip

Let $S_{\pm}\in \mathfrak P$ with $P_{\pm}(G_{T,\gamma})\in S_{\pm}$ and $S'_{\pm}\in \mathfrak P'$ with $P_{\pm}(G_{T',\gamma})\in S'_{\pm}$. Thus, $\pi'(S'_{+})=S_{+}$ and $\pi'(S'_{-})=S_{-}$ by Lemma \ref{maxtomax1}.

\medskip

As applications of Lemma \ref{induction}, we have the following statements.

\medskip

\begin{Proposition}\label{initial}

With the same notation as above, we assume $\tau'\neq\gamma\notin T$. Then

\begin{enumerate}[$(1)$]

  \item $\sum_{P'\in S'_{+}}q^{v'_{+}(P')/2}X^{T'}(P')=\sum_{P\in S_{+}}q^{v_{+}(P)/2}X^T(P)$.

  \item $\sum_{P'\in S'_{-}}q^{v'_{-}(P')/2}X^{T'}(P')=\sum_{P\in S_{-}}q^{v_{-}(P)/2}X^T(P)$.

\end{enumerate}

\end{Proposition}

\begin{proof}

The proof is the same as that of \cite[Proposition 6.10]{H}.
\end{proof}

\medskip

\begin{Proposition}\label{non-a}

With the same notation as above, we assume $\tau'\neq \gamma\notin T$. Let $Q$ be a perfect matching of $G_{T,\gamma}$ which can twist on a tile $G(p)$ with the diagonal not labeled $a_1,a_2,a_3,a_4$. Let $S_1,S_2\in \mathfrak P$ so that $Q\in S_1$ and $\mu_{p}Q\in S_2$. Then

\begin{enumerate}[$(1)$]

  \item $\sum_{P'\in \pi(S_1)}q^{v'_{+}(P')/2}X^{T'}(P')=\sum_{P\in S_1}q^{v_{+}(P)/2}X^{T}(P)$ holds if and only if\\
   $\sum_{P'\in \pi(S_2)}q^{v'_{+}(P')/2}X^{T'}(P')=\sum_{P\in S_2}q^{v_{+}(P)/2}X^{T}(P)$ holds.

  \item $\sum_{P'\in \pi(S_1)}q^{v'_{-}(P')/2}X^{T'}(P')=\sum_{P\in S_1}q^{v_{-}(P)/2}X^{T}(P)$ holds if and only if\\
   $\sum_{P'\in \pi(S_2)}q^{v'_{-}(P')/2}X^{T'}(P')=\sum_{P\in S_2}q^{v_{-}(P)/2}X^{T}(P)$ holds.

\end{enumerate}

\end{Proposition}

\begin{proof}

The proof is similar to that of \cite[Lemma 8.10]{H1}.
\end{proof}

\medskip

\begin{Lemma}\label{special-a}

With the same notation as above, we assume $\tau'\neq\gamma\notin T$. Let $Q$ be a perfect matching of $G_{T,\gamma}$ which can twist on a tile $G(p)$ with diagonal labeled $a=a_q$ for $q=1,2,3,4$. We assume $S_1,S_2\in \mathfrak P$ so that $Q\in S_1$ and $\mu_{p}Q\in S_2$. If all the $\tau$-mutable edges pairs in $Q$ are labeled $a_{q-1},a_{q+1}$ (addition in $\mathbb Z_4$), then

\begin{enumerate}[$(1)$]

  \item $\sum_{P'\in \pi(S_1)}q^{v'_{+}(P')/2}X^{T'}(P')=\sum_{P\in S_1}q^{v_{+}(P)/2}X^{T}(P)$ holds if and only if\\
   $\sum_{P'\in \pi(S_2)}q^{v'_{+}(P')/2}X^{T'}(P')=\sum_{P\in S_2}q^{v_{+}(P)/2}X^{T}(P)$ holds.

  \item $\sum_{P'\in \pi(S_1)}q^{v'_{-}(P')/2}X^{T'}(P')=\sum_{P\in S_1}q^{v_{-}(P)/2}X^{T}(P)$ holds if and only if\\
   $\sum_{P'\in \pi(S_2)}q^{v'_{-}(P')/2}X^{T'}(P')=\sum_{P\in S_2}q^{v_{-}(P)/2}X^{T}(P)$ holds.

\end{enumerate}

\end{Lemma}

\begin{proof}

The proof is similar to that of \cite[Lemma 8.11]{H1}.
\end{proof}

\medskip

\begin{Proposition}\label{general-a}

With the same notation as above, we assume $\tau'\neq\gamma\notin T$. Let $Q$ be a perfect matching of $G_{T,\gamma}$ which can twist on a tile $G(p)$ with diagonal labeled $a_q$ for $q=1,2,3,4$. Let $S_1,S_2\in \mathfrak P$ such that $Q\in S_1$ and $\mu_{p}Q\in S_2$. Then

\begin{enumerate}[$(1)$]

  \item $\sum_{P'\in \pi(S_1)}q^{v'_{+}(P')/2}X^{T'}(P')=\sum_{P\in S_1}q^{v_{+}(P)/2}X^{T}(P)$ holds if and only if\\
   $\sum_{P'\in \pi(S_2)}q^{v'_{+}(P')/2}X^{T'}(P')=\sum_{P\in S_2}q^{v_{+}(P)/2}X^{T}(P)$ holds.

  \item $\sum_{P'\in \pi(S_1)}q^{v'_{-}(P')/2}X^{T'}(P')=\sum_{P\in S_1}q^{v_{-}(P)/2}X^{T}(P)$ holds if and only if\\
   $\sum_{P'\in \pi(S_2)}q^{v'_{-}(P')/2}X^{T'}(P')=\sum_{P\in S_2}q^{v_{-}(P)/2}X^{T}(P)$ holds.

\end{enumerate}

\end{Proposition}

\begin{proof}

The result follows by Proposition \ref{non-a} and Lemma \ref{special-a}.
\end{proof}

\medskip

Summaries Propositions \ref{non-a} and \ref{general-a}, we obtain the following.

\medskip

\begin{Proposition}\label{all}

With the same notation as above, we assume $\tau'\neq \gamma\notin T$. Let $Q$ be a perfect matching of $G_{T,\gamma}$ which can twist on a tile $G(p)$. Let $S_1,S_2\in \mathfrak P$ so that $Q\in S_1$ and $\mu_{p}Q\in S_2$. Then

\begin{enumerate}[$(1)$]

  \item $\sum_{P'\in \pi(S_1)}q^{v'_{+}(P')/2}X^{T'}(P')=\sum_{P\in S_1}q^{v_{+}(P)/2}X^{T}(P)$ holds if and only if\\
   $\sum_{P'\in \pi(S_2)}q^{v'_{+}(P')/2}X^{T'}(P')=\sum_{P\in S_2}q^{v_{+}(P)/2}X^{T}(P)$ holds.

  \item $\sum_{P'\in \pi(S_1)}q^{v'_{-}(P')/2}X^{T'}(P')=\sum_{P\in S_1}q^{v_{-}(P)/2}X^{T}(P)$ holds if and only if\\
   $\sum_{P'\in \pi(S_2)}q^{v'_{-}(P')/2}X^{T'}(P')=\sum_{P\in S_2}q^{v_{-}(P)/2}X^{T}(P)$ holds.

\end{enumerate}

\end{Proposition}

\medskip

\begin{Theorem}\label{mainpre}

With the same notation as above, for any $S\in \mathfrak P$,

\begin{enumerate}[$(1)$]

  \item $\sum_{P\in S}q^{v_{+}(P)/2}X^{T}(P)=\textstyle\sum_{P'\in \pi(S)}q^{v'_{+}(P')/2}X^{T'}(P')$.

  \item $\sum_{P\in S}q^{v_{-}(P)/2}X^{T}(P)=\textstyle\sum_{P'\in \pi(S)}q^{v'_{-}(P')/2}X^{T'}(P')$.

\end{enumerate}

\end{Theorem}

\begin{proof}

We shall only prove (1) because (2) can be proved similarly. If $\gamma\neq \tau,\tau'$, for any $P\in S$, by Lemma \ref{transitive}, $P$ can be obtained from $P_{+}(G_{T,\gamma})$ by a sequence of twists. In this case, the equality follows from Proposition \ref{initial} and Proposition \ref{all}. If $\gamma=\tau$ or $\tau'$, then the equality becomes to
$$X^{T}_{\tau}=(X^{T'})^{-e_{\tau'}+(b^{T'}_{\tau'})_{+}}+(X^{T'})^{-e_{\tau'}+(b^{T'}_{\tau'})_{-}}\;\;\text{or}\;\; X^{T'}_{\tau'}=(X^{T})^{-e_{\tau}+(b^T_{\tau})_{+}}+(X^T)^{-e_{\tau}+(b^T_{\tau})_{-}},$$
which clearly holds. The result follows.
\end{proof}

\medskip

We now give the proof of Theorem \ref{mainthm}.

\medskip

{\bf Proof of Theorem \ref{mainthm}:} By Theorem \ref{mainpre}, it suffices to show $v_{+}=v_{-}$. We prove this by induction on $N(\gamma,T)$, the number of crossing points of $\gamma$ with $T$. When $N(\gamma,T)=0$ or $1$, it is easy to see that $v_{+}=v_{-}=0$. We assume $v_{+}=v_{-}$ for all $N(\gamma,T)<d$. When $N(\gamma,T)=d>1$, there exists $\tau\in T$ such that the $N(\gamma,T')<d$. We denote by $v'_{\pm}$ the maximal and minimal valuation maps on $\mathcal P(G_{T',\gamma})$. By induction hypothesis, $v'_{+}=v'_{-}$. By Lemma \ref{maxtomax1}, $P_{+}(G_{T,\gamma})\in \pi'(P_{+}(G_{T',\gamma}))$. We assume $P_{+}(G_{T,\gamma})\in \mathcal P^{\tau}_{\nu}(G_{T,\gamma})$ for some $\nu$.

If $\sum\nu_l\leq 0$, by the dual version of Proposition \ref{2^n} and Lemma \ref{max-min}, $P_{+}(G_{T,\gamma})=P(\lambda)$ with $\lambda=(0,0,\cdots,0)$ or $(1,1,\cdots,1)\in \{0,1\}^{-\sum\nu_l}$, where $P(\lambda)$ is the perfect matching in $\pi'(P_{+}(G_{T',\gamma}))$ which is determined by $\lambda$. By Theorem \ref{mainpre} and Lemma \ref{induction} (2), we have $v_{-}(P_{+}(G_{T,\gamma}))=v'_{-}(P_{+}(G_{T',\gamma})).$ Since $v'_{+}=v'_{-}$, $v'_{-}(P_{+}(G_{T',\gamma}))=0$, and hence $v_{-}(P_{+}(G_{T,\gamma}))=0$. Consequently, $v_{+}=v_{-}$.

If $\sum\nu_l\geq 0$, by Proposition \ref{2^n} and Lemma \ref{max-min}, $P_{+}(G_{T',\gamma})=P'(\lambda)$ with $\lambda=(0,0,\cdots,0)$ or $(1,1,\cdots,1)\in \{0,1\}^{\sum\nu_l}$, where $P'(\lambda)$ is the perfect matching in $\pi(P_{+}(G_{T,\gamma}))$ which is determined by $\lambda$. By Theorem \ref{mainpre} and the dual version of Lemma \ref{induction} (2), we have $v_{-}(P_{+}(G_{T,\gamma}))=v'_{-}(P_{+}(G_{T',\gamma})).$ Since $v'_{+}=v'_{-}$, $v'_{-}(P_{+}(G_{T',\gamma}))=0$, and hence $v_{-}(P_{+}(G_{T,\gamma}))=0$. Consequently, $v_{+}=v_{-}$.

The proof of Theorem \ref{mainthm} is complete.

\medskip

{\bf Acknowledgements:}\;  {\em The author is thankful to S. Liu, I. Assem, T. Br$\ddot{u}$stle and D. Smith for financial support.}

\end{document}